\documentclass{article}

\usepackage{graphicx}
\usepackage{xcolor}
\usepackage{amsfonts}
\usepackage{setspace}
\usepackage{amsmath}
\usepackage{amssymb}
\usepackage{comment}
\usepackage{apptools}
\usepackage[round]{natbib}
\usepackage{float}
\usepackage{subfigure}
\usepackage{hyperref}
\hypersetup{
	colorlinks=true,
	linkcolor=blue,
	citecolor=blue,
	filecolor=blue,
	urlcolor=blue,
}

\definecolor{verde}{rgb}{0.2, 0.9, 0.6}

\newcommand{\bl}{\color{blue}}

\DeclareRobustCommand{\FIN}{%
  \ifmmode 
  \else \leavevmode\unskip\penalty9999 \hbox{}\nobreak\hfill
  \fi
  $\bullet$ \vspace{5mm}}

\newcommand{\calC}{\mbox{${\cal C}$}}
\newcommand{\calG}{\mbox{${\cal G}$}}

\newcommand{\cvec}{\mbox{\bf c}}
\newcommand{\zvec}{\mbox{\bf z}}

\newcommand{\xvec}{\mbox{\bf x}}

\newcommand{\ox}{\overline{x}}
\newcommand{\oy}{\overline{y}}
\newcommand{\uvec}{\mbox{\bf u}}
\newcommand{\vvec}{\mbox{\bf v}}
\newcommand{\Tvec}{\mbox{\bf T}}

\newcommand{\Xvec}{\mbox{\bf X}}
\newcommand{\Uvec}{\mbox{\bf U}}

\newcommand{\Zvec}{\mbox{\bf Z}}

\newcommand{\disp}{\displaystyle}

\newcommand{\gammavec}{\mbox{\boldmath $\gamma$}}

\newcommand{\Nat}{\mbox{$\mathbb{N}$}}
\newcommand{\Rea}{\mbox{$\mathbb{R}$}}
\newcommand{\Hil}{\mbox{$\mathbb{H}$}}
\newcommand{\Prob}{\mbox{$\mathbb{P}$}}

\newcommand{\ProbZ}{\mbox{$\Prob_{\bf Z}$}}
\newcommand{\SigmaZ}{\mbox{$\Sigma^{\bf Z}$}}

\newcommand{\Esp}{\mbox{$\mathbb{E}$}}

\newcommand{\nin}{\mbox{$n \in \mathbb{N}$}}

\newcommand{\conv}{\rightarrow}

\newcommand{\convp}{\stackrel{P}{\to}}

\newtheorem {Prop}{Proposition} [section]
\newtheorem {Lemm}[Prop] {Lemma}

\newtheorem {Theo}[Prop]{Theorem}

\newtheorem {Nota}{Remark}[Prop]

\newtheorem {Defi}[Prop]{Definition}

\numberwithin{table}{section}

\begin{document}

\thispagestyle{empty} \baselineskip=28pt \vskip 5mm
\begin{center} 
{\Huge{\sc On Perfect Classification and Clustering for Gaussian Processes}}
\end{center}

\baselineskip=12pt \vskip 5mm

\begin{center}\large
{Juan A. Cuesta-Albertos
\footnote{This author has been partially supported by the
Spanish Ministerio de Ciencia y Tecnolog\'{\i}a, grant MTM2017-86061-C2-2-P.} \ and 
Subhajit Dutta
\footnote{This author has been partially supported by the
DST-SERB grant ECR/2017/000374.}\\
\vspace{0.1in}
Department of Mathematics, Statistics and Computation,
University of Cantabria, Spain\\
Department of Mathematics and Statistics, IIT Kanpur, India.
}
\end{center}

\baselineskip=17pt \vskip 3mm \centerline{\today} \vskip 3mm

\begin{center}
{{\bf Abstract}}
\end{center}
In this paper, we propose a data based transformation for infinite-dimensional Gaussian processes and derive its limit theorem.
For a classification problem, this transformation induces complete separation among the associated Gaussian processes. The misclassification probability of any simple classifier when applied on the transformed data asymptotically converges to zero. In a clustering problem using mixture models, an appropriate modification of this transformation asymptotically leads to perfect separation of the populations.
Theoretical properties are studied for the usual $k$-means clustering method when used on this transformed data.
Good empirical performance of the proposed methodology is demonstrated using simulated as well as benchmark data sets, when compared with some popular parametric and nonparametric methods for such functional data.

\baselineskip=14pt \vskip 5mm

\par\medskip\noindent

\baselineskip=26.5pt

\large

\section{Introduction} \label{Intro}

Let us consider processes, \Zvec, defined on a bounded real interval, which without loss of generality, we identify with the unit interval $[0,1]$. We assume that their trajectories belong to the Hilbert space of square integrable functions $\Hil$ which is defined as follows: 
{$$\Hil: \mbox{ set of real functions } f(t) \mbox{ with } t \in [0,1] \mbox{ such that } \int_0^1 f^2(t)dt <\infty.\vspace{-0.15in}$$}
The inner product  in \Hil \ is $\langle f,g \rangle = \disp \int_0^1 f(t)g(t) dt$. 

The keystone of this paper is Theorem \ref{Thm1}. It states that under appropriate assumptions, if \Zvec \ is a Gaussian process (GP) with values in  \Hil \ and $b \in \Hil$, then the limit of a sequence of scaled Mahalanobis distances between some finite-dimensional projections of $\Zvec$ and $b$ converges in probability to a non-random limit. Scaling is done using the dimension of the projection, and this convergence holds as the dimension goes to infinity.

From the point of view of applications, the interest of this result lies in the fact that this limit depends on the distribution of \Zvec \ (say, \ProbZ). Let us assume that we have an observation $\zvec \in \Hil$ and two Gaussian distributions $\Prob_1$ and $\Prob_2$ on \Hil \ with known parameters. We need to decide which Gaussian distribution (GD) generated \zvec. Assume that $\Prob_1$ and $\Prob_2$ produce different limits in Theorem \ref{Thm1} (more precisely, see Theorem \ref{Prop1}, which is an immediate consequence of this fact). Now, if we fix a large enough dimension and compute the value of the corresponding Mahalanobis distance by replacing the generic \Zvec \ with \zvec, then we will obtain a value close to the limit associated with the same probability which generated \zvec \ with a high probability. 

This theoretical result can be extended to a clustering problem in the mixture setting. Let us now assume that we have a probability distribution \Prob \ such that $\Prob= \sum_{h=1}^J\pi_h\Prob_{h}$, where $0<\pi_h<1$, $\sum_{h=1}^J\pi_h=1$ and $\Prob_{h}$ are GDs on $\mathbb{H}$ for $h=1,\ldots,J$. Additionally, we assume \Prob \ to be known, but have no information on the precise values of $J$, $\pi_h$ and $\Prob_h$ for $h=1,\ldots,J$.
According to the model, every function \zvec \ produced by \Prob \ was in fact  produced by one of the $\Prob_h$'s. Consequently, under appropriate conditions, this function would give different limits in Theorem \ref{Prop1} depending on the value of $h$. 
So, if we have at least one observation from each $\Prob_h$, we can identify the value of $J$ as well as the subsets of observations produced by the same $\Prob_h$.

The main problem with this idea is that computation of the involved limits in Theorem \ref{Prop1} requires knowledge of the means and the covariance operators of the probabilities $\Prob_h$ for $h=1,\ldots,J$, which is not available in this setting.  
We fix this issue in Theorem \ref{Prop2}, which also provides a modification of Theorem \ref{Thm1}. Given two random functions $\Zvec_1$ and $\Zvec_2$ produced by $\Prob_{h_1}$ and $\Prob_{h_2}$, respectively, the limit in probability of another sequence of scaled Mahalanobis distances between some finite-dimensional projections of $\Zvec_1$ and $\Zvec_2$ is non-random, and again depends on the indices ${h_1}$ and $h_2$.
Under appropriate assumptions, if we have a random sample taken from \Prob \ which contains at least an observation from each population, then the limits (on the dimension) of those Mahalanobis distances between pairs of points in the sample allows one to determine the value of $J$ as well as the points in the sample which were generated by the same $\Prob_h$, without any possibility of mistake.
Additionally, computation of the Mahalanobis distances involves only the covariance operator of \Prob .

In summary, Theorem \ref{Thm1}  (and its above mentioned consequences) identifies situations in which two GDs are {\it `mutually singular'}. The singularity of GDs is related to the Hajek and Feldman  property, HFp, (established independently by \cite{Hajek} and \cite{Feldman}). The HFp states that if $\Prob_1$ and $\Prob_2$ are GDs, then they are either equivalent, or mutually singular. In other words, for every measurable set $A$, $\Prob_1(A)=0$ if and only if $\Prob_2(A)=0$, or else there exist two disjoint measurable sets $S_{1}$ and $S_{2}$  with
$$\Prob_1(S_{1})=1, \Prob_2(S_{1})=0 \mbox{ and } \Prob_1(S_{2})=0,\Prob_2(S_{2})=1.$$

Mutual singularity is not very interesting in finite dimensions because it happens only when at least one of the covariance matrices is singular. However, in the functional case, this singularity appears in non-trivial situations. Just to mention an example, it was shown in \cite{Rao_Varadarajan_1963} that if the covariance operators of $\Prob_1$ and $\Prob_2$, say, $\Sigma_1$ and $\Sigma_2$ satisfy $\Sigma_2 = a \Sigma_1$ for some $a \neq 1$, then $\Prob_1$ and $\Prob_2$ are mutually singular.

It seems that the HFp should have attracted the attention of researchers in  classification and clustering for functional data, the orthogonality case apparently being more attractive because it would allow one to obtain `{\it perfect classification}' and  `{\it perfect clustering}'. 
To the best of our knowledge nobody used these results in practice till \cite{Baillo_2011}, where the authors derived a classification procedure using likelihood ratios. They focus on the equivalence case
and hence, do not obtain perfect classification.
Optimal classification of GPs was analyzed in \cite{TorrecillaEtA_l2020} from the HFp viewpoint. In this paper, the optimal (Bayes') classifier of equivalent GPs was derived and a procedure to obtain asymptotically perfect classification of mutually singular GPs was described as well. The results covered both homoscedastic and heteroscedastic cases.

Additionally, \cite{Delaigle_Hall_2012} and \cite{Delaigle_Hall_2013} investigated conditions under which a perfect classification procedure for GPs was possible and developed related classifiers. The paper by \cite{DMY_2017} proposed a functional classifier based on ratio of density functions, which also leads to perfect classfication. These papers contain no reference to the HFp. In fact, the relationship between \cite{Delaigle_Hall_2012} and the HFp was analyzed in \cite{BCT_2018}, where the authors presented an expression of the optimal Bayes' rule in some classification problems.  
Research on perfect  clustering of functional data is not so abundant. 
As far as we know, \cite{Delaigle_Hall_Pham_2019} is the only available paper with results in this field.



The methods we propose in this paper differ from \cite{TorrecillaEtA_l2020}, where one needs to have some information on the mutual singularity of $\Prob_1$ and $\Prob_2$, and as far as we know, have not been applied to clustering. However, our procedure  can be applied  with no previous information. In exchange, our results do not cover all possible mutually singular measures, but only those which satisfy the appropriate limit conditions, which in turn are general enough as to cover a wide range of situations.

In \cite{Rao_Varadarajan_1963} and \cite{Shepp_2} the authors obtain characterizations of the singularity, or equivalence of Gaussian measures in functional spaces. Their results also involve increasing sequences of subspaces. For equivalent GDs, the limit obtained in \cite{Rao_Varadarajan_1963} includes a term which is the exponential of an expression involving the difference of the means of $\Prob_1$ and $\Prob_2$; curiously, the logarithm of this term is related with our limit expression.
Similarities between our proposal and those in \cite{Rao_Varadarajan_1963} and \cite{Shepp_2} end here because the other three involved terms are different. Additionally, we handle Mahalanobis distances {\it between data points}, while those papers use Hellinger and Jeffreys functionals to measure discrepancy {\it between distributions}. 
As a consequence, the characterizations they obtain are not applicable in practice to classify (or, cluster) data points because they depend on the full distribution, and it is not straight forward to compute the related functionals using data points.

%
%
%
%
%
%

In this paper, we first analyze the limit of the above mentioned scaled Mahalanobis' distances by assuming the underlying parameters of the GPs to be {known} in Section \ref{Trans.GP}. We begin with a general concentration result (Theorem \ref{Thm1}).
Based on this limit theorem, we translate the mutually singular case in the HFp to the idea of perfect classification (Theorem \ref{Thm2}) and perfect clustering of GPs (Theorem \ref{Theo.Cluster.M1}) by proposing some transformations for classification as well as clustering that asymptotically yield `{\it perfect separation among the classes or clusters}' for $J$ class problems. This transformation can also be used to find the unknown number of clusters (Proposition \ref{Prop.Number.Clusters}).
In Section \ref{Sec.Estimation}, we estimate the underlying parameters (the means and covariance operators of the involved distributions in the classification case, and the covariance operator of the mixture distribution in the clustering problem) from data, and state related asymptotic results of the proposed transformations. More specifically, we propose a simple transformation to directly address the general $J$ class classification problem of GPs, and prove related consistency results for the empirical case (Theorem \ref{Theo.Classif.2classes}).
We also prove uniform (on the sample points) consistency of the empirical version (Theorem \ref{Theo.Cluster}) for the transformation associated with GP clustering.
It is surprising that our GP clustering method fails to discriminate `location only' scenarios, but yields perfect clustering if differences in scales (see Remark \ref{Remark.UselessLocation}) exist. A possible way to fix this is discussed in Appendix II.

We have also compared our work both theoretically (see Section \ref{DellaigleHall} of Appendix II) and numerically (see Sections \ref{Num.Res.Simulation} and \ref{Sec.BenchData}) with the existing literature on perfect classification and clustering for functional data.
All proofs are deferred to Appendix I. 
Appendix II 
also contains some additional material including a possible extension to non-Gaussian distributions; 
theoretical comparisons of our results with those obtained in the papers by \cite{Delaigle_Hall_2012}, \cite{Delaigle_Hall_2013}, \cite{Delaigle_Hall_Pham_2019} and \cite{TorrecillaEtA_l2020}; the proof of a technical lemma and some additional simulation results.

In this paper, we will use the following notation. The distribution of the random process ${\bf Z}$ will be denoted as $\Prob_{\bf Z}$, its mean function by $\mu^{\bf Z}$ and its covariance operator (referred to simply as covariance) by $\Sigma^{\bf Z}$. We will write $\Sigma^{\bf Z}(s,t)$ to denote the covariance between $\Zvec(s)$ and $\Zvec(t)$ for $s,t \in [0,1]$.
Further, we will assume that all involved random quantities are defined on a common probability space $(\Omega, \mathcal{A}, \Prob)$. Given a square matrix $A$, $trace(A)$ will denote its trace. The usual Euclidean norm on $\Rea^d$ is denoted by $\|\cdot\|$. 


\section{Transformation with Known Distributions} \label{Trans.GP}

%

Let $\{V_d\}_{d\in \mathbb{N}}$ (with $V_d \subset \Hil$) be an increasing sequence of subspaces. Here, the dimension of $V_d$ is $d$. This restriction is not necessary for the development which follows as long as the dimension of $V_d$ goes to infinity with increasing $d$, but it simplifies the notation. 
Given the subspace $V_d$, let $\mu^{\bf Z}_d$ and $\Sigma^{\bf Z}_d$ represent the $d$-dimensional mean and the $d \times d$ covariance matrix of the projection of $\Zvec$ on this subspace. 
If $\bf u \in \Hil$, we denote ${\bf u}_d$ to be its projection on $V_d$. 

Fix $b\in \Hil$. Theorem \ref{Thm1} analyses the behaviour of the limit of squared Mahalanobis norm of the random vector $(\Zvec - b)_d$ for $d \in \mathbb{N}$. For every positive definite $d\times d$ matrix $A_d$, we define the map

\begin{equation}\label{Eq.def.D_d}
D_{d}^{A_d}({\uvec},\vvec) =  \frac 1 d \left \|A_d^{-1/2}({\uvec}-\vvec)_d \right \|^2 \mbox{ for} \ \uvec,\vvec \in \Hil .
\end{equation}

\noindent
In this section, the underlying distributions are assumed to be known.
After stating Theorem \ref{Thm1} and some remarks related to it, we will look into two applications inspired from this result. We will take advantage of the fact that the limit in this theorem is not random, but it may depend on the underlying probability distribution $\Prob_{\bf Z}$.


\begin{Theo} \label{Thm1}
Let $\{A_d\}$ be a sequence of $d\times d$ symmetric, positive definite matrices and $\alpha_1^d, \ldots , \alpha_d^d$ be the eigenvalues of the matrix $S_d = (A_d)^{-1/2}  \Sigma^{\bf Z}_d (A_d)^{-1/2}$ for $d \in \mathbb{N}$. We define $\alpha_d=(\alpha_1^d,\ldots,\alpha_d^d)^T$ and $\|\alpha_d\|_{\infty}=\max (\alpha_1^d,\ldots,\alpha_d^d)$ is the supremum norm.
Let $b \in \Hil$ such that  there exist constants $L_\mu$ and $L_S$ (finite, or not) with
\begin{eqnarray} \label{A1.Thm1}
L_\mu&= &\lim_{d \to \infty} D_{d}^{A_d} (\mu^{\bf Z} ,b),
\\
\label{A2.Thm1}
L_S & = &\lim_{d \to \infty} \frac 1 {d} \mbox{ trace}(S_d), \mbox{ and}
\\
\label{A3.Thm1}
0 &=&\lim_{d \to \infty} \frac{\|\alpha_d\|_{\infty}}{d}.
\end{eqnarray}

\noindent
Then 
$
D_d^{A_d}({\bf Z},b)  \stackrel{P}{\conv}  L:= L_\mu+L_S  \mbox{ as } d \conv \infty.
$
\end{Theo}

\begin{Nota} \label{R1}
{\rm
A condition in Theorem \ref{Thm1} is required to ensure that no single component is extremely influential. For instance, it may happen that we take a sequence such that $\alpha^d_{1} = d$ and $\alpha^d_{i} = o(d^{-1})$ for every $2 \leq i \leq d$. Under this condition, no limit is possible in Theorem \ref{Thm1}. However, this possibility is excluded by assumption (\ref{A3.Thm1}). 
}
\end{Nota}

\begin{Nota} \label{R2}
{\rm
We allow both  constants in Theorem \ref{Thm1} to be infinite. When $L_S$ is finite, Lemma \ref{Lemm2} (in Appendix I) shows that assumption (\ref{A3.Thm1}) follows from assumption (\ref{A2.Thm1}). 
}
\end{Nota}

\begin{Nota} \label{R3}
{\rm
Let $\Zvec_{1}$ and $\Zvec_2$ be independent observations generated from GDs $\Prob_{1}$ and $\Prob_{2}$. Thus, ${\bf Z}_1 - {\bf Z}_2$ is a GP with mean $\mu_1-\mu_2$ and covariance $\Sigma^{\bf Z_1}+\Sigma^{\bf Z_2}$. Consider the matrix $S_d = (A_d)^{-1/2} (\Sigma^{\bf Z_1}_d+\Sigma^{\bf Z_2}_d) (A_d)^{-1/2}$ with $d \in \mathbb{N}$. Take ${\bf Z}={\bf Z}_1-{\bf Z}_2$ and $b=0$ in Theorem \ref{Thm1}. Then, the following convergence result holds:
$$
D_d^{A_d}(\Zvec_1,\Zvec_2)=D_d^{A_d}(\Zvec_1-\Zvec_2,0) \stackrel{P}{\conv} L:=L_{\mu}+L_S \mbox{ as } d \conv \infty.
$$
Here, $L_{\mu}=\lim_{d \to \infty} D_{d}^{A_d} (\mu_1, \mu_2)$ and $L_S$ is as defined in equation (\ref{A2.Thm1}) of Theorem \ref{Thm1}.
}
\end{Nota}

\begin{Nota} \label{R4}
{\rm 
In general, the fact that $V_d \subset V_{d+1}$ does not guarantee the  existence of any relationship between the sets $\{\alpha^d_{1}, \ldots, \alpha^d_{d}\}$ and $\{\alpha^{d+1}_{1}, \ldots, \alpha^{d+1}_{d+1}\}$. However, in some cases $\{\alpha^d_{1}, \ldots, \alpha^d_{d}\} \subset \{\alpha^{d+1}_{1}, \ldots, \alpha^{d+1}_{d+1}\}$ (see, for instance, Section \ref{Class.GP.Examples}, where $V_d$ is generated by the first $d$ eigenfunctions of $\Sigma$ and $A = a\Sigma$ for some $a>0$). 
}
\end{Nota}
%
%

\subsection{Application I: Supervised Classification} \label{Class.GP}

Assume that we have two $\Hil$-valued Gaussian processes associated with probability distributions $\Prob_1$ and $\Prob_2$. We denote their means by
$\mu_1$ and $\mu_2$, and their covariances to be $\Sigma_1$ and $\Sigma_2$. 
An observation $\zvec \in \Hil$ has been drawn either from $\Prob_1$ or $ \Prob_2$, and our aim is to discriminate between those two possibilities. 
We will write $D_{d}^{i}({\uvec},\vvec)$ instead of $D_{d}^{\Sigma_{id}}({\uvec},\vvec)$  for $i=1,2$.
The classification procedure that we propose is based on the behavior of $D_{d}^{1}(\zvec,\mu_{1})$ and $D_{d}^{2}(\zvec,\mu_2)$ depending on the distribution of $\zvec$ (i.e., either $\Prob_1$ or $\Prob_2$). This analysis is carried out in Theorem \ref{Prop1} below. 

\begin{Theo} \label{Prop1}
Assume that $\ProbZ= \Prob_1$. Then, it happens that
\begin{equation}\label{Eq.limit1}
D_{d}^{1} (\Zvec,\mu_{1}) \stackrel{P}{\conv} 1~\mbox{as}~d \conv \infty.
\end{equation}
Define the matrix $S_d^{12}= (\Sigma_{2d})^{-1/2}  \Sigma_{1d}$ $(\Sigma_{2d})^{-1/2}$ with $d \in \mathbb{N}$. Moreover, if there exist constants 
$L_\mu^{12}= \lim_d  D_{d}^{2} (\mu_{1},\mu_{2})$
and
$L_S^{12} = \lim_d \frac 1 {d}\mbox{ trace}(S_d^{12})$, and assumption (\ref{A3.Thm1}) holds for the sequence of matrices $S_d^{12}$ with $d \in \mathbb{N}$, 
then 
\begin{equation}\label{Eq.limit2.2}
D_{d}^{2} ({\bf z},\mu_{2})  \stackrel{P}{\conv} L^{12}:= ( L^{12}_\mu +  L^{12}_S) ~\mbox{as}~ d \conv \infty.
\end{equation}
\end{Theo}

\begin{Nota} \label{R5}
{\rm Clearly, Theorem \ref{Prop1} holds if we replace $\Prob_1$ by $\Prob_2$. Consider the matrix $S_d^{21}= (\Sigma_{1d})^{-1/2}  \Sigma_{2d}$ $(\Sigma_{1d})^{-1/2}$ with $d \in \mathbb{N}$ and $L^{21}=L_\mu^{21} +L_S^{21}$, where 
$L_\mu^{21}= \lim_d D_{d}^{1} (\mu_{1},\mu_{2})$ and $L_S^{21} = \lim_d \frac 1 {d}\mbox{ trace}(S_d^{21})$.
Under the corresponding assumptions, for $\ProbZ=\Prob_2$, we get
$$D_{d}^{1} (\Zvec,\mu_{1}) \stackrel{P}{\conv} L^{21} ~\mbox{and}~ D_{d}^{2} (\Zvec,\mu_{2})\stackrel{P}{\conv} 1 ~\mbox{as}~ d \conv \infty.$$
}
\end{Nota}

In order to apply Theorem \ref{Prop1} to classification problems, we define a sequence of two-dimensional transformations as follows:
\begin{equation} \label{Transf.Class}
{\bf T}_d(\zvec)=(D_{d}^{1}(\zvec,\mu_{1}), D_{d}^{2}(\zvec, \mu_{2}))^T \mbox{ for } d \in \Nat.
\end{equation}
Based on $\Tvec_d$, let us consider the following classifier
\[
\Psi_{d}(\zvec) =
\left\{
\begin{array}{rl}
1, & \mbox{ if } ~\| \Tvec_{d}(\zvec) - (1,L^{12})^T \| \leq \| \Tvec_{d}(\zvec) - (L^{21},1)^T \|,\\[2mm]
2, & \mbox{ else}.
\end{array}
\right.
\]

\begin{Defi}
For a fixed value of $d \in \mathbb{N}$, we define $p_d$ to be the total misclassification probability of the classifier $\Psi_d$, i.e., $p_d=\pi_{1} \Prob_1[\Psi_d = 2]+\pi_{2}  \Prob_2[\Psi_d= 1]$ with $\pi_{1}$  (respectively, $\pi_{2}$) being the prior probability corresponding to $\Prob_1$ (respectively, $\Prob_2$). 

\noindent
We assume that $0<\pi_{1} <1$ with $\pi_{1}+ \pi_{2}=1$.
\end{Defi}

\begin{Theo} \label{Thm2}
If $L^{12}$ and $L^{21}$ exist, with $L^{12} \neq 1$ or $L^{21} \neq 1$, and assumption (\ref{A3.Thm1}) holds for the sequences of matrices $S_d^{12}$ and $S_d^{21}$ with $d \in \mathbb{N}$, then the misclassification probability $p_{d} \to 0$ as $d \to \infty$.
\end{Theo}

A consequence of Theorem \ref{Thm2} is that under appropriate assumptions (which are stated in this result), we can now identify the distribution which generated ${\bf z}$ without a possibility of mistake as $d \to \infty$, i.e., we obtain asymptotic `{\it perfect classification}'.
For a $J$ class problem, the $2$-dimensional transformation ${\bf T}_d(\zvec)$ proposed in (\ref{Transf.Class}) can be easily generalized to $J$ dimensions by computing the Mahalanobis distances w.r.t. the $J$ competing class distributions. If $\ProbZ= \Prob_j$ with $1 \leq j \leq J$, then the limit of ${\bf T}_d(\Zvec)$ as $d \to \infty$ can be derived in a similar way.
The classifier $\Psi_{d}(\zvec)$ can be constructed by classifying $\zvec$ to the class which corresponds to the minimum distance between ${\bf T}_d(\zvec)$ and the associated vector of constants. Further, the consistency result in Theorem \ref{Thm2} will continue to hold under analogous conditions on the related constants.

\subsubsection{Some Example of GPs} \label{Class.GP.Examples}

We now analyze the limiting behavior of the transformation ${\bf T}_d(\Zvec)$ in the particular case when $\mu_{1}=0$ and $\Sigma_2 = a \Sigma_{1}$, where $a \in \Rea^+$. Let $V_d$  be the subspace generated by the first $d$ eigenfunctions of the covariance  $\Sigma_{1}$. Assume that there exists  $\nu=\lim_{d \to \infty} D_d^1(0,\mu_2)=\lim_{d \to \infty} \frac{1}{d} {\|({\Sigma_{1d}})^{-1/2}\mu_{2d}\|^2}$. This  implies that we have definite expressions for the limiting constants mentioned in equations (\ref{Eq.limit1}) and (\ref{Eq.limit2.2}), and Remark \ref{R5} as follows:
$$L_\mu^{12}=\nu/a, L_S^{12}=1/a~ \mbox{and} ~L_\mu^{21} =\nu, L_S^{21}=a.$$
The following simplified limits hold (as $d \to \infty$):
\[
 \Tvec_d(\Zvec) \stackrel{P}{\to}
\left\{
\begin{array}{ll}
(1,(\nu+1)/{a})^T,  & \mbox{ if }\ProbZ= \Prob_1,
\\
[2mm]
(\nu+a, 1)^T, & \mbox{ if }  \ProbZ= \Prob_2.
\end{array}
\right.
\]
Let us consider two special cases.

\subsubsection{Homoscedastic case} \label{Class.GP.Traslation}

Let us assume that $\mu_{2} \neq 0$, and $a=1$. So, we have a common covariance for both classes. In such a case, we obtain the following (as $d \to \infty$):
\[
 \Tvec_d(\Zvec) \stackrel{P}{\to}
\left\{
\begin{array}{rl}
(1, \nu + 1)^T, & \mbox{ if } \ProbZ= \Prob_1,
\\
[2mm]
(\nu+1, 1)^T, & \mbox{ if } \ProbZ= \Prob_2.
\end{array}
\right.
\]

\noindent 
If $\nu > 0$, this transformation allows us to identify the distribution which produced $\Zvec$. However, when $\nu=0$, then it happens that $\Tvec_d(\Zvec) \stackrel{P}{\to} (1,1)^T$ independently of $\ProbZ= \Prob_1$, or $\ProbZ= \Prob_2$. In other words, the proposed transformation allows one to decide the distribution which produced $\zvec$ without any possibility of mistake, or alternatively, it is completely useless. 

\subsubsection{Equality in means, and difference only in scale}  \label{Class.GP.Scale}

We will now assume that $\mu_{2}={0}$ and $a \neq 1$. Therefore, $\mu_{1}=\mu_{2}$ and this implies that $\nu=0$. In this setting, we obtain the following limits (as $d \to \infty$):
\[
\Tvec_d(\Zvec) \stackrel{P}{\to}
\left\{
\begin{array}{ll}
(1, a^{-1})^T,  & \mbox{ if } \ProbZ= \Prob_1,
\\
[2mm]
(a, 1)^T, & \mbox{ if } \ProbZ= \Prob_2.
\end{array}
\right.
\]
This is the simplest case in which two different covariances give way to a perfect classification problem, but one may easily construct more involved situations.

\subsection{Application II: Unsupervised Classification} \label{Cluster.GP}

In this subsection, we handle a random {
function} $\Zvec$ whose distribution is a two component mixture distribution of the form: $\ProbZ= \pi_1\Prob_{1} + \pi_2\Prob_{2}$, where $0<\pi_1<1$ and $\pi_1+\pi_2=1$. Here, $\Prob_{h}$ denotes the GD on $\mathbb{H}$ with mean function $\mu_h$ and covariance $\Sigma_h$ for $h=1,2$. The  mean function and the  covariance of the mixture satisfy that $\mu^{\bf Z}(t)=\pi_1 \mu_1(t)+\pi_2 \mu_2(t)$ with $t \in [0,1]$ and 
\begin{equation} \label{Pooled.Sig}
\SigmaZ(s,t)= \pi_1 \Sigma_1(s,t) + \pi_2 \Sigma_2(s,t) + \pi_1 \pi_2 [\mu_1(s) - \mu_2(s)][\mu_1(t) - \mu_2(t)], ~s,t \in [0,1].
\end{equation}

\noindent
Given a random sample $\Zvec_1,\ldots,\Zvec_N$ from \ProbZ, consider the following set: 
\[
{\cal C}_h = \{j: \Zvec_j \mbox{ was obtained from } \Prob_h \mbox{ for } j=1,\ldots,N \}
\]
with $h \in \{1,2\}$. The components of the mixture distribution \ProbZ \ and the sets ${\cal C}_h$ for $h=1,2$ are unknown, and the problem we are dealing with is the estimation of these sets. However, we assume \ProbZ \ and the sets ${\cal C}_h$ for $h=1,2$ to be known in this section to build the fundamental idea behind using the proposed transformation for GP clustering.

Let $V_d$ with $d \in \mathbb{N}$ denote the sequence of $d$-dimensional subspaces generated by the $d$ eigenfunctions associated with the $d$ largest eigenvalues of $\SigmaZ$ (recall our framework in Section \ref{Class.GP} and the discussion in Remark \ref{R4}).
In the following result, $\Zvec_1$ and $\Zvec_2$ are assumed to be independent and $\Prob_{{\bf Z}_1}= \Prob_{h}$ and $\Prob_{{\bf Z}_2}= \Prob_{k}$ with $h,k \in \{1,2\}$. The clustering procedure that we propose is based on the behavior of the transformation $D_{d}^{\Sigma^{\bf Z}}(\Zvec_1,\Zvec_2)$, which is stated below in Proposition \ref{Prop2}. Recall the notation introduced in Theorem \ref{Prop1}.

The structure of the covariance  $\SigmaZ$ stated in equation (\ref{Pooled.Sig}) imposes some restrictions on the associated constants as stated in part (c) of Theorem \ref{Prop2} below. In particular, the fact that $L_S^{h}$ and $L_S^{hk}$ are finite implies that assumption (\ref{A3.Thm1}) in Thoerem \ref{Prop2} always holds for the sequence of matrices $\{S_d^h\}_{d \in \mathbb{N}}$ and $\{S_d^{hk}\}_{d \in \mathbb{N}}$ for $h,k \in \{1,2\}$.

\begin{Theo} \label{Prop2}
(a) Assume that $h=k \in \{1,2\}$. Define $S^{h}_d := (\Sigma_d)^{-1/2} (2\Sigma_{{h}d}) (\Sigma_d)^{-1/2}$ for $d \in \Nat$, and assume that $L_S^{h}  = \lim_d \frac 1 {d} trace(S^{h}_d)$ exists.
Then,
\begin{equation} \label{Cluster.GP.E1}
D^{\Sigma_d^{\bf Z}}_{d} (\Zvec_1, \Zvec_2) \stackrel{P}{\to} L_S^{h} \mbox{ as } d \to \infty.
\end{equation}
(b) Assume that $h \neq k \in \{1,2\}$. Define $S^{hk}_d:=(\Sigma_d)^{-1/2}  (\Sigma_{{h}d}+\Sigma_{{k}d}) (\Sigma_d)^{-1/2}$ for $d \in \Nat$, and assume that $L_S^{hk}  = \lim_d \frac 1 {d} trace(S_d^{hk})$  exists.
Then,
\begin{equation} \label{Cluster.GP.E2}
D^{\Sigma^{\bf Z}_d}_{d} (\Zvec_1, \Zvec_2) \stackrel{P}{\to}L^ {hk}:= L_S^ {hk} \mbox{ as } d \to \infty.
\end{equation}
(c) If $h \neq k \in \{1,2\}$,  then $L_\mu^{hk}  = \lim_d D_d^{\Sigma^{\bf Z}_d}(\mu_{{h}}, \mu_{{k}})=0$ and both $L_S^ {h} $ and $L_S^ {hk} $ are finite.
\end{Theo}


\begin{Nota}\label{Remark.UselessLocation}
{\rm
It follows from (c) in Theorem \ref{Prop2} that the statistic we propose is useless for unsupervised classification in the homoscedastic case (independently of the difference between  $\mu_1 $ and $ \mu_2$) because if $\Sigma_1 = \Sigma_2$, then $L^{12}=L^{1}_S=L^{2}_S$.
A possibility is to replace the statistic $D^{\Sigma_d^{\bf Z}}_{d} ({\bf z}_1, {\bf z}_2)$ by another one in which the value of the transformation $D^{\Sigma_d^{\bf Z}}_{d}({\mu}_1, {\mu}_2)$ increases with $d \in \mathbb{N}$. Our proposal is to use  
\[
D^{\Sigma_d,r}_{d} ({\bf u}, {\bf v}) := \frac 1 d \left \|(\Sigma_d^{-1/2})^r({\uvec}-\vvec)_d \right \|^2 
=\frac 1 d \sum_{i=1}^d \frac{(u_i-v_i)^2}{\lambda_i^r} \ \mbox{with } r \in \mathbb{N}.
\]
Discussion of this transformation, and some numerical results are included in Appendix~II.
}
\end{Nota}

To simplify notation and avoid technicalities with empty classes, we additionally assume
that the observations whose indices belong to the sets $\mathcal{C}_1^N=\{1,\ldots,N_1\}$ and $\mathcal{C}_2^N=\{N_1+1,\ldots,N\}$ with $N=N_1+N_2$ and $ N_1,N_2>0$, were generated by $\Prob_{1}$ and $\Prob_{2}$, respectively. 
In practice, these sets are unknown and in fact our aim is their estimation. We begin with this simplifying assumption for ease of notation, and to obtain a clear exposition of the proposed methodology.

Define the $N\times N$ matrix 
$\Gamma_{d}$ whose $(i,j)$-th element is
\begin{equation} \label{Dist2.Cluster}
\Gamma_{d}(\Zvec_i,\Zvec_j)=\gamma_{ij}^d
=
\frac{1}{N-2}\sum_{t=1, ~t\neq i,j}^{N} [D^{\Sigma_d^{\bf Z}}_{d}(\Zvec_t,\Zvec_i) - D^{\Sigma_d^{\bf Z}}_{d}(\Zvec_t,\Zvec_j)]^2
\end{equation}
for $1 \leq i,j \leq N$.
Theorem \ref{Prop2} and the fact that $t\neq i,j$ in \eqref{Dist2.Cluster} give us that
\begin{equation} \label{Eq.Conv.Deltas}
\gamma_{ij}^d \stackrel{P}{\to}
\left \{ \begin{array}{cl}
0 & ~~\mbox{if}~ i,j \in \mathcal{C}_h~\mbox{for}~h=1,2,\\[1mm]
\gamma_{hk} & ~~\mbox{if}~ i \in \mathcal{C}_{h}~\mbox{and}~j \in \mathcal{C}_{k}, \mbox{ with } h\neq k \in \{1,2\},
\end{array} \right.
\end{equation}
as $d \to \infty$, where 
$$\gamma_{hk}=\frac{N_{h}-1}{N-2}(L_{S}^{h}-L_S^{hk})^2 + \frac{N_{k}-1}{N-2}(L_{S}^{k}-L_S^{kh})^2.$$


\noindent
Combining the fact stated above in (\ref{Eq.Conv.Deltas}), as $d \to \infty$, we obtain
\begin{equation} \label{Eq.Estr_Delta}
\Gamma_{d} \stackrel{P}{\to} \Gamma :=
\left[ \begin{array}{cc}
\vspace{.1in}
{\bf 0}_{N_1}{\bf 0}_{N_1}^T & \gamma_{12} {\bf 1}_{N_1}{\bf 1}_{N_2}^T \\
\gamma_{21} {\bf 1}_{N_2}{\bf 1}_{N_1}^T & {\bf 0}_{N_2}{\bf 0}_{N_2}^T\\
\end{array} \right].
\end{equation}


\noindent
Let $\beta_i^d$ and $\beta_i$ (with $1 \leq i \leq N$) denote the eigenvalues corresponding to the matrices $\Gamma_{d}$ and $\Gamma$, respectively. Define the following quantities
\begin{equation} \label{Eq.Cluster.Number}
K_d = \sum_{i=1}^N I(|\beta_i^d| > a_d) \mbox{ and } K_0=\sum_{i=1}^N I(|\beta_i| > 0), 
\end{equation}
with $\{a_d\}_{d \in \mathbb{N}}$ decreasing to $0$ as $d \to \infty$ at an appropriate rate, and $I$ is the indicator function. The constant $K_0$ clearly equals $2$ for the limiting $N \times N$ matrix $\Gamma$ stated in (\ref{Eq.Estr_Delta}), and hence, correctly identifies the true underlying number of clusters. 

\begin{Prop} \label{Prop.Number.Clusters}
Assume $N_1, N_2 \geq 1$ are fixed. Under the assumptions of Theorem \ref{Prop2}, with $L_S^{12} \neq L_S^1$ and $L_S^{21} \neq L_S^2$, there exists a sequence $\{a_d\}_{d \in \mathbb{N}} \subset \Rea^+$ such that $a_d \to 0$ and $K_d \stackrel{P}{\to} 2 \mbox{ as } d \to \infty$.
\end{Prop}

This now implies that we can correctly identify the true number of clusters asympototically, as $d \to \infty$. Note that the structure of the martix $\Gamma$ in \eqref{Eq.Estr_Delta} is straight forward because of the simplifying assumption on the sets $\mathcal{C}_1$ and $\mathcal{C}_2$. However, this is not a requirement and we will drop it.
Proposition \ref{Prop.Number.Clusters} holds more generally for {\it any permutation} of the data points $\Zvec_1,\ldots,\Zvec_N$. In fact, if the sets $\mathcal{C}_1$ and $\mathcal{C}_2$ are unknown, then the rows/columns of the $\Gamma$ matrix will be permuted accordingly. But, the underlying structure will remain the same and Proposition \ref{Prop.Number.Clusters} will continue to hold.
As a followup, we now prove that if any standard clustering method is used on the $\Gamma_d$ matrix, then we can perfectly cluster all the observations asympototically (as $d \to \infty$) because of the structure of the $\Gamma$ matrix stated in \eqref{Eq.Estr_Delta}.

\begin{Defi}
A clustering method can be defined as a map from $\Hil$ to the set $\{1,\ldots, J\}$. Consider two such maps $\psi_{d}$ for a fixed $d \geq 1$ and $\phi$. A measure of distance between two clusterings based on the Rand index (see p. 847 of \cite{Rand_1971}) is defined as follows:
$$\mathbb{R}_{d,N} = \frac{1}{\binom{N}{2}} \sum_{1 \leq i<j \leq N} I \biggl [ I[\psi_{d}(\zvec_i) = \psi_{d}(\zvec_j)] + I[\phi(\zvec_i) = \phi(\zvec_j)] = 1 \biggr ],
$$
for a fixed $N \geq 2$.
\end{Defi}

\noindent
Let $\phi$ be the map giving the true labels of the data, i.e., $\phi(\xvec_j)=h$ for $j \in \calC_h$ and $h \in \{1,2\}$. We can now construct $\psi_{d}$ based on the data by directly applying any clustering technique. Here, we use the $k$-means algorithm on the rows (or, columns) of the matrix~$\Gamma_d$. 

Mathematically, the $k$-means algorithm finds $J$ groups (say, $\calG_1,\ldots,\calG_J$) with centers $\cvec_1,\ldots,\cvec_J$ such that $\phi(\calG_1,\ldots,\calG_J)=\sum_{h=1}^J \sum_{i:\small \xvec_i \in \calG_h} \|\xvec_i-\cvec_h\|^2$ is minimized. 
The asymptotic properties of the matrix ${\Gamma}_{d}$ as $d \to \infty$ stated above in equation (\ref{Eq.Estr_Delta}) imply that differences in the constants should yield perfect clustering.
Our next result proves label consistency for this $k$-means algorithm when $J=2$.

\begin{Theo} \label{Theo.Cluster.M1}
Assume $J=2$ and $\gamma_{12}>0$. Further, assume that the conditions in Theorem \ref{Prop2} and Proposition \ref{Prop.Number.Clusters} hold. 
Then, the clusters will be perfectly identifiable, i.e., $\mathbb{R}_{d,N} \stackrel{P}{\to} 0$ as $d \to \infty$.
\end{Theo}

\begin{Nota}\label{RemarkOnJ}
{\rm
The structure of the $N \times N$ matrix $\Gamma$ stated in equation (\ref{Eq.Estr_Delta}) continues to hold, and will lead us to perfect clustering for every value of $J \ge 2$. Moreover, the procedure described in Proposition \ref{Prop.Number.Clusters} also works fine, with the limit equal to the rank of $\Gamma$. 
However, generalizing this idea to $J(>2)$ clusters is not trivial. 
		
\noindent
The quantity $K_0$ in equation \eqref{Eq.Cluster.Number} is the rank of $\Gamma$, and one may be tempted to think that it generally coincides with $J$. But, this is true only for $J \leq 3$ and may be different for $J \geq 4$ (as shown in Lemma \ref{J.clust} in Appendix I).
The proof of Lemma \ref{J.clust} further shows that the condition under which Rank$(\Gamma)< J$ is quite restrictive. Thus, in practice, our proposal is to estimate $J$ with $K_d$. 
}
\end{Nota}

\subsubsection{Example with GPs} \label{Gaussian.Cluster.Example}

Statement (c) in Theorem \ref{Prop2} implies that the matrix $\Gamma$ will be null for the homoscedastic case.
However, if we assume that $\Sigma_2=a\Sigma_1$ with $a>0$, then we have the following expressions for the scale constants stated in Theorem \ref{Prop2}:
$$L_{S}^{1}=\frac{2}{\pi_1+\pi_2a},~ L_{S}^{2}=\frac{2a}{\pi_1+\pi_2a}~ \mbox{and}~L_S^{12}=\frac {1+a}{\pi_1+a\pi_2}.$$
Thus, it is possible to identify perfectly the clusters as long as $a\neq 1$ since this implies that $\gamma_{12}$ and $\gamma_{21}$ both are positive quantities.


\subsubsection{Uniform Convergence}

In Theorem \ref{Prop2}, we have proved consistency for finite sets of  data points of the function $D_d^{{\Sigma}_d}(\Zvec_1,\Zvec_2)$ defined in (\ref{Eq.def.D_d}). In Theorem \ref{Theo.Conssit.Clusytering} below, we prove the uniform convergence of this function on the random sample as $N \to \infty$. This result will be useful in establishing a second result on uniform convergence, which is stated in the next section.

\begin{Theo} \label{Theo.Conssit.Clusytering}
Assume the conditions in Theorem \ref{Prop2}, and let $\{d_N\} \subset \Nat$ be such that $d_N \to \infty$. Then,
\begin{itemize}

\item[a)] For $h \in \{1,2\}$, 
let $\alpha_{d_N}=(\alpha^{d_N}_{1}, \ldots, \alpha^{d_N}_{{d_N}})^T$ be the eigenvalues of $S_{d_N}^h$ with ${d_N} \in \Nat$. 
If 
\begin{equation}\label{Hip.Theo.Conssit.Clusytering.Ass.a}
\log N = o \left(\frac{d_N}{\|\alpha_{d_N}\|_{\infty}}\right), 
\end{equation}
then it happens that
\begin{equation} \label{Eq.Consistency.Clust_1}
\sup_{{\bf Z}_{1},  {\bf Z}_{2} \in {\mathcal C}^N_h} \left|D_{d_N}^{\Sigma_{d_N}}({\bf Z}_{1} ,{\bf Z}_{2})- L_S^h\right|  \stackrel{P}{\to}  0 \mbox{ as } N \to \infty.
\end{equation}
		
\item[b)] For any $h\neq k \in \{1,2\}$, 
let $\alpha_{d_N}=(\alpha^{d_N}_{1},\ldots,\alpha^{d_N}_{{d_N}})^T$ be the eigenvalues of $S_{d_N}^{hk}$ with ${d_N} \in \mathbb{N}$. If 
\begin{equation}\label{Hip.Theo.Conssit.Clusytering.Ass.b}
\log N =o\left(\frac{d_N}{\|\alpha_{d_N}\|_{\infty}}\right),
\end{equation}
then it happens that
\begin{equation} \label{Eq.Consistency.Clust_2}
\sup_{{\bf Z}_{1} \in {\mathcal C}^N_h, {\bf Z}_{2} \in {\mathcal C}^N_k} \left|D_{d_N}^{\Sigma_{d_N}}({\bf Z}_{1} ,{\bf Z}_{2}) - L^{hk}\right|  \stackrel{P}{\to} 0 \mbox{ as } N \to \infty.
\end{equation}
\end{itemize}
\end{Theo}


\begin{Nota} \label{Clust.R1}
{\rm
Assumption (\ref{A3.Thm1}) holds here, so $\displaystyle \frac{\|\alpha_{d_N}\|_{\infty}}{d_N} = \frac 1 {d_N} {\max_{1\leq i \leq d_N}\alpha^{d_N}_i}  \to 0$. Thus, if we take $d_N$ growing fast enough, then it is assured that assumptions (\ref{Hip.Theo.Conssit.Clusytering.Ass.a}) and (\ref{Hip.Theo.Conssit.Clusytering.Ass.b}) hold. The structure of the matrices $S^h_d$ and $S^{hk}_d$ for $d \in \Nat$ with $h \neq k \in \{1,2\}$ implies that a sufficient condition is  $\log N =o(d_N)$ ({see Proposition \ref{Cluster.C1} in Appendix I}).
}
\end{Nota}

\section{Transformations with Estimated Distributions} \label{Sec.Estimation}

In this section, we will discuss the first steps to implement the procedure described in Section \ref{Trans.GP}. 
In practice, the involved distributions and all the associated quantities need to be estimated from the data. Here, $\Zvec$ will denote indistinctly a random element with distribution $\Prob_1$ or $\Prob_2$ in the supervised problem, and the mixture $\pi_1 \Prob_1+\pi_2 \Prob_2$ in the unsupervised one.

For $j \in \mathbb{N}$, let $\phi_{j}^{\bf Z}(t)$ with $t \in [0,1]$ and $\lambda_{j}^{\bf Z}$ be the eigenfunctions and eigenvalues of \SigmaZ, respectively. We will now make the following assumptions:


\begin{itemize}
  
\item[\it A.1] $\sup_{t \in [0,1]} E[(\Zvec(t))^4] <\infty$. 
 
\item[\it A.2] It  happens that $\lambda_{1}^{\bf Z}> \lambda_{2}^{\bf Z}> \cdots >0$ satisfying $\sum_{j=1}^{\infty} \lambda_{j}^{\bf Z} < \infty$.

\end{itemize}
It is well known that Assumption A.2 implies $\{ \phi_{j}^{\bf Z}\}_{j \in \mathbb{N}}$ forms an orthonormal basis of $\Hil$.

To estimate $\Sigma^{\bf Z}$ and its eigenvalues and eigenfunctions, we will use the corresponding empirical quantities. Suppose that we have a simple random sample $\Zvec_1,\ldots,\Zvec_{N}$ taken from $\Prob_{\zvec}$. Given $s,t \in [0,1]$, we  define
\[
\hat \Sigma^{\bf Z}(s,t)= \frac 1 N \sum_{i=1}^N [\Zvec_i(s) - \overline{\Zvec}_N(s)][\Zvec_i(t) - \overline{\Zvec}_N(t)],
\]
where $\bar{\Zvec}_N(t) = \frac 1 N \sum_{i=1}^N \Zvec_{i}(t)$. Consider the corresponding estimated families $\hat \lambda_{1}^{\bf Z}\geq \hat \lambda_{2}^{\bf Z}\geq \cdots$ and $\hat \phi_{1}^{\bf Z},\hat \phi_{2}^{\bf Z},\ldots$ of its eigenvalues and eigenvectors, respectively. Note that $\hat \Sigma^{\bf Z}$ as well as all the $\hat \lambda_{j}^{\bf Z}$'s and $\hat \phi_{j}^{\bf Z}$'s depend on $N$. 
Given $\uvec \in \Hil$, we denote
\[
\hat u_j^{{\bf Z}} =  \langle \uvec, \hat \phi_j^{\bf Z} \rangle= \int_0^1 \uvec (t) \hat \phi_{j}^{\bf Z} (t) dt \mbox{ for } \ j \in \mathbb{N}.
\]

With a finite sample, we cannot estimate all the infinite eigenvalues and eigenvectors. 
Thus, we follow the work of \cite{Delaigle_Hall_2012} and \cite{Hall_2006}, and {select a non-random decreasing sequence $\eta_N$ going to zero slowly enough as to satisfy }$\lim_N N^{1/5}\eta_N= \infty$. We take
\begin{equation}\label{Eq.Rn}
\hat R^{\bf Z}_N = \inf \{ j: \hat \lambda_{j}^{\bf Z} - \hat \lambda_{j+1}^{\bf Z}< \eta_N \}-1.
\end{equation}
This definition implies that $\hat \lambda_j^{\bf Z} \geq \eta_N$ for every $j \leq \hat R_N^{\bf Z}$. 
Moreover, we will also need that the theoretical eigenvalues are reasonably well separated. To obtain this, given $\delta >0$, we also define
\begin{equation}\label{Eq.Rn_theor}
R^{\bf Z}_N = \inf \{ j:  \lambda_{j}^{\bf Z} -  \lambda_{j+1}^{\bf Z}< (1 + \delta) \eta_N \}-1,
\end{equation}
and later we will assume that $\log N =o(R^{\bf Z}_N)$. 
We will now state the empirical analogue of the results stated in Sections \ref{Class.GP} and \ref{Cluster.GP}.

\subsection{Consistency of Classification} \label{Class.Estimation}

Let ${\bf X}_1^1,\ldots,{\bf X}_n^1$ and ${\bf X}_1^2,\ldots,{\bf X}_m^2$ be sequences of independent observations taken from $\Prob_1$ and $\Prob_2$, respectively, which constitutes the training data. Denote $\bar{\Xvec}^1_n$ and $\bar {\bf X}_m^2$ to be the associated empirical means. Using the training data, we construct estimates of the underlying mean and covariance as described above.
Now, the estimated map $\hat D^{i}_{\hat R^{i}}({\uvec},\vvec)$ is defined using $D^{\hat \Sigma_{id}}_{d}({\uvec},\vvec)$ (recall the expression in (\ref{Eq.def.D_d})) with $d=\hat R_{n}^{1}$ or $\hat R_{m}^{2}$.

\begin{Theo} \label{Theo.Classif.2classes}
Let assumptions A.1 and A.2, and those in Theorem \ref{Thm2} hold. 
Let $\bf Z$ be an observation independent from both samples. 

\vspace{0.05in}
\noindent 
If $\ProbZ= \Prob_1$, then 

\vspace{0.05in}
{\hspace{0.25in} $\hat D^{1}_{\hat R_n^{1}}({\bf Z},\bar{\Xvec}^1_n) \stackrel{P}{\to} 1$ as $n \to \infty$ and $\hat D^{2}_{\hat R_m^{2}}({\bf Z},\bar {\bf X}_m^2) \stackrel{P}{\to} L^{12}$ as $m \to \infty$.}
\vspace{0.05in}

\vspace{0.05in}
\noindent 
If $\ProbZ= \Prob_2$, then

\vspace{0.05in}
{\hspace{0.25in} $\hat D^{1}_{\hat R_n^{1}}({\bf Z},\bar{\Xvec}^1_n) \stackrel{P}{\to} L^{21}$ as $n \to \infty$ and $\hat D^{2}_{\hat R_m^{2}}({\bf Z},\bar {\bf X}_m^2) \stackrel{P}{\to} 1$ as $m \to \infty$.}
\end{Theo}

\begin{Nota} \label{R1.Class.Emp}
{\rm
Let us assume the conditions of Theorem \ref{Thm2}. Now, if we follow the steps of Theorem \ref{Thm2}, then it is straight forward to prove that the misclassification probability of the estimated classifier goes to $0$ as $\min\{n,m\} \to \infty$, i.e., we obtain  perfect classification. 
}
\end{Nota}

\noindent 
Practical interest of Theorem \ref{Theo.Classif.2classes} is clear from the comments after Theorem \ref{Prop1}. Consequently, it happens that this result provides a procedure which asymptotically, as $\inf(n,m) \to \infty$, allows one to classify observations without a possibility of mistake. 
The modification to deal with $J (> 2)$ classes adds no special difficulty to this theoretical result. 

\subsection{Consistency of Clustering}

Let $\Zvec_1,\ldots,\Zvec_N$ be a simple random sample taken from $\Prob_{\bf Z}$ as described at the beginning of Section \ref{Cluster.GP}. Now, $\Prob_{\bf Z}$ and the sets $\mathcal{C}_1$ and $\mathcal{C}_2$ (containing information on the class labels) are unknown.
The extension of Theorems \ref{Prop2} and \ref{Theo.Conssit.Clusytering} to Theorems \ref{Theo.Cluster.2classes} and \ref{Theo.Cluster} is presented below. 
The following results will be based on the analysis of the map $\hat D_{\hat R_N}({\bf u},{\bf v})$, which is the function $D_{d}^{\Sigma}({\bf u},{\bf v})$ defined in equation (\ref{Eq.def.D_d}) with $d=\hat R_N$ (stated above in (\ref{Eq.Rn})), and the pooled covariance matrix $\Sigma_{\hat R_N}$ which is estimated by $\hat \Sigma_{\hat R_N}$ (sample covariance of the full sample). 
The first result is related to the consistency of the transformation on finite~sets.

\begin{Theo} \label{Theo.Cluster.2classes}
Let assumptions A.1 and A.2, and those in Theorem \ref{Theo.Cluster.M1} hold. \\

\noindent
(a) If $h=k \in \{1,2\}$, then
\begin{equation} \label{Cluster.GP.Est.E1}
\hat D_{\hat R_N} ({\bf Z}_1, {\bf Z}_2) \stackrel{P}{\to} L_S^{h} \mbox{ as } N \to \infty.
\end{equation}
(b) If $h \neq k \in \{1,2\}$, then 
\begin{equation} \label{Cluster.GP.Est.E2}
\hat D_{\hat R_N} ({\bf Z}_1, {\bf Z}_2) \stackrel{P}{\to} L_S^ {hk} \mbox{ as } N \to \infty.
\end{equation}
\end{Theo}

\noindent
In the context of clustering, we need an increasing sample size in order to estimate the parameters consistently. 
Thus, it is desirable to be able to cluster the increasing number of data points, asymptotically without error. The only way to achieve this is to get some kind of uniform convergence in \eqref{Cluster.GP.Est.E1} and  \eqref{Cluster.GP.Est.E2} when the sample size increases. 
This is the purpose of Theorem \ref{Theo.Cluster}, which gives us clear evidence that using this transformation would lead to asymptotic perfect separation in the empirical case as well.

\begin{Theo} \label{Theo.Cluster}
Let us assume all hypothesis in  Theorem \ref{Theo.Conssit.Clusytering} with $\log N = o(R_N^{\bf Z})$ in (\ref{Eq.Rn_theor}). 

\begin{itemize}

\item[(a)] For $h \in \{1,2\}$, it happens that 
\[
\sup_{{\bf Z}_{1},  {\bf Z}_{2} \in {\mathcal C}^N_h} \left|\hat D_{\hat R_N}({\bf Z}_{1} ,{\bf Z}_{2} ) - L_S^h\right|  \stackrel{P}{\to} 0 \mbox{ as } N \to \infty.
\]

\item[(b)]For any $h,k \in \{1,2\}$ with $h\neq k$, we have that
\[
\sup_{{\bf Z}_{1} \in {\mathcal C}^N_h, {\bf Z}_{2} \in {\mathcal C}^N_k} \left|\hat D_{\hat R_N}({\bf Z}_{1} ,{\bf Z}_{2} )- L_S^{hk} \right|  \stackrel{P}{\to} 0 \mbox{ as } N \to \infty.
\]
\end{itemize}

\end{Theo}


\begin{Nota} \label{R2.Clust.Emp}
{\rm
Clearly, Theorem \ref{Theo.Cluster.2classes} follows from Theorem \ref{Theo.Cluster}. But, the conditions required for proving the former are weaker, and hence we state it as a separate result.
}
\end{Nota}

{
\begin{Nota} [\it Asymptotic perfect identification of clusters]\label{R3.Clust.Emp}
{\rm
Recall the matrix $\Gamma_d$ from equation \eqref{Dist2.Cluster} with $d \in \mathbb{N}$. Now, consider the matrix $\hat \Gamma_{\hat R_N}$, which is obtained by replacing $\gamma_{ij}^d$s in the matrix $\Gamma_{\hat R_N}$ with their estimated values $\hat \gamma_{ij}^{\hat R_N}$ (computed using the quantities $\hat D_{\hat R_N}(\Zvec_i,\Zvec_j)$) with $1 \leq i \neq j \leq N$. Define $v_{12}=\pi_1\left|L_S^{1}-L_S^{12}\right|^2 + \pi_2\left|L_S^{2}- L_S^{21}\right|^2$. Fix $\epsilon >0$. Theorem \ref{Theo.Cluster} implies that with probability converging to one as $N\to \infty$, we have

\begin{itemize}

\item[-]
if $ \Zvec_i,\Zvec_j \in \mathcal{C}_h$ for $h \in \{1,2\}$, then
$
\left| \hat \gamma_{ij}^d \right| \leq 4\epsilon^2,
$

\item[-]

if $ \Zvec_i \in \mathcal{C}_h,\Zvec_j \in \mathcal{C}_k$ for $h\neq k \in \{1,2\}$, then 
$
\left| 
\hat \gamma_{ij}^d - v_{12}
\right|
\leq H\epsilon,
$
\end{itemize}
for a suitable $H>0$.
Consequently, if $v_{12} > 0$, then the elements in $\hat \Gamma_d$ will be clustered into two well-separated clusters: one around $0$ and another one around $v_{12}$ with probability converging to one.

Similarly, if $\Prob_{\bf Z}$ is a mixture of $J(>2)$ components and the $L_S^{h}$'s and $L_S^{hk}$'s satisfy 
$$v_{hk}:=\pi_h\left|L_S^{h}-L_S^{hk}\right|^2 + \pi_k\left|L_S^{k}- L_S^{kh}\right|^2$$ 
with $1 \leq h\neq k \leq J$. For positive and distinct $v_{hk}$s, the elements in the matrix $\hat \Gamma_d$ will be perfectly clustered into $1 + \left(
\begin{array}{c}
J\\2
\end{array}
\right)$ well-separated clusters: one of them around $0$ and the remaining around the values $v_{hk}$ (for $h<k$) with probability converging to one as $N\to \infty$.
Therefore, asymptotically, the sequence of matrices $\{\hat \Gamma_{\hat R_N}\}_{N \in \mathbb{N}}$ will contain enough information to perfectly cluster all the data points. 
}
\end{Nota}
}

\subsection{Implementation in Practice} \label{Implement}

In the supervised setting, we are given a sample of functional data points with known labels.
Implementation of our classification procedure is quite straight forward here because we only need to estimate the transformations using the sample, and apply any existing classification method on this transformed data. More specifically, we do not need to estimate the unknown constants $L^{12}$ and $L^{21}$ (stated in Theorem \ref{Theo.Classif.2classes}) for the implementation of our classification procedure. Given ${\bf Z}$, we compute the statistic ${\hat \Tvec}_{n,m}({\bf Z}) = (\hat D^{1}_{\hat R_n^{1}}({\bf Z},\bar{\Xvec}^1_n), \hat D^{2}_{\hat R_m^{2}}({\bf Z},\bar {\bf X}_m^2))^T$ (which is just the empirical version of the transformation ${\Tvec}_{d}$ based on the sample means and covariances) and apply any classification method on the $2$-dimensional data. 
The statistic ${\hat \Tvec}_{n,m}({\bf Z})$ involves $\hat R^{1}_n$ and $\hat R^{2}_m$. After estimating the eigenvalues of the covariances, we may use the respective expressions related with $\hat R^{1}_n$ and $\hat R^{2}_m$ (also see \eqref{Eq.Rn}). However, in practice, we prefer an approach which directly relates to the misclassification probability (see Section \ref{Class.Practice} for more details).

In the unsupervised setting, we are given a sample of functional data points without the labels. Here, we only need to consider the $N \times N$ estimated matrix ${\hat \Gamma}_N$ with the $(i,j)$-th element as $\hat D_{\hat R_N} ({\bf Z}_i, {\bf Z}_j)$ (which is just the empirical version of $D^{\Sigma_d^{\bf Z}}_{d} (\Zvec_i, \Zvec_j)$ based on the pooled sample covariance) for $1 \leq i,j \leq N$ and apply any clustering procedure on its rows (or, columns).
Again, note that we do not need to estimate the unknown constants $L_S^{h}$ and $L_S^{hk}$ for $h, k \in \{1,2\}$ (stated in Theorem \ref{Theo.Cluster.2classes}) for the implementation of our clustering procedure.
The expression related with $\hat R_N$ is also not used in practice (see Section \ref{Cluster.Practice} for more details).


\section{Analysis of Simulated Datasets} \label{Num.Res.Simulation}

For our simulation study, we consider two ($J=2$) class problems. 
We generated data on a discrete grid of $100$ equi-spaced points in the unit interval $[0,1]$ from four different simulation examples which are described below. Fix $s>0$.

\begin{itemize}
	\item[I.] Let $X_{h}(t) = \sum_{j=1}^{40} (\lambda_{hj}^{1/2} Z_{hj} + \mu_{hj}) \phi_j(t)$ with $t \in [0,1]$ and $h=1,2$. Here, the $Z_{hj}$s were independent standard normal (i.e., $N(0,1)$) random variables, $\phi_j(t) = \sqrt{2} \sin(\pi jt)$ with $t \in [0,1]$ and $j=1,\ldots,40$. Also, $\mu_{hj} = 0$ for $j > 6$, and we set the other components equal to $(0,-0.5, 1,-0.5, 1,-0.5)^T$ and $(0,-0.75, 0.75,-0.15, 1.4, 0.1)^T$ for $k = 1,2$, and $\lambda_{1j} = 1/j^{2}$ and $\lambda_{2j} = s/j^{2}$ for $j=1,\ldots,40$.\\
	{\tt This model is from the paper \cite{Delaigle_Hall_2012}.}

	\item[II.] 
	In this example, $X_1 \sim B$ and $X_2 \sim \mu+sB$ with $\mu(t)=Gt$ for $t \in [0,1]$ and $G \sim N(0,4)$ independent of $B$. Here, $B$ is the standard Brownian bridge, i.e., a centered Gaussian process with $\sigma_{ij}=\min(t_i,t_j)-t_it_j$ with $t_i, t_j \in [0,1]$ for $i,j \in \mathbb{N}$.
	
	Since $E[X_2(t)]=E[Gt]=0$ for $t \in [0,1]$, the {\it differences in mean never appear} in this setting. In fact, the inclusion of $\mu$ modifies the covariances because if $0<t_i < t_j<1$, then the independence between $G$ and the $B$ yields the following:
	\[
	E[X_2(t_i)X_2(t_j)]= 4t_it_j+s^2t_i(1-t_j). 
	\]
	
	{\tt This model is from the paper \cite{BCT_2018}.}
	
	
	\item[III.] Let $X_{h} = \mu_h + \sum_{j=1}^{50} \xi_{hj} {\lambda_{hj}}^{1/2} \phi_{j}$ for $h=1,2$. Here, $\xi_{hj}$s are i.i.d. $N(0,1)$, $\mu_{1}=0$ and $\mu_{2}(t) = t$ with $t \in [0,1]$, $\lambda_{1j} = e^{-j/3}$ and $\lambda_{2j} = \sqrt{s}e^{-j/3}$ for $j=1,\ldots,50$, and $\phi_{2i-1} = \sqrt{2} \sin (2i\pi t)$ and $\phi_{2i} = \sqrt{2} \cos (2i\pi t)$ for $i=1,\ldots,25$ with $t \in [0,1]$.\\
	{\tt This model is from the paper \cite{DMY_2017}.}
	
	\item[IV.] This two class problem consists of two Brownian motions defined in the closed interval $[0,1]$ with means $\mu_1(t)=
	20t^{1.1}(1 - t)$ and $\mu_2(t) = 20t(1 - t)^{1.1}$, respectively, for $t \in [0,1]$. For the first class, the eigenfunctions are $\phi_{j}(t) = \sqrt{2} \sin((j - 0.5) \pi t)$ and associated eigenvalues are $\lambda_{1j} = 1/(\pi (j - 0.5))^2$ for $j=1,\ldots,15$. The second class is similar to the first one, but the eigenvalues are multiplied by $\sqrt{s}$ (i.e., $\lambda_{2j} = \sqrt{s}\lambda_{1j} = \sqrt{s}/(\pi (j - 0.5))^2$) for $j=1,\ldots,15$.\\
	{\tt This model is from the paper \cite{Galeano}.}
\end{itemize}

We set $s=1$ for {\tt location only} problems. In {\tt location and scale} problems, we fixed $s=3$, while for {\tt scale only} problems the mean functions $\mu_1$ and $\mu_2$ were set to be the constant function $0$ and we retained $s=3$.

\subsection{Choice of $d$} \label{d_choice}

A critical issue is selection of the optimal dimension of the projected space for a given a set of data points (i.e., a fixed value of $(n,m)$, or $N$).
In the context of classification, we expect values of the estimated statistic ${\hat \Tvec}_{d}$ to form two clearly separated clusters depending on the class label of the observation for large values of $d$. Moreover, ${\hat \Tvec}_{d}$ should lie closer to the cluster formed by the transformed observations with the same class label as the test observation (also see Theorem \ref{Theo.Classif.2classes}). To demonstrate this, we construct a sequence of images and show how the separation varies with increasing values of $d$.
We generated samples of size $50$ from each of the two classes for the `scale case' of Example II, and used a pooled estimate of the scatter matrix to construct ${\hat \Tvec}_{d}$.

\begin{figure}[htp]
	\vspace{-0.1in}
	\centering
	\includegraphics[width=12cm, height=4cm]{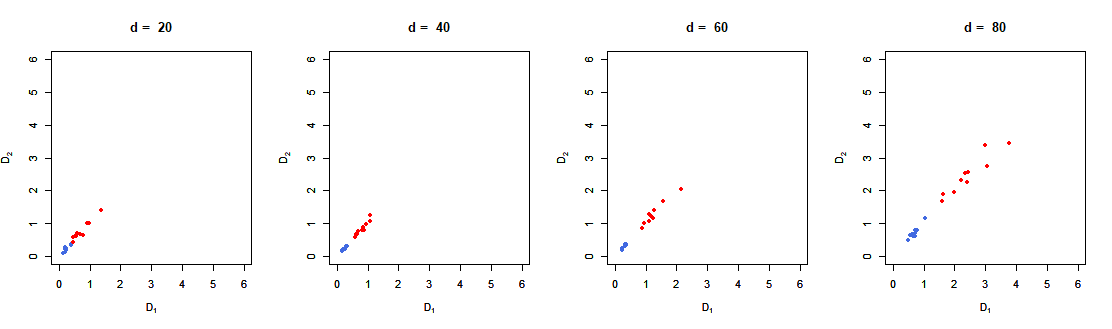}
	\vspace*{-0.1in}
	\caption{Transformed ($\hat \Tvec_d$) data points for varying values of $d$. The two colors correspond to observations from two different classes.}
	\label{Fig.1}
\end{figure}

Figure \ref{Fig.1} clearly shows that separation between the two classes increases with $d$. Observe that the transformed data points are concentrated about two distinct points corresponding to the two classes till $d=40$. The data clouds start to disperse when $d=60$, and this dispersion increases for $d=80$. This can be explained by the numerical instability in dimensions higher than $d=40$, with a fixed sample size of $N=100$.

On a related note, let us recall Theorem \ref{Theo.Cluster.2classes}. 
We observe a similar phenomena for the estimated matrix ${\hat \Gamma}_d$ in the clustering approach, when $250$ observations generated from each GD for the same example (as stated above). For the purpose of demonstration, the first $250$ observations correspond to the first GD, while the next $250$ observations to the second.
Note that this information is not a requirement for the implementation of our procedure.
Figure \ref{Fig.2} below shows the heatmap for increasing values of $d$, and we observe the best concentration at $d=80$. However, some noise in the off-diagonal submatrices for $d=80$ (compared to $d=60$) makes us to consider that the optimum could be somewhere between the values $60$ and $80$.

\begin{figure}[h]
\centering
\subfigure[$d=20$]{\includegraphics[width=0.2\linewidth,height=0.22\linewidth]{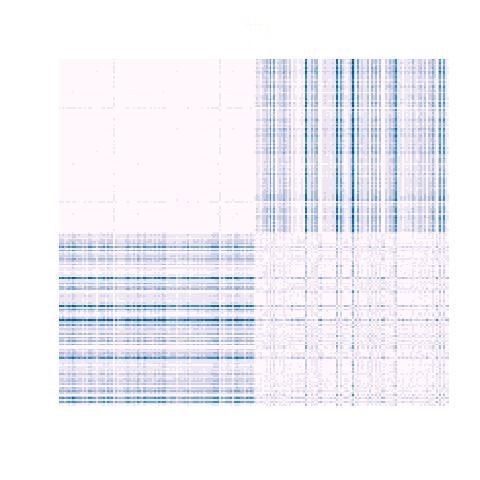}}
\subfigure[$d=40$]{\includegraphics[width=0.2\linewidth,height=0.22\linewidth]{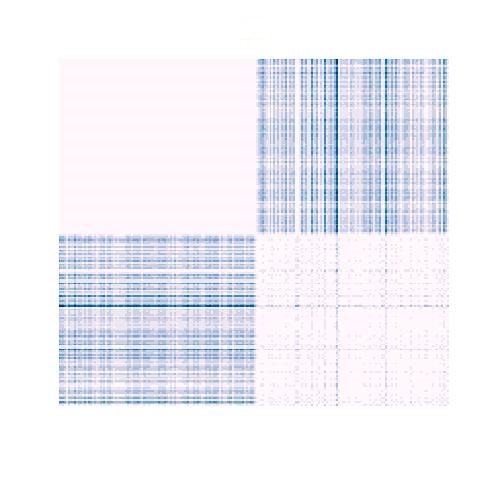}}
\subfigure[$d=60$]{\includegraphics[width=0.2\linewidth,height=0.22\linewidth]{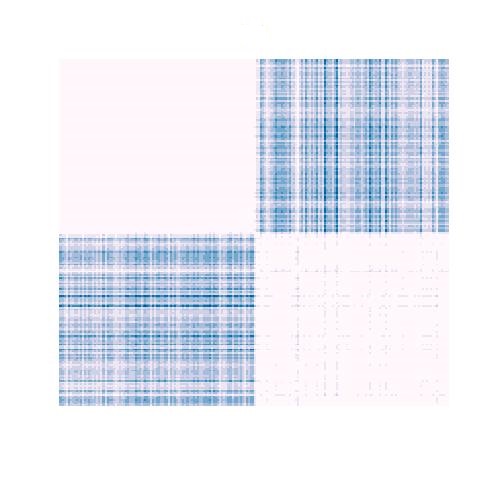}}
\subfigure[$d=80$]{\includegraphics[width=0.2\linewidth,height=0.22\linewidth]{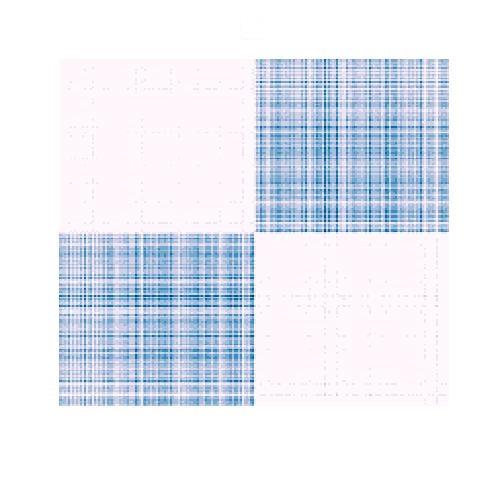}}
\vspace*{-0.1in}
\caption{Heatmap of ${\hat \Gamma}_d$ for varying values of $d$.}
\label{Fig.2}
\end{figure}

Clearly, the choice of $d$ is quite important as $d$ is the dimension of the space where we project our observations for a fixed sample size, and we can observe from Figures \ref{Fig.1} and \ref{Fig.2} that its estimation is quite crucial.
Subsections \ref{Class.Practice} and \ref{Cluster.Practice} contain further details on the choice of $d$ in practice for GP classification and clustering, respectively.

\subsection{Classification Procedure} \label{Class.Practice}

We propose to choose a unique dimension $d_{n,m}$ for both samples using cross-validation (CV). 
It is estimated by minimizing a CV estimate of the misclassification rate, and call it ${\hat d}_{CV}$ (defined below). For CV, we have used a common value of $d$ over both classes  in all situations to reduce the computational burden. To estimate the misclassification rate that a fixed value of $d$ produces, we randomly create a partition of the training sample into two subsets, $S_b$ and its complement for $b=1,\ldots,B$. 
The size of $S_b$ is $0.9$ times the size of the initial training sample, and the subsamples from each distribution are proportional to the original sample sizes (see \cite{HTF_2009} for more details). We take $B=50$. For each value of $b$, we treat the points in $S_b$ as the training set and its complement as the test set. With this split, we obtain $B$ estimates of the misclassification rate and compute their average to obtain $\hat{p}^{CV}_d$. Define $\hat{d}_{CV} = \arg \min_{2 \leq d \leq N} \hat{p}^{CV}_d$, where $N=n+m$ or $N=\min\{n,m\}$ if we use a pooled or individual estimates of the covariance, respectively.
The classifier ${\Psi}_d$ introduced in Section \ref{Class.GP} is not data adaptive. So, we use some data-driven classifiers which separate the two clusters (see Figure \ref{Fig.1}). We have included the {\it simple} centroid classifier (CD-CENT), support vector machine (CD-SVM) with a {\it linear} kernel, and the {\it non-linear} $1$ nearest neighbor (CD-$1$NN) classifier on the transformed data ${\hat \Tvec}_d$. One may refer to the book by \cite{HTF_2009} for more details on these three popular classifiers. 
Further, we considered two choices for estimating $\Sigma$, namely, the common pooled covariance or different covariances for each class (say, PC) and the common partial least squares (say, PLS) covariance of \cite{PCL_2007}.
{\tt R} codes for our methods are available here: \href{https://www.dropbox.com/sh/dug1n4ufxubqplr/AADxA1myR3K-krvZsAEh-KYwa?dl=0}{GP-classification}.
The minimum misclassification rate is reported as CD in the tables below. The complete results are available in Appendix II.

We considered several popular methods of functional classification for comparison with our proposed method. The first method is a centroid based classifier developed in \cite{Delaigle_Hall_2012}. The authors considered two variants, the first was based on principal components (DH-PC), while the second used PLS of \cite{PCL_2007} (DH-PLS). We report the best result among these two, and call it DH. The second set of methods are non-parametric approaches developed by \cite{Ferraty_2006}. The authors developed two classifiers, which are available as the functions {\tt funopare.knn.gcv} and {\tt funopadi.knn.lcv} from this link: {\tt http://www.math.univ-toulouse.fr/staph/npfda/}. 
The first classifier performs functional prediction of a scalar response from a sample of curves using the functional kernel estimator (NP1), 
while the second classifier performs functional discrimination of curves when a categorical response is observed (NP2). 
We report the best result among these two (with {\tt nknots} set to be $20$), and call it NP. 
Recently, some new methods have been proposed for functional classification. \cite{Galeano} propose several (about $20$) classifiers using the Mahalanobis distance and its variants, we report the minimum misclassification rate among them (GJL). A sequence of nonparametric one-dimensional density estimates are constructed, and the density ratio is used for classification by \cite{DMY_2017}. The best result is reported as DMY. \cite{BCT_2018} propose functional Fisher rules using reproducing kernels, which is denoted by BCT. Except DMY, all methods were implemented using the {\tt R} software. The codes for DH and BCT were kindly shared by the authors, while those for GJL was available from the journal webpage. The first author of DMY kindly shared the {\tt Matlab} codes.

The training sample size of each class is $50$, while the test is $100$. We replicated our experiment $200$ times, and the results are reported in tables below. The minimum misclassification rate is marked in {\bf bold}, while the second lowest is stated in {\it italics}. 

\begin{table}[!ht]
	\begin{center}
		\caption{Misclassification rates for different GPs with difference only in locations (with standard error in brackets).} \label{Table.4.1}
		\vspace{0.1in}
		
		\small
		\setlength{\tabcolsep}{10pt}
		\begin{tabular}{|c|c|c|c|c|c|c|} \hline
			Ex. & DH & NP & GJL & BCT & DMY & CD \\ \hline 
			
			I & 0.0007 & 0.0865 & {\bf 0.0001} & 0.0143 & 0.0006 & {\it 0.0003} \\ 
			& (0.0001) & (0.0021) & (0.0000) & (0.0006) & (0.0001) & (0.0001) \\ 
			
			{II} & 0.0706 & 0.1501 & {\bf 0.0000} & {\it 0.0185} & 0.0409 & 0.0343 \\ 
			& (0.0082) & (0.0046) & (0.0000) & (0.0043) & (0.0077) & (0.0063) \\ 
			
			III & 0.0561 & 0.1941 & {\it 0.0418} & 0.1842 & {\bf 0.0006} & 0.0589 \\ 
			& (0.0012) & (0.0042) & (0.0010) & (0.0023) & (0.0001) & (0.0013) \\ 
			
			IV & 0.0356 & 0.1846 & {\it 0.0080} & 0.0375 & {\bf 0.0001}  & 0.0351\\ 
			& (0.0009) & (0.0034) & (0.0004) & (0.0011) & (0.0000) & (0.0010) \\ 
			\hline
		\end{tabular}
	 \label{table_1_class}
	 \vspace{-0.2in}
	\end{center}
\end{table}

In the first setting, we considered classification problems with differences only in their {\tt location} parameters. It is clear from Table \ref{table_1_class} that GJL consistently yields the best result across all four examples with low misclassification rates. DMY is also competitive with the minimum misclassification rate in Examples III and IV. Our proposal yields good performance in Example I (which is mainly due to the use of PLS method), while BCT performs quite well in Example II. Overall, most of the methods lead to a stable performance in location problems, which supports the perfect classification property for GPs.

\begin{table}[!ht]
	\begin{center}
		\caption{Misclassification rates for different GPs with difference in locations and scales (with standard error in brackets).} \label{Table.4.2}
		\vspace{0.1in}
		
		\small
		\setlength{\tabcolsep}{10pt}
		\begin{tabular}{|c|c|c|c|c|c|c|} \hline
			Ex. & DH & NP & GJL & BCT & DMY & CD \\ \hline 
			
			I & 0.0256 & 0.0469 & 0.0122 & 0.0395 & {\bf 0.0030} & {\it 0.0042} \\ 
			& (0.0009) & (0.0014) & (0.0006) & (0.0012) & (0.0006) & (0.0003) \\ 
			
			{II} & 0.1110 & 0.0475 & {\bf 0.0000} & 0.0259 & 0.0064 & {\it 0.0046} \\ 
			& (0.0104) & (0.0014) & (0.0000) & (0.0039) & (0.0009) & (0.0007) \\ 
			
			III & 0.1576 & 0.0495 & 0.1141 & 0.1326 & {\bf 0.0038} & {\it 0.0103} \\ 
			& (0.0021) & (0.0014) & (0.0016) & (0.0017) & (0.0006) & (0.0006) \\ 
			
			IV & 0.1304 & 0.0510 & 0.0572 & 0.0596 & {\bf 0.0030} & {\it 0.0103} \\ 
			& (0.0020) & (0.0015) & (0.0011) & (0.0013) & (0.0008) & (0.0005) \\ 
			\hline
		\end{tabular}
	 \label{table_2_class}
 	 \vspace{-0.2in}
	\end{center}
\end{table}

In the second setting, we dealt with classification problems having differences in both {\tt location and scale} parameters. In Examples I, II and IV, DMY leads to the minimum misclassification rate, while GJL yields zero misclassification in Example II (see Table \ref{table_2_class}). Our proposed method is persistent, and holds the second position across all four examples. The performance of GJL and BCT deteriorates in Example III, while DH yields misclassification rates more than $10\%$ in Examples II, III and IV.

\begin{table}[!ht]
	\begin{center}
		\caption{Misclassification rates for several GPs with difference only in scales (with standard error in brackets).} 
		\vspace{.1in}
		
		\small
		\setlength{\tabcolsep}{10pt}
		\begin{tabular}{|c|c|c|c|c|c|c|} \hline
			Ex. & DH & NP & GJL & BCT & DMY & CD \\ \hline 
			
			I & 0.4564 & 0.0494 & 0.2138 & 0.1882 & {\it 0.0190} & {\bf 0.0114} \\ 
			& (0.0024) & (0.0014) & (0.0026) & (0.0021) & (0.0025) & (0.0005) \\ 
			
			II & 0.4884 & 0.0479 & 0.1295 & 0.1325 & {\bf 0.0037} & {\it 0.0185} \\ 
			& (0.0026) & (0.0015) & (0.0022) & (0.0023) & (0.0010) & (0.0006) \\ 
			
			III & 0.4542 & 0.0465 & 0.2994 & 0.2207 & {\it 0.0356} & {\bf 0.0118} \\ 
			& (0.0024) & (0.0014) & (0.0031) & (0.0020) & (0.0040) & (0.0006) \\ 
			
			IV & 0.4570 & 0.0518 & 0.1790 & 0.2064 & {\bf 0.0034} & {\it 0.0125} \\ 
			& (0.0024) & (0.0016) & (0.0022) & (0.0020) & (0.0006) & (0.0006) \\ 
			\hline
		\end{tabular}
	 \label{table_3_class}
	 \vspace{-0.2in}
	\end{center}
\end{table}

In the final setting, we considered classification problems with differences only in their {\tt scale} parameters. Our proposal based on the ${\hat \Tvec}_d$ transformation lead to the minimum misclassification rate in Examples I and III, while it attains the second best performance in Example II and IV. The DMY method performs quite well in this scenario, and interestingly, the ranks switch between these two methods in respective examples (see Table \ref{table_3_class}). The method NP consistently  holds the third position.
Performance of GJL and BCT is clearly not very promising in scale only problems. 
The DH method, which is a linear classifier, performs quite poorly in this scenario with a misclassification rate higher than 45\% across all examples.

Generally, we observe that the transformation CD combined with PLS works better than CD combined with PC when changes in location are involved. 
The opposite happens with scenarios involving differences in scales.
Among the three classifiers that we used on observations transformed using ${\hat \Tvec}_d$, both 1NN and linear SVM yield competitive misclassification rates (see the complete numerical results in Appendix II).

\subsection{Clustering Procedure} \label{Cluster.Practice}

Implementation of the clustering method is quite similar to the classification procedure described in the previous section. Again, one needs to choose the dimension $d$ suitably, and we use cross-validation (CV) here too. We use the idea developed by \cite{Wang_2010}, and state it briefly here.
Split the data into three random subsets (say, $S_{1b}$, $S_{2b}$ and $S_{3b}$) each of equal size for $b=1,\ldots,B$. For each value of $b$, we treat the points in $S_{1b}$ and $S_{2b}$ as training sets, and $S_{3b}$ as the validation set. For a fixed value of $d$ and given a clustering algorithm, the two training sets $S_{1b}$ and $S_{2b}$ are used to construct two cluster assignments. An appropriate distance between these two cluster assigments (say, $\mathbb{D}$) is computed based on the validation set $S_{3b}$ (see Section 2 of \cite{Wang_2010} for more details).
We repeat this partition $B=50$ times and average it over these $B$ samples to get $\hat{\mathbb{D}}^{CV}_d$. Define $\hat{d}_{CV} = \arg \min_{2 \leq d \leq N} \hat{\mathbb{D}}^{CV}_d$.

Recall the structure of the $\Gamma$ matrix in \eqref{Eq.Estr_Delta}, and also see Figure \ref{Fig.2}. As mentioned in Section \ref{Implement}, the number of clusters were estimated using the method described in Section \ref{Cluster.GP} 
(see equation \eqref{Eq.Cluster.Number}). To implement the procedure in practice, one needs to estimate the sequence $\{a_d\}_{d \in \mathbb{N}}$. However, we have used the function {\tt optishrink} available in the {\tt R} package {\tt denoiseR}.
This function extracts a {\it low-rank signal from Gaussian noisy data using the Optimal Shrinker of the singular~values}. 
The low rank structure of the $\Gamma$ matrix motivates us to directly apply this function on ${\hat \Gamma}_d$. The overall implementation yields desired results in our numerical study (see Tables \ref{table_1_clust} and \ref{table_2_clust} below).
We can apply any clustering method on the transformed data ${\hat \Gamma}_d$. 
In addition to the $k$-means algorithm (CD-$k$-means) discussed in Theorem \ref{Theo.Cluster.M1}, we considered spectral clustering (CD-Spectral) and Gaussian mixture models (CD-mclust). One may refer to the book by \cite{HTF_2009} for details on these three popular clustering methods. 
{\tt R} codes for our methods are available here: \href{https://www.dropbox.com/sh/ont3ggvz44g5j07/AAC1PRuzIWx9_yFUiR5mk4Zna?dl=0}{GP-clustering}.

We consider several methods for comparison. The first method is the classical $k$-means algorithm for functional data. Several competent methods for functional clustering using functional mixed mixture models are implemented in the function {\tt funcit} from the {\tt R} package {\tt funcy}. 
We report this method as funclust.
The methodology developed by \cite{CL_2007} is available in the function {\tt FClust} from the {\tt R} package {\tt fdaspace} 
using two clustering techniques `EMcluster' (CL1) and `kCFC' (CL2). We have reported the best result, and stated it as the CL method. 
In \cite{Delaigle_Hall_Pham_2019}, the authors developed functional clustering based on the $k$-means using basis functions. We implemented this method for two choices of the basis functions, namely, Haar and PC, and reported the best result among these two (we call it DHP). We have not used the DB2 basis for our comparisons because it requires the grid points to be of a power of $2$. The DHP method is available from the journal website, and we used those {\tt Matlab} codes for our comparisons.

Simulations are done based on models I to IV introduced earlier. We did not consider the {\tt location only} scenario as our proposed method is useless in such cases (recall (c) in Theorem \ref{Prop2}). However, we have some discussion and additional results in Appendix II 
for this scenario.
The sample size of each class was set to be $250$. Our experiment was replicated $100$ times, and the results are reported in Tables \ref{table_1_clust} and \ref{table_2_clust} below. We have computed adjusted Rand index using the function {\tt RRand} in the {\tt R} package {\tt phyclust}. One minus the adjusted Rand index are reported in tables below, where the minimum is marked in {\bf bold} and the second lowest is in {\it italics}. 

\begin{Nota}
{\rm
It is worth noting that all the competing methods require the number of clusters as an input variable, and we have run these methods with $k=2$ (the true number of clusters). However,  when applying the CD procedure we have estimated the number of clusters following the procedure described above. We obtained the correct value in more than 99\% of the cases (across all four examples for both scenarios) in our simulation study.
}
\vspace{-.1in}
\end{Nota}

\begin{table}[!ht]
	\begin{center}
		\caption{One minus adjusted Rand indices for different GPs with difference in location and scales (with standard error in brackets).} \label{Table.4.4}
		\vspace{.1in}
		
		\small
		\setlength{\tabcolsep}{10pt}
		\begin{tabular}{|c|c|c|c|c|c|} \hline
			Ex. & $k$-means & funclust & CL & DHP & CD \\ \hline
			
			I & 0.0632 & 0.1541 & {\it 0.0239} & 0.0818 & {\bf 0.0001} \\ 
			& (0.0007) & (0.0017) & (0.0007) & (0.0025) & (0.0001) \\ 
			
			II & 0.9445 & 0.8222 & {\it  0.5767} & 0.5149 & {\bf 0.4240} \\ 
			& (0.0036) & (0.0027) & (0.0045) & (0.0049) & (0.0030) \\ 
			
			III & 0.4250 & 0.3858 & {\it 0.2891} & 0.4137 & {\bf 0.0625} \\ 
			& (0.0017) & (0.0003) & (0.0000) & (0.0054) & (0.0006) \\ 
			
			IV & 0.4945 & 0.3975 & 0.1833 & {\it 0.1379} & {\bf 0.0000} \\ 
			& (0.0005) & (0.0011) & (0.0000) & (0.0033) & (0.0000) \\ 
			\hline
			
		\end{tabular}
	 \label{table_1_clust}
	 \vspace{-0.2in}
	\end{center}
\end{table}

In the first setting, we considered clustering problems with differences in their {\tt location and scale} parameters. Usefulness of the proposed transformation is clear from Table \ref{table_1_clust}. Our method attains the first position across all examples, while in Example IV we obtain perfect clustering. 
Although there is no location difference in Example II, sub-optimal performance of our method is probably due to low signal from the difference between the two covariance structures. CL attains the second best performance in the first three examples among the competing methods. DHP performs better than CL in Example IV.

In the next setting, we dealt with differences only in {\tt scale} parameters. It is clear from Table \ref{table_2_clust} that the separation in scatters is captured very well by the proposed transformation ${\hat \Gamma}_{d}$.
Moreover, our method again leads to perfect clustering (with a significant improvement in Example II compared to Table \ref{table_1_clust}). The method funclust (respectively, CL) gets second position in Examples II and III (respectively, Examples I and IV). The performances of $k$-means and DHP are similar, and quite bad in this scenario. Generally, the results in Table \ref{table_2_clust} suggest that existing methods fail to judiciously capture information if it is present only in the scale parameters.

\begin{table}[!ht]
	\begin{center}
		\caption{One minus adjusted Rand indices for different GPs with difference only in scales (with standard error in brackets).} 
		\vspace{.1in}
		
		\small
		\setlength{\tabcolsep}{10pt}
		\begin{tabular}{|c|c|c|c|c|c|} \hline
			Ex. & $k$-means & funclust & CL & DHP & CD \\ \hline
			
			I & 1.0019 & 0.9776 & {\it 0.8269} & 0.9966 & {\bf 0.0000} \\ 
			& (0.0000) & (0.0006) & (0.0000) & (0.0005) & (0.0000) \\ 
			
			II & 1.0006 & {\it 0.5004} & 0.9065 & 0.9999 & {\bf 0.0084} \\ 
			& (0.0001) & (0.0049) & (0.0007) & (0.0003) & (0.0003) \\ 
			
			III & 0.9990 & {\it 0.9956} & 0.9994 & 0.9967 & {\bf 0.0856} \\ 
			& (0.0002) & (0.0000) & (0.0000) & (0.0007) & (0.0006) \\ 
			
			IV & 0.968 & 1.0006 & {\it 0.8464} & 0.9980 & {\bf 0.0005} \\ 
			& (0.0001) & (0.0000) & (0.0000) & (0.0006) & (0.0004) \\ 
			\hline
			
		\end{tabular}
		\label{table_2_clust}
		\vspace{-0.2in}
	\end{center}
\end{table}

After applying the transformation ${\hat \Gamma}_{d}$, we used three methods for clustering the transformed observations. Overall, it seems that the Gaussian mixture model (i.e., mclust) achieves better results than the other two procedures (see the complete numerical results in Appendix II).

\section{Analysis of Benchmark Datasets} \label{Sec.BenchData}

We have applied our proposed methods to some benchmark data sets, {\tt Wheat} (from the {\tt R} package {\tt fds}), {\tt Control Chart} from the UCI Machine Learning Repository ({\tt https://arch ive.ics.uci.edu/ml/datasets.html}), {\tt Phoneme} and {\tt Satellite} (both datasets are available from the link: {\tt https://www.math.univ-toulouse.fr/$\sim$ferraty/SOFTWARES/NPFDA /index.html}), and {\tt Cars} (kindly provided by the authors of \cite{TorrecillaEtA_l2020}).

In classification problems, we have a sample of $\sum_j m_j$ observations for each set, where $m_j$ data points are from the $j$-th population for $1 \leq j \leq J$. Data were randomly split $M=100$ times to construct a training sample of size $n_j$, and a test sample of size $N_j$ from $j$-th population with $n_j+N_j=m_j$ for all $1 \leq j \leq J$. Define $n=\sum_j n_j$ and $N=\sum_j N_j$. We built the empirical classifiers using the training data, and used them to classify the test data points. For each classifier, the corresponding  proportion of test observations that were misclassified is reported below in Table \ref{Table.BenchData1} for the {\tt Wheat} and {\tt Phoneme} (with `aa' and `ao' class labels following \cite{BCT_2018}) data sets. 

\begin{table}[!ht]
	\begin{center}
		\caption{Misclassification rates for different classifiers (with standard error in brackets).} 
		
		\hspace*{-0.25in}
		\small
		\setlength{\tabcolsep}{10pt}
		\begin{tabular}{|c|c|c|c|c|c|c|c|c|} \hline
			Data & $(n,N)$ & $J$ & DH & NP & GJL & BCT & DMY & CD \\ \hline 
			
			Wheat & $(50,50)$ & $2$ & 0.1568 & 0.0899 & {\bf 0.0000} & 0.0440 & {\bf 0.0000} & {\bf 0.0000} \\ 
			& & & (0.0147) & (0.0026) & (0.0000) & (0.0082) & (0.0000) & (0.0000) \\ 
			
			Phoneme & $(100,1617)$ & $2$ & 0.2972 & 0.2200 & {0.2164} & 0.2387 & {\bf 0.2076} & {\it 0.2122} \\ 
			& & & (0.0063) & (0.0154) & (0.0032) & (0.0058) & (0.0045) & (0.0025) \\ 
			
			Control Chart & $(300,300)$ & $6$ & * & 0.2818 & * & * & * & {\bf 0.0145} \\ 
			& & & * & (0.0015) & * & * & * & (0.0005) \\ 
			
			\hline
		\end{tabular}
		\label{Table.BenchData1}
		\vspace{-0.2in}
	\end{center}
\end{table}

Table \ref{Table.BenchData1} shows that several methods like GJL, DMY and our proposal (CD) leads to {\it perfect classification} in the {\tt Wheat} data. The {\tt Phoneme} data poses a difficult scenario (see \cite{BCT_2018}), and the DMY method leads to the best performance closely followed by CD. 
The {\tt Control Chart} set comprises of $6$ classes. Most of the competing methods are for {\tt two} class problems, and one may solve $\binom{6}{2}$ $2$-class problems followed by majority voting. However, we did not opt for this option because our method can directly handle multi-class problems. In fact, our method CD exhibits {\it perfect classification} for this time series data set as well.

We ran a single execution of a data set (without splitting) for the {\tt Wheat}, {\tt Satellite} and {\tt Cars} data sets. Class assigments are already available for the {\tt Wheat} dataset. The {\tt Satellite} data has been analyzed in detail in the paper \cite{DFV_2007}, where the authors split the curves into two clusters `unimodal' and `multimodal'. 
The authors of this paper kindly shared the exact cluster assignments for this data set with us.
The {\tt Cars} data contains asset log-returns of the car companies Tesla, General Motors and BMW (see \cite{TorrecillaEtA_l2020} for more details). However, the rank of the estimated ${\hat \Gamma}_d$ matrix was $2$ for this data set, and our method detected only two distinct clusters. This is in accordance with \cite{TorrecillaEtA_l2020}, where the authors noted that assets of General Motors and BMW are very similar and quite difficult to distinguish. 
So, we merged General Motors with BMW while assigning the class labels for this data set. Consequently, the number of clusters was set to be $2$ for the competing methods.
We report one minus the adjusted Rand index for these data sets in Table \ref{Table.BenchData2}. Superiority of our proposal w.r.t. the competing methods is clear from the results given~below.

\begin{table}[!ht]
	\begin{center}
		\caption{One minus adjusted Rand indices for different clustering methods.} 
		\vspace{.05in}
		
		\small
		\setlength{\tabcolsep}{10pt}
		\begin{tabular}{|c|c|c|c|c|c|c|} \hline
			Data & $M$ & $2$-means & funclust & CL & DHB & CD \\ \hline 
			Wheat & $100$ & 0.6960 & 0.6960 & 0.8058 & {\it 0.5730} & {\bf 0.3644} \\
			Satellite & $472$ & 0.6072 & 0.6072 & {\it 0.6060} & 0.7253 & {\bf 0.4448} \\
			Cars & $90$ & {\it 0.8856} & {\it 0.8856} & 0.9650 & 0.9088 & {\bf 0.4680} \\
			\hline
		\end{tabular}
		\label{Table.BenchData2}
		\vspace{-0.1in}
	\end{center}
\end{table}

\noindent
To get a better understanding of the performance of our proposed method, we further computed the well-known average purity function. A value of average purity function close to one indicates good performance of a method.
We obtained the values as $0.90$, $0.8622$ and $0.8666$ for the {\tt Wheat} data, the {\tt Satellite} data and the {\tt Cars} data, respectively. Overall, our proposed method CD yields quite promising results in all three data sets.

\newpage

\appendix

\section*{Appendix I: Proofs and Mathematical Details} \label{Apendix}


\setcounter{page}{1}



\section{Proof of Theorem \ref{Thm1}}
Fix $d \in \mathbb{N}$. The $d$-dimensional random vector $({\bf Z} - b)_d $ has a Gaussian distribution with mean equal to $(\mu-b)_d$ and covariance matrix equal to $\Sigma_d$.
Now, $\|(A_d)^{-1/2} ({\bf Z} - b)_d\|^2$ is equal to  the square of the norm of a $d$-dimensional normal variable with mean ${\bf m}_d = (A_d)^{-1/2} (\mu-b)_d$ and covariance matrix $S_d = (A_d)^{-1/2} \Sigma_d (A_d)^{-1/2}$. Therefore, if 
$\uvec_d$ is a $d$-dimensional vector with centered normal distribution and covariance matrix equal to $S_d$, then
\begin{equation} \label{Eq.Des.Esc.prod}
D^{A_d}_d ({\bf Z} ,b)\sim
 \frac 1 {d} \langle {\bf m}_d + {\bf \uvec}_d, {\bf m}_d + {\bf \uvec}_d \rangle = \frac 1 {d}\left(\|{\bf m}_d\|^2 + \|{\bf \uvec }_d\|^2 + 2 \langle {\bf m}_d, {\bf \uvec}_d \rangle \right).
 \end{equation}
By assumption (\ref{A1.Thm1}), we have
\[
\lim_{d \to \infty} \frac 1 d \|{\bf m}_d\|^2 = L_\mu.
 \]

Let us consider the {\it second term} in (\ref{Eq.Des.Esc.prod}). Fix a basis in $V_d$ spanned by the eigenvectors of $S_d$. 
Note that this term is not dependent on $L_\mu$.
Denote $\uvec_d= (u_{d,1},\ldots,u_{d,d})^T$ and ${\bf m}_d= (m_{d,1},\ldots,m_{d,d})^T$ in this basis. Therefore, the random variables $(u_{d,i})^2$ with $1 \leq i \leq d$ are independent with means equal to $\alpha_i^d$ for $1 \leq i \leq d$ and $\sum_{i=1}^d (u_{d,i})^2 \sim\sum_{i=1}^d \alpha_i^d (u_{i})^2$. Here, $\{u_i\}_{1 \leq i \leq d}$ is a sequence of independent and identically distributed (i.i.d.) real variables with the standard normal distribution. We split the proof into two cases.

\subsubsection{$L_S$ is finite} 

Fix $\epsilon >0$. Taking into account that the variance of a $\chi^2$ distribution with one degree of freedom is two and using Tchebychev's inequality, we have that
\begin{eqnarray*}
\nonumber
\Prob \left[ \frac 1 d
\left|
\|{\bf \uvec }_d\|^2 - \mbox{trace}(S_d)
\right|
\geq \epsilon
\right]
&=&
\Prob \left[ \frac 1 d
\left|\sum_{i=1}^d 
\left(
(u_{d,i})^2 -\alpha_i^d
\right)
\right|
\geq \epsilon
\right]
\\
\nonumber
&\leq &
\frac{2 }{\epsilon^2 d^2} \sum_{i=1}^d (\alpha_i^d)^2 
\\
&\leq&
\frac{2}{\epsilon^2 d^2}  \| \alpha^d \|_{\infty} \sum_{i=1}^d \alpha_i^d , 
\end{eqnarray*}
which converges to zero by assumptions (\ref{A2.Thm1}) and (\ref{A3.Thm1}). Consequently, we have shown that 
\[
\frac 1 d \|{\bf \uvec }_d\|^2 - \frac 1 d \mbox{trace}(S_d) \stackrel{P}{\conv} 0\mbox{ as } d \conv \infty,
\]
and assumption (\ref{A2.Thm1})  gives 
\[
\frac 1 d \|{\bf \uvec }_d\|^2 \stackrel{P}{\conv} L_S\mbox{ as } d \conv \infty .
\]

\subsubsection{$L_S$ is infinite} 
We have that
\begin{eqnarray*}
\Prob \left[ \frac 1 {\sum_{i=1}^d \alpha_i^d}
\left|\sum_{i=1}^d 
\left(
(u_{d,i})^2 -\alpha_i^d
\right)
\right|
\geq \epsilon
\right]
&=&
\Prob \left[ 
\left|\sum_{i=1}^d \frac{\alpha_i^d} {\sum_{i=1}^d \alpha_i^d}
\left(
(u_{i})^2 - 1
\right)
\right|
\geq \epsilon
\right]
\\
&\leq &
\frac{2 }{\epsilon^2}
\sum_{i=1}^d \left(\frac{\alpha_i^d} {\sum_{i=1}^d \alpha_i^d}\right)^2\\
&\leq&
\frac{2}{\epsilon^2 }  \frac{\| \alpha^d \|_{\infty}}{ \sum_{i=1}^d \alpha_i^d }, 
\end{eqnarray*}
which converges to zero because $L_S=\infty$  and assumption (\ref{A3.Thm1}). Thus, we have shown that
\begin{equation} \label{Eq.equivalence}
\frac 1 {\frac1 d { \sum_{i=1}^d \alpha_i^d }}
\left(
\frac1 d \|{\bf \uvec }_d\|^2
-
\frac 1 {d} \mbox{ trace}(S_d)
\right)
\stackrel{P}{\conv} 0.
\end{equation}
Consequently, $\frac 1 d \|{\bf \uvec }_d\|^2$ converges to $\infty$ at the same rate as $\frac 1 {d} \mbox{ trace}(S_d)$.

Concerning  the {\it last term} in (\ref{Eq.Des.Esc.prod}), we have $\langle {\bf m}_d, {\bf \uvec}_d \rangle = \sum_{i=1}^d m_{d,i} u_{d,i}$. We split the proof into cases. 

\subsubsection{$L_\mu$ is finite} 
Fix $\epsilon >0$, and define $\alpha^d=(\alpha_1^d,\ldots,\alpha_d^d)^T$. Using Tchebychev's inequality again, we get 
\[ 
\Prob \left[
\frac 1 d |
\langle
{\bf m}_d,{\bf u}_d
\rangle
| > \epsilon
\right]
\leq
\frac{1}{\epsilon^2 d^2} \sum_{i=1}^d (m_{d,i})^2 \alpha_i^d
\leq
\frac{1}{\epsilon^2 d^2} \| \alpha^d \|_{\infty} \| {\bf m}_d \|^2,
\]
which converges to zero by assumptions (\ref{A1.Thm1}) and (\ref{A3.Thm1}), and the proposition is proved in this case. 

\subsubsection{$L_\mu$ is infinite} 
The result follows from equation (\ref{Eq.Des.Esc.prod}) and the previous results, if we are able to show that the sequence of real valued random variables
\[
w_d = \frac{\langle {\bf m}_d, {\bf \uvec}_d \rangle} { \max (\|{\bf m}_d\|^2 ,\|{\bf \uvec }_d\|^2)}
\]
converges to zero in probability as $d \to \infty$. In turn, this will be fixed if we show that every subsequence of $\{w_d\}$ contains a new subsequence which satisfies this property. Thus, let $\{w_{d_k}\}$ be a subsequence of $\{w_d\}$ and let us consider the associated subsequences $\{\|{\bf m}_{d_k}\|\}$  and $\{\|{\bf u}_{d_k}\|\}$. Obviously, there exists a further subsequence $\{d_{k^*}\}$ such that one of the following holds:
\begin{itemize}
\item[(i)]
$\lim_{d_{k^*}} \disp \frac{\|{\bf m}_{d_{k^*}}\|^2}{\mbox{trace}(S_{d_{k^*}})} = 0 $.

\item[(ii)]
$\lim_{d_{k^*}} \disp \frac{\|{\bf m}_{d_k}\|^2}{\mbox{trace}(S_{d_{k^*}})} = \infty $.

\item[(iii)]
There exists a finite $C>0$ such that $\lim_{d_{k^*}} \disp \frac{\|{\bf m}_{d_{k^*}}\|^2}{\mbox{trace}(S_{d_{k^*}})} = C $.

\end{itemize}

Notice that in cases i) and iii) $L_S=\infty$. 
To simplify notation, we denote the sequence $\{S_{d_{k^*}}\}$ by $\{S_h\}$, and similarly for the remaining ones. In case (i), since equation (\ref{Eq.equivalence}) shows that 
\begin{equation}\label{Eq.equivalence2}
\frac{\|{\bf \uvec }_h\|^2}
{\mbox{ trace}(S_h)}
\stackrel{P}{\conv} 1 \mbox{ as } h \to \infty,
\end{equation}
we have $\disp \frac{\|{\bf m}_h\|} {\|{\bf u}_h\|} \stackrel{P}{\conv} 0$ as $h \to \infty$. Consequently,
\[
\lim_h 
\left|w_h\right|  = \lim_h \frac{\left| \langle {\bf m}_h, {\bf \uvec}_h\rangle \right|} { \|{\bf \uvec }_h\|^2}
\leq
\lim_h 
\frac{ \|{\bf m}_h\|} { \|{\bf \uvec }_h\|} = 0 \mbox{ in probability.}
\]

If (ii) holds, we have that
$
\left| w_h \right|  
\leq
\frac{ \|{\bf \uvec }_h\|} { \|{\bf m}_h\|} 
$. Since $E\left[ \|{\bf \uvec }_h\|^2 \right] = \mbox{ trace}(S_d)$, we have that $\frac{ \|{\bf \uvec }_h\|^2} { \|{\bf m}_h\|^2} \stackrel{P}{\conv} 0$, and, also in this case $ w_h  \stackrel{P}{\conv} 0$ as $h \to \infty$.

In case (iii),  taking into account that equation (\ref{Eq.equivalence2}) now holds, it is enough to show that
\[
\frac{\langle {\bf m}_h, {\bf \uvec}_h \rangle} {C\mbox{ trace}(S_h)}
\stackrel{P}{\conv} 0 \mbox{ as } h \to \infty,
\]
Fix $\epsilon >0$. We have that 
\begin{eqnarray*}
\Prob \left[
\left|
\frac{\langle {\bf m}_h, {\bf \uvec}_h \rangle} {C\mbox{ trace}(S_h)}
\right|
> \epsilon
\right]
&\leq &
\frac 1 {C^2\epsilon^2} \sum_{i=1}^h \frac{m_{h,i}^2 \alpha_i^h}{\left(\sum_{i=1}^h \alpha_i^h\right)^2}
\leq
\frac 1 {C^2\epsilon^2} \frac {\|\alpha ^h\|_\infty }{ \sum_{i=1}^h \alpha_i^h} \frac{\|{\bf m }_h\|^2}
{\mbox{ trace}(S_h)},
\end{eqnarray*}
which converges to zero by assumptions (\ref{A2.Thm1}) and (\ref{A3.Thm1}). \FIN

\section{On assumptions (\ref{A2.Thm1}) and (\ref{A3.Thm1}) }

Next lemma shows that if $L_S< \infty$, then   assumption (\ref{A2.Thm1}) implies assumption (\ref{A3.Thm1}).
\begin{Lemm}\label{Lemm2}
Let $\{a_d\}_{d \geq 1}$ be a sequence of real positive numbers such that $\lim_d \frac 1 d \sum_{i=1}^d a_i$ exists, and it is finite. Then, it happens that 
$\lim_d \frac {1}d \|a^d\|_{\infty}=0.$

\end{Lemm}
\noindent
{\it Proof:} Fix $d \in \mathbb{N}$, and denote $A_d=\sum_{i=1}^d a_i$. We have that
\[
\frac {a_d}d = \frac {A_d}d -\frac {A_{d-1}}{d-1} \frac {d-1}d 
\]
and consequently, $0= \lim_d \frac {a_d}d$. Given $\epsilon >0$, there exists $d_1>0$ such that if $d>d_1$, then $\frac {a_d}d \leq \epsilon$ and $d_2\geq d_1$ such that
\[
\sup_{i=1,\ldots, d_1} \frac {a_i}{d_2} \leq \epsilon.
\]
Let $d>d_2$ and take $1 \leq i \leq d$. So, we have that if $i\leq d_1$, then $\frac{a_i} d < \frac{a_i} {d_2} \leq \epsilon$ and if $i>d_1$, then $\frac{a_i} d \leq \frac{a_i} i \leq \epsilon$. This completes the proof. \FIN

\section{Proof of Theorem \ref{Prop1}}

This proposition follows immediately from Theorem \ref{Thm1}. Statement (\ref{Eq.limit1}) follows by taking $b=\mu_{1}$, and $A_d = \Sigma_{1d}$. We obtain (\ref{Eq.limit2.2}) by taking $b=\mu_{2}$, and $A_d = \Sigma_{2d}$. \FIN

\section{Proof of Theorem \ref{Thm2}}

Let $\ProbZ=\Prob_{1}$. Under the assumptions of Theorem \ref{Prop1}, we get  
\begin{eqnarray*}
	\| \Tvec_{d}(\Zvec) - (1,L^{12})^T \|^2 &\stackrel{P}{\conv}& 0 \mbox{ and}
	\\
	\| \Tvec_{d}(\Zvec) - (L^{21},1)^T \|^2 &\stackrel{P}{\conv}& (1-L^{21})^2+(L^{12}-1)^2 ,
\end{eqnarray*}
as $ d \to \infty$.  
The latter expression equals $0$ iff $L^{21}=1$ and $L^{12}=1$. Thus, the assumption in the statement gives $I{[\Psi_d({\Zvec}) = 2]} \stackrel{P}{\conv} 0$ as $d \to \infty$. Here, $I[\cdot]$ is the indicator function. Applying the Dominated Convergence Theorem, we now get $\Prob_1[\Psi_d = 2] \to 0$ as $d \to \infty$.

When $\ProbZ= \Prob_2$, Remark \ref{R5} implies that the first limit is zero iff $L^{12}=1$ and $L^{21}=1$. A similar argument to the previous one gives that $\Prob_2[\Psi_d = 1] \to 0$ as $d \to \infty$.  \FIN

\section{Proof of Proposition \ref{Prop2}}

First, note that
$(\Zvec_1 - \Zvec_2)_d $ is a $d$-dimensional normal  vector, with mean $(\mu_{h}-\mu_{k})_d$ and covariance $\Sigma_{hd}+\Sigma_{kd}$.

\noindent
To prove (a) and (b) we will assume that (c) holds. Statement (c) is proved later, and its proof is independent of (a) and (b).

In case $(a)$, we have $h=k$. So, $(\mu_{k}-\mu_{h})_d= {0}_d$ and $\Sigma_{hd}+\Sigma_{kd} = 2\Sigma_{hd}$. If we take $A_d=S_d^h$, according to Remark \ref{R2}, (c) gives that assumption (\ref{A3.Thm1}) holds for this selection of $A_d$. Therefore, \eqref{Cluster.GP.E1} follows from Theorem \ref{Thm1}  because in this case $L^h_{\mu}=0$.

In case $(b)$, we have $h \neq k$. We take  $A_d=S_d^{hk}$ and $b=\mu_{h}-\mu_{k}$. Similarly as in (a) we have that assumption (\ref{A3.Thm1}) also holds in this case and Theorem  \ref{Thm1} gives that 
  \[
  D^{\Sigma^{\bf Z}_d}_{d} (\Zvec_1, \Zvec_2) \stackrel{P}{\to}L^ {hk}:= L_{\mu}^{hk}+L_S^ {hk} \mbox{ as } d \to \infty
\]
And \eqref{Cluster.GP.E2} follows because (c) gives that  $L_{\mu}^{hk}=0$.

To prove (c), let us denote $\Sigma^*= \pi_1 \Sigma_1 + \pi_2\Sigma_2$, $\mu=(\mu_1-\mu_2)$ and $\pi_{12} = \pi_1 \pi_2$, from (\ref{Pooled.Sig}), we have that
$$\Sigma_d =\Sigma^*_{d} + \pi_{12} \mu_{d}\mu_{d}^T.
$$
From here, the Sherman-Morrison formula gives
$$
\Sigma_d^{-1}=(\Sigma_{d}^*)^{-1} - 
\frac
{ \pi_{12} (\Sigma_{d}^*)^{-1} \mu_d \mu_d ^T (\Sigma_{d}^*)^{-1}}{1+\pi_{12} \mu_d^T (\Sigma_{d}^*)^{-1} \mu_d}. 
$$
Since $(\Sigma_{d}^*)^{-1}$ is positive definite for all $d \in \mathbb{N}$, this now implies that 
$$
0\leq \mu_d^T \Sigma_d^{-1}\mu_d = \mu_d^T (\Sigma_{d}^*)^{-1} \mu_d 
- \frac
{\pi_{12}(\mu_d^T (\Sigma_{d}^*)^{-1} \mu_d)^2}
{1+ \pi_{12}(\mu_d^T (\Sigma_{d}^*)^{-1} \mu_d)} = 
\frac
{\mu_d^T (\Sigma_{d}^*)^{-1} \mu_d}
{1+ \pi_{12}\mu_d^T (\Sigma_{d}^*)^{-1} \mu_d}
\leq \frac 1 {\pi_{12}},
$$
and the proof that $L_{\mu}^{hk} = 0$ trivially ends from definition of $L_\mu^{hk}$.

To handle the terms $L_S^h$ and $L_S^{hk}$, recall the Woodbury matrix identity: 
$$(U+V)^{-1}=U^{-1}-(U+UV^{-1}U)^{-1}.$$
Using this identity, we have 
\[
\Sigma_d^{-1}=1/\pi_1 \Sigma_{1d}^{-1} - B_d,
\]
where $B_d=(\pi_1 \Sigma_{1d}+\pi_1^2 \Sigma_{1d}(\pi_2\Sigma_{2d} + \pi_{12} \mu_{d} \mu_{d}^T)^{-1}\Sigma_{1d})^{-1}$.

If $U$ and $V$ are positive definite (p.d.), then $U^TVU$ is p.d.
In $B_d$, both the matrices $\Sigma_{1d}$ and $(\pi_2\Sigma_{2d} + \pi_{12} \mu_{d} \mu_{d}^T)$ are symmetric and p.d., and this implies that $B_d$ is also p.d. Further, $\Sigma_{1d}^{1/2}$ and $B_d$ are p.d. which now implies that $\Sigma_{1d}^{1/2} B_d \Sigma_{1d}^{1/2}$ is p.d. Recall that $trace$ is a linear map.
Now,
\begin{eqnarray*}
L_{S}^1 &=& \lim_d \frac 1 d trace(\Sigma_d^{-1/2} (2\Sigma_{1d}) \Sigma_d^{-1/2})\\
&=& \lim_d \frac 2 d trace(\Sigma_{1d} \Sigma_d^{-1})\\
&=& \lim_d \frac 2 d trace(1/\pi_1 I_d) - \lim_d \frac 2 d trace(\Sigma_{1d} B_d) \\
&=& \lim_d \frac 2 d trace(1/\pi_1 I_d) - \lim_d \frac 2 d trace(\Sigma_{1d}^{1/2} B_d \Sigma_{1d}^{1/2}) \\
& \leq & \lim_d \frac 2 d trace(1/\pi_1 I_d) = \frac{2}{\pi_1}.
\end{eqnarray*}
Similarly, we can also prove that $\displaystyle L_{S}^2<\frac{2}{\pi_2}$.
Again, 
\begin{eqnarray*}
L_{S}^{12} &=& \lim_d \frac 1 d trace(\Sigma_d^{-1/2} (\Sigma_{1d}+\Sigma_{2d}) \Sigma_d^{-1/2})\\
&=& \lim_d \frac 1 d trace(\Sigma_d^{-1/2}\Sigma_{1d} \Sigma_d^{-1/2}) + \lim_d \frac 1 d trace(\Sigma_d^{-1/2} \Sigma_{2d} \Sigma_d^{-1/2})\\
&<& \frac{1}{\pi_1} + \frac{1}{\pi_2} = \frac{1}{\pi_1\pi_2}.
\end{eqnarray*}
\FIN

\section{Proof of Proposition \ref{Prop.Number.Clusters}}

Under the conditions of Proposition \ref{Prop.Number.Clusters}, the number  of significant (unique) eigenvalues of the matrix $\Gamma$ is $2$. Recall that $N$ is fixed here. 

Consider the standardized distance matrix $D_d$ with the $(i,j)$-th element as $D^{\Sigma_d^{\bf Z}}_{d}({\bf z}_i,{\bf z}_j)$ for $1 \leq i, j \leq N$ and $d \in \mathbb{N}$.
We have a sequence of matrices $D_d \stackrel{P}{\to} D_0$ as $d \to \infty$ (componentwise). Since the map $D$ to $\Gamma$ is clearly continuous w.r.t. this convergence, we have that $\Gamma_d \stackrel{P}{\to} \Gamma$  as $d \to \infty$. Let us denote the eigenvalues of $\Gamma_d$ (respectively, $\Gamma$) to be $\beta_1^d,\ldots,\beta_N^d$ (respectively, $\beta_1,\ldots,\beta_N$). Since eigenvalues are continuous functions of the respective matrices, we have $\beta_j^d \stackrel{P}{\to} \beta_j$ as $d \to \infty$ for all $1 \leq j \leq N$. 

Let us now look into the following:
$$ \sum_{i=1}^N I(|\beta_i^d| > a_d) \stackrel{P}{\to} \sum_{i=1}^N I(|\beta_i^0| > 0) \mbox{ as } d \to \infty$$
with $a_d \downarrow 0$ as $d \to \infty$ at an appropriate rate. Recall that the limiting quantity on the right should give us the correct number of clusters. 
Consider the sequence $\{1/m\}_{m \in \mathbb{N}}$. Let us take $i$ such that $\beta_i^d \stackrel{P}{\to} 0$ as $d \to \infty$. Thus, for every $\epsilon,\delta>0$ there exists $D^i_{\delta,\epsilon}$ 
such that if $d \ge D^i_{\delta,\epsilon}$ then
\[
\Prob [ |\beta_i^d| > \delta ] < \epsilon.
\]
In particular, if we take $\delta=\epsilon = 1/m$, there exists $D_m^i$ such that if $d \geq D_m^i$:
\[
\Prob \left[ |\beta_i^d| > \frac 1 m \right] < \frac 1 m .
\]
Without loss of generality, we can assume that $D_1^i < D_2^i < \cdots$, and consider the sequence
\[
a_d^i = \left\{
\begin{array}{ll}
2 & \mbox{ if } 1 \leq i < D_1^i,
\\
[1mm]
\frac 1 m  & \mbox{ if } D_m^i \leq i < D_{m+1}^i, \mbox{ for some } m \geq 1.
\end{array}
\right.
\]
Then, obviously $a_d^i \to 0$, and 
\[
\Prob\left[ I(|\beta_i^d| > a_d^i) >0\right]
=
\Prob\left[ |\beta_i^d| > a_d^i \right] < a_d^i .
\]
If we define $a_d=\sup \{a_d^i: \beta_i^0=0\}$, and $i$ satisfies that $\beta_i^0=0$, then $I(|\beta_i^d| > a_d) \stackrel{P}{\to} 0$ as $d \to \infty$. 
A similar reasoning allows us also to conclude that if $|\beta_i^0| > 0$, then $I(|\beta_i^d| > a_d) \stackrel{P}{\to} 1$ as $d \to \infty$.
\FIN

\section{Proof of Theorem \ref{Theo.Cluster.M1}}



In this proof, we use the superindex $d$ in $\calG_i^d$ to emphasize that the groupings can change with the dimension $d \in \mathbb{N}$. Proposition \ref{Prop.Number.Clusters} implies that $K_d=2$ with probability converging to one. 

Note that $\phi(\calG_1,\ldots,\calG_J)$ has an alternative mathematical expression as 
\begin{equation} \label{clust.phi}
\sum_{h=1}^J \frac{1}{2|\calG_h|} \sum_{\small \uvec, \vvec \in \calG_h} \|\uvec-\vvec\|^2, 
\end{equation}
where $|\calG|$ denotes the cardinality of the set $\calG$.
Let us denote the rows/columns of $\Gamma_d$ as $\gammavec_{1}^d,\ldots,\gammavec_{N}^d$. The structure of $\Gamma^d$ implies that $\|\gammavec_{i}^d -\gammavec_{j}^d \|^2 \stackrel{P}{\to} 0$ as $d \to \infty$ iff $i,j \in \calC_h$ for $h \in \{1,2\}$. So, if each $\calG_h^d ~(h = 1,2)$ contains observations from the same population, then $\phi_d(\calG_1^d,\calG_2^d) \stackrel{P}{\to} 0$ as $d \to \infty$.

Let us assume that, on the contrary, there exists a subsequence of dimensions $\{d_k\}$ such that, for every $k$ there exists at least a couple of points $i_k,j_k$ such that $i_k \in \calG_{1}^d$ and  $j_k \in \calG_{2}^d$, say. Since the number of points is finite, there exists a further subsequence $\{d_{k^*}\}$ such that both sequences $\{i_{k^*}\}$ and $\{j_{k^*}\}$ are constant. Therefore, for those subsequences, \eqref{clust.phi} implies that
\[
\liminf_d \phi_d(\calG_1,\calG_2) \geq \lim_d \|\gammavec_{i_{k^*}}^{d_{k^*}}-\gammavec_{j_{k^*}}^{d_{k^*}}\|^2 \stackrel{P}{\to} \gamma_{12} >0.
\]
So, for the minimization of $\phi_d(\calG_1^d,\calG_2^d)$, each $\calG_h^d $ must contain all observations from a single population with probability converging to one as the dimension increases. This proves the convergence in probability of the Rand index $\mathbb{R}_{d,N}$ to zero as $d \to \infty$. 
\FIN

\section{Rank of the $\Gamma$ matrix}
{
Identifying number of clusters from the matrix $\Gamma$ is not equivalent to finding the rank of the matrix $\Gamma$.

\begin{Lemm} \label{J.clust}
The rank of the matrix $\Gamma$ is less or equal than $J$. However, the equality is only guaranteed when $J \leq 3$.
\end{Lemm}

\noindent
{\it Proof:} Trivially, $rank(\Gamma) \leq J$.
Let us denote the reduced Echelon form of $N \times N$ matrix $\Gamma$ as $\Gamma^{\circ}$. 
Thus, the matrix $\Gamma^{\circ}$ is a $J \times J$ symmetric matric with $\gamma_{ij}>0$ and distinct when $i \neq j$, while $\gamma_{ii}=0$. 


Moreover, for $J=3$, we have
\begin{equation*}
det(\Gamma^{\circ}) = det
\begin{pmatrix}
0 & \gamma_{12} & \gamma_{{13}}\\ 
\gamma_{12} & 0 & \gamma_{23}\\ 
\gamma_{13} & \gamma_{23} & 0 
\end{pmatrix}
=2 \gamma_{12} \gamma_{13} \gamma_{23} \neq 0.
\end{equation*}

In the case $J = 4$, a simple computation gives that if it happens that if $\disp \gamma_{12} = {\frac {\gamma_{{13}} \gamma_{{24}} + \gamma_{{14}} \gamma_{{23}} + 2\,\sqrt {\gamma_{{13}} \gamma_{{14}} \gamma_{{23}} \gamma_{{24}}}}{\gamma_{{34}}}}$, then $det(\Gamma^{\circ}) = 0$.
This happens, for instance, if we consider the following matrix (with all {\it positive and distinct} off-diagonal entries):
\begin{equation*}
\begin{pmatrix}
0 & t & 1 & 2\\ 
t & 0 & 3 & 4\\ 
1 & 3 & 0 & 5\\ 
2 & 4 & 5 & 0
\end{pmatrix},
\end{equation*}
where $t = 2 + 4 \sqrt{6}/5 > 0$.
}

\FIN

\section{Proof of Theorem \ref{Theo.Conssit.Clusytering}}

In order to simplify the writing, we will write $d$ instead of $d_N$.  We will use
the notation $\|\alpha_d\|_2:= \left( \sum_{i=1}^d (\alpha^d_i)^2\right)^{1/2}$. The real r.v.'s $\{u_i\}$ are assumed to be i.i.d. with standard normal distribution.

The following lemma is deduced  from Lemma 1 in \cite{LaurentM} on p. 1325, after some simple computations, taking into account that $\|\alpha_d\|_2 \geq \|\alpha_d\|_\infty$. We state it here for further reference.

\begin{Lemm} \label{Lem.Aux.Theo.Conssit.Clusytering}
If  $Z_d := \sum_{i\leq d} \alpha^d_i(u_i^2-1)$, and $x \geq 1$,  then
\[
\Prob\left[|Z_d| \geq  4 x \|\alpha_d\|_\infty\right] \leq 2\exp(-x).
\]
\end{Lemm}
We will also employ the following well known bound for the tail of the standard normal distribution:
\begin{equation}\label{Ineq.tail.normal}
\Prob [|N(0,1)| \geq t ] \leq \sqrt{\frac{ 2}{ \pi}} \exp(-t^2/2) \mbox{ for all } t \geq 1.
\end{equation}
\vspace{0.05in}

\noindent
{\it Proof of Theorem \ref{Theo.Conssit.Clusytering}} :
Let us show part b). The proof of (\ref{Eq.Consistency.Clust_1}) is similar to that of (\ref{Eq.Consistency.Clust_2}). We use the notation ${\bf m}_d= (\Sigma_d)^{-1/2}(\mu_1 - \mu_2)_d$ and ${\bf u}^{i}_d = (\Sigma_d)^{-1/2}({\bf Z}_i- \mu_i)_d$ with $d \in \mathbb{N}$, {
where $\Zvec_i$ is a generic observation with distribution $\Prob_i$ for $i=1,2$. Moreover, with an obvious abuse of notation, we will often write ${\bf u}_d^{i} \in {\cal C}_i^N$} with $d \in \mathbb{N}$ for $i=1,2$. 

Recall that $L_\mu=0$ and $L_S < \infty$ (see (c) in Proposition \ref{Prop2}). 
Repeating the first steps in the proof of Theorem \ref{Thm1}, we  have that 
\begin{eqnarray}
\nonumber
\lefteqn{
\sup_{{\bf Z}^1 \in {\mathcal C}^N_1, {\bf Z}^2 \in {\mathcal C}^N_2}  
\left| D_{d}^{\Sigma}({\bf Z}_{1}, {\bf Z}_{2} )- \frac 1 d \mbox{trace}(S_d^{12})\right|}
\\
\label{Eq2.3_1}
&\hspace*{12mm}\leq&
\left| 
\frac 1 d
\|{\bf m}_d\|^2
\right|
+
\sup_{{\bf u}^1 \in {\mathcal C}^N_1, {\bf u}^2 \in {\mathcal C}^N_2} 
 \left|\frac 1 d \|{\uvec}^1_d-{\uvec}_d^{2}\|^2- \frac 1 d \mbox{trace}(S_d^{12})\right|
\\
\label{Eq2.3_2}
&&+~
2 \sup_{{\bf u}^1 \in {\mathcal C}^N_1, {\bf u}^2 \in {\mathcal C}^N_2}  
\frac 1 d |\langle {\bf m}_d,{\uvec}^1_d-{\uvec}^2_d\rangle|,
 \end{eqnarray}
 and it is enough to prove that the terms in (\ref{Eq2.3_1}) and  (\ref{Eq2.3_2})  converge to zero in probability.

The {\tt first term} in (\ref{Eq2.3_1}) converges to zero by first part of (c) in Proposition \ref{Prop2}. Concerning the {\tt second one}, let $N_1, N_2$ be the number of elements in ${\cal C}_1^N$ and  ${\cal C} _2^N$, respectively. Since $N_1+N_2=N$, it is clear that $N_1\times N_2 \leq N^2/4$. 
Let $\varepsilon >0$. We have that
\begin{eqnarray}
\nonumber
P_N & := &
\Prob\left[
\sup_{{\bf u^1} \in {\mathcal C}^N_1, {\bf u^2} \in {\mathcal C}^N_2}  
\left|\frac 1 d \|{\bf u}^1_{d} -{\bf u}^2_{d} \|^2- \frac 1 d \mbox{trace}(S_d^{12})\right| > \varepsilon
\right]
\\
\nonumber
&=&
\Prob \left[\bigcup_{{\bf u}^1 \in {\mathcal C}^N_1, {\bf u}^2 \in {\mathcal C}^N_2} 
\left\{
 \left|\frac 1 d \|{\uvec}^{1}_d-{\uvec}^{2}_d \|^2- \frac 1 d \mbox{trace}(S_d^{12})\right|
 > \varepsilon
 \right\}
\right]
\\
\label{Eq.Cota.Uniforme}
&\leq&
\frac {N^2} 4\Prob \left[
 \left|\frac 1 d \|{\uvec}^{1}_d-{\uvec}^{2}_d \|^2-  \frac 1 d \mbox{trace}(S_d^{12})\right|
 > \varepsilon
\right],
\end{eqnarray}
where $\uvec^1$ and $\uvec^2$ are associated with some ${\bf Z}_1 \in {\mathcal C}^N_1$ and ${\bf Z}_2 \in {\mathcal C}^N_2$, respectively. However, it is clear that 
\[
\frac 1 d \|{\uvec}^{1}_d-{\uvec}^{2}_d \|^2-  \frac 1 d \mbox{trace}(S_d^{12}) \sim \frac 1 d \sum_{i=1}^d \alpha^d_i(u_i^2-1).
\]
Take $x=\varepsilon d/(4\|\alpha_d\|_\infty)$. By assumption (\ref{Hip.Theo.Conssit.Clusytering.Ass.b}), we have  $d/\|\alpha_d\|_\infty\to \infty$ and eventually $x\geq 1$. So, from Lemma \ref{Lem.Aux.Theo.Conssit.Clusytering}, we obtain
\[
P_N
\leq
\frac {N^2} 4 
\Prob \left[
\left|
\sum_{i\leq d} \alpha^d_i(u_i^2-1)
\right|
 >  \varepsilon d
\right]
\leq
\frac 1 2
\exp\left(- \frac {\varepsilon d}{4 \|\alpha^d\|_\infty} + 2 \log N\right),
\]
which converges to zero by assumption (\ref{Hip.Theo.Conssit.Clusytering.Ass.b}).

For the {\tt  term in equation }(\ref{Eq2.3_2}), we have that
\begin{eqnarray}
\nonumber
P_N^*&:=&
\Prob\left[
\sup_{{\bf u}^1 \in {\mathcal C}^N_1, {\bf u}^2 \in {\mathcal C}^N_2}  
\frac 1 {d} |\langle {\bf m}_d,{\uvec}^1_d-{\uvec}^2_d\rangle|
>\varepsilon
\right]
\\
\nonumber
&\leq&
\frac {N^2} 4 \Prob\left[ 
\frac 1 {d} |\langle {\bf m}_d,{\uvec}^1_d-{\uvec}^2_d\rangle|
>\varepsilon
\right]
\\
\nonumber
&=&
\frac {N^2} 4 \Prob\left[ 
\frac 1 {d} \left|\sum_{i\leq d}m_{di}(\alpha^d_i)^{1/2} u_i\right|
>\varepsilon
\right]
\\
\nonumber
&=&
\frac {N^2} 4 \Prob\left[ 
\left|N(0,1)\right|
>
\varepsilon\frac{d}
{\sqrt{\sum_{i\leq d}(m_{di})^2\alpha^d_i}}
\right]
\\
\label{Eq.Theo.Conssit.Clusytering.Ass.1}
&\leq&
\frac 1 {2^{3/2}\pi^{1/2}} 
\exp\left(-
\frac {\varepsilon^2} 2\frac{d^2}{\sum_{i\leq d}(m_{di})^2\alpha^d_i}+2\log N\right)
\\
\nonumber
&\leq&
\frac 1 {2^{3/2}\pi^{1/2}} 
\exp\left(-
\frac {\varepsilon^2} 2\frac{d^2}{\|\alpha^d\|_\infty\sum_{i\leq d}(m_{di})^2}+2\log N\right),
\end{eqnarray}
which converges to $0$ because of the fact that $L_\mu=0$ (see (c) in Proposition \ref{Prop2}) and (\ref{Hip.Theo.Conssit.Clusytering.Ass.b}). The same assumption allows us to apply inequality (\ref{Ineq.tail.normal}) to equation (\ref{Eq.Theo.Conssit.Clusytering.Ass.1}).
\FIN

\section{Result related to Remark \ref{Clust.R1}}

\begin{Prop} \label{Cluster.C1}
Under assumptions of Proposition \ref{Prop2}, if we assume that $\displaystyle \frac{\log N}{d_N} \to 0$, then conditions \eqref{Hip.Theo.Conssit.Clusytering.Ass.a} and \eqref{Hip.Theo.Conssit.Clusytering.Ass.b} hold.
\end{Prop}

\noindent
{\it Proof}:
Fix $h \in \{1,2\}$, and recall that 
$$\Sigma_d^{-1}=(\pi_h \Sigma_{hd})^{-1} - P_d,$$ 
where $P_d=(\pi_h \Sigma_{hd}+\pi_h \Sigma_{hd} (T_d^h)^{-1}\pi_h \Sigma_{hd})^{-1}$ is a positive definite matrix. Further, 
$$I_d + \Sigma_{d}P_d = \frac{1}{\pi_h} \Sigma_{d}\Sigma_{hd}^{-1}.$$
From here, Weyl's inequality gives
\begin{eqnarray} \label{E31}
1 \leq \alpha_{min}(\frac{1}{\pi_h} \Sigma_{d}\Sigma_{hd}^{-1})=\frac{1}{\pi_h} \alpha_{min}(\Sigma_{d}\Sigma_{hd}^{-1}).
\end{eqnarray}
Note the fact that the eigenvalues of the matrices $AB$ and $BA$ are same. So, the matrices $S^h_d$ and $\Sigma_{hd} \Sigma_{d}^{-1}$ will have the same eigenvalues. Furthermore, the eigenvalues of $S^h_d$ are the inverses of the eigenvalues of $\Sigma_{d} \Sigma_{hd}^{-1}$. Thus, (\ref{E31}) gives that

\begin{align} 
\alpha_{max}(S_{d}^h) < \frac{2}{\pi_h} (\mbox{free of } d). \label{E3}
\end{align}
\vspace{0.1in}
We now have
$$\log N =o\left(\frac{d_N}{\alpha^{d_N}_{1}}\right) \Leftrightarrow
\frac{\alpha^{d_N}_{1} \log N}{d_N} \to 0.$$
Equation (\ref{E3}) now implies that 
condition (\ref{Hip.Theo.Conssit.Clusytering.Ass.a}) holds if we assume $\displaystyle \frac{\log N}{d_N} \to 0$.
\vspace{0.1in}

\noindent
Fix $h \neq k \in \{1,2\}$. Our second matrix of interest is $$S^{hk}_d=(\Sigma_d)^{-1/2}(\Sigma_{hd}+\Sigma_{kd}) (\Sigma_d)^{-1/2}.$$
Since the matrices are symmetric, we have
$$\alpha_{max}(S_d^{hk}) \leq \alpha_{max}(\Sigma_d^{-1/2}\Sigma_{hd} \Sigma_d^{-1/2}) + \alpha_{max}(\Sigma_d^{-1/2}\Sigma_{kd} \Sigma_d^{-1/2}).$$
Again, the eigenvalues of $\Sigma_d^{-1/2}\Sigma_{id} \Sigma_d^{-1/2}$ and of $\Sigma_{id}\Sigma_d^{-1}$ will be equal for $i=h,k$. So,
\begin{eqnarray} \nonumber
\alpha_{max}(S_d^{hk}) &\leq& \alpha_{max}(\Sigma_{hd} \Sigma_d^{-1}) + \alpha_{max}(\Sigma_{kd} \Sigma_d^{-1})\\ \nonumber
&=& \frac{1}{\alpha_{min}(\Sigma_{d} \Sigma_{hd}^{-1})} + \frac{1}{\alpha_{min}(\Sigma_{d} \Sigma_{kd}^{-1})}
\\ \nonumber
& \leq & \frac{1}{\pi_h} + \frac{1}{\pi_k} = \frac{1}{\pi_h \pi_k} ~(\mbox{using equation } (\ref{E31})).
\end{eqnarray}
From here, similarly as before, we would obtain that $\displaystyle \frac{\log N}{d_N} \to 0$ implies (\ref{Hip.Theo.Conssit.Clusytering.Ass.b}) holds.
\FIN

\section{Proof of Theorem \ref{Theo.Classif.2classes}} \label{Proof.Theo.Classif.2classes}

Recall that in this theorem,  we use the subspaces generated by the estimates of the first $d$ eigenfunctions of the covariance of the random process $\bf Z$.

We begin with some notation and preliminary results which have been taken from \cite{Delaigle_Hall_2012} and \cite{Hall_2006}, or follow directly from the results there. Then, we will give the proof of Theorem \ref{Theo.Classif.2classes}. For every \nin, let us consider 
\begin{eqnarray}
\nonumber
\hat \Delta_{\bf Z}^2 & = & \int_0^1 \int_0^1 (\hat \Sigma(s,t) - \Sigma(s,t))^2 ds dt,
\\
\nonumber
\delta^{\bf Z}_{j} & = & \min_{k \leq j} (\lambda_{k} - \lambda_{k+1}).
\end{eqnarray}
In \cite{Delaigle_Hall_2012} and  \cite{Hall_2006}, it is shown that if $j\geq 1$, then
\begin{equation}
\label{Eq.Eigenv}
|\hat \lambda_{j} - \lambda_{j}|  \leq  \hat \Delta_{\bf Z},
\end{equation}
and that, if $j \leq \hat R^{\bf Z}_N$ (recall the definition of $\hat R^{\bf Z}_N$ in (\ref{Eq.Rn})), then
\begin{eqnarray}
\label{Eq.Eigenf}
\| \hat \phi_{j} -  \phi_{j}\|&\leq & 8^{1/2} \hat \Delta_{\bf Z} (\delta^{\bf Z}_{j})^{-1},
\\
\label{Eq.OrdenDelta}
\hat \Delta_{\bf Z} & = & O_p(N^{-1/2}),
\\
\label{Eq.OrdenRn}
R^{\bf Z}_N \to \infty \mbox{ and }  \hat R^{\bf Z}_N & \leq & \hat \lambda_1^{\bf Z} \eta_N^{-1}.
\end{eqnarray}

Moreover,  if $j \leq \hat R ^{\bf Z}_{N}$, there exists a $k \leq j$ such that
\begin{equation} \label{Eq.OrdenDeltan}
\delta^{\bf Z}_{j} =  \lambda_{k} -  \lambda_{k+1} \geq \hat \lambda_{k} - \hat \lambda_{k+1}-  2\hat \Delta_{\bf Z} \geq \eta_N -  2\hat \Delta_{\bf Z} = \eta_N + o_P(\eta_N),
\end{equation}
where we have applied (\ref{Eq.Eigenv}) and (\ref{Eq.Rn}) and that, from (\ref{Eq.OrdenDelta}) and the assumption on $\eta_N$, we can conclude that $\eta_N > 2 \hat \Delta_{\bf Z}$ from an index onward. Thus,  (\ref{Eq.OrdenDeltan})  and (\ref{Eq.Eigenf}) yield
\begin{equation}
\label{Eq.Eigenf2}
\| \hat \phi_{j} -  \phi_{j}\| \leq  8^{1/2} \frac{\hat \Delta_{\bf Z} }{\eta_N  - 2 \hat \Delta_{\bf Z}}.
\end{equation}
From (\ref{Eq.Eigenv}), (\ref{Eq.Rn}) and (\ref{Eq.OrdenDelta}), we obtain that
\begin{equation}
\label{Eq.CotaLambda}
 \lambda_{ j} \geq  \hat  \lambda_{j} - \hat \Delta_{\bf Z} \geq \eta_N  - \hat \Delta_{\bf Z} = \eta_N + o_P(\eta_N).
\end{equation} 
Now, we are in a position to prove {Theorem \ref{Theo.Classif.2classes}}.

\vspace{0.1in}
\noindent{\it Proof of Theorem \ref{Theo.Classif.2classes}:}
Let us assume that $\ProbZ=\Prob_1$. W.l.o.g., we assume that $\mu_{1}=0$.
We  split the proof into two lemmas. In Lemma \ref{Lemm.1Theo.Classif.2classes}, we  prove that $\left| \hat{D}^{1}_{\hat R_n^{1}}(\Zvec,\bar{\Xvec}_n^1) - {D}^{1}_{\hat R_n^{1}}(\Zvec,\bar{\Xvec}_n^1) \right| \stackrel{P}{\rightarrow} 0$ as $n \to \infty$. 
The proof that $\left| \hat{D}^{2}_{\hat R_m^{2}}(\Zvec,\bar{\Xvec}^2_m) - {D}^{2}_{\hat R_m^{2}}(\Zvec,\bar{\Xvec}^2_m) \right| \stackrel{P}{\rightarrow} 0$ as $m \to \infty$ is identical. 
Then, we will show in Lemma \ref{Lemm.2Theo.Classif.2classes} that the limits of ${D}^{1}_{\hat R_n^{1}}(\Zvec,\bar{\Xvec}^1_n)$ and  ${D}^{2}_{\hat R_m^{2}}(\Zvec,\bar{\Xvec}^2_{m})$  coincide with those of ${D}^{1}_{\hat R_n^{1}}(\Zvec,\mu_{1})$ and  ${D}^{2}_{\hat R_m^{2}}(\Zvec,\mu_{2})$, respectively.
Combining these two facts, the proof will be complete.

\begin{Lemm} \label{Lemm.1Theo.Classif.2classes}
Under the assumptions in Theorem \ref{Theo.Classif.2classes}, it happens that
\[
\left| \hat{D}^{1}_{\hat R_n^{1}}(\Zvec,\bar{\Xvec}^1_n)
-
{D}^{1}_{\hat R_n^{1}}(\Zvec,\bar{\Xvec}^1_n)
\right|
\stackrel{P}{\rightarrow}0 \mbox{ as } n \to \infty.
\]
\end{Lemm}
{\it Proof.}
For a fixed $\Zvec$, let us denote $\uvec = \Zvec - \bar{\Xvec}^1_n$. Note that $\uvec$ depends on $n$, but by the Strong Law of Large Numbers (SLLN), we have that $\|\uvec\| \leq \|\Zvec\| + \|\bar{\Xvec}^1_n\| = O(1)$ a.s. Let us denote $ (u_1,\ldots,u_{\hat{R}_n^1})^T$ and $(\hat u_1,\ldots,\hat u_{\hat{R}_n^1})^T$ to be the projections of ${\bf u}$ on the subspaces generated by the first ${\hat R_n^1}$ eigenvectors of the matrices $\Sigma_{\hat R_n^1}$ and $\hat \Sigma_{\hat R_n^1}$, respectively, when written in the basis generated by those eigenvectors.
Let \nin, and take $j \leq \hat R^{1}_{n}$. We now have
\begin{eqnarray*}
\left|
\frac{\left(  u_j\right)^2}{  \lambda^{1}_{j}}
 - 
 \frac{\left(\hat  u_j \right)^2}{\hat \lambda^{1}_{j}}
 \right|
 & = & 
 \left|
\frac{  u_j }{ (\lambda^{1}_{j})^{1/2}}
 - 
 \frac{\hat  u_j}
 { (\hat \lambda^{1}_{j})^{1/2}}
 \right|
\left|
\frac{  u_j }{  (\lambda^{1}_{j})^{1/2}}
 + 
 \frac{\hat  u_j }{ (\hat \lambda^{1}_{j})^{1/2}}
 \right|
 \\
 &\leq&
\left(
  \left|
\frac{  u_j -\hat  u_j}{ (\lambda^{1}_{j})^{1/2}}
 \right|
+
  \left|
 \hat  u_j
 \frac{  (\lambda^{1}_{j})^{1/2}-(\hat \lambda^{1}_{j})^{1/2}}{ (\lambda^{1}_{j}\hat \lambda^{1}_{j})^{1/2}}
  \right|
 \right)
\left|
\frac{  u_j }{  (\lambda^{1}_{j})^{1/2}}
 + 
 \frac{\hat  u_j }{ (\hat \lambda^{1}_{j})^{1/2}}
 \right|.
\end{eqnarray*}
We analyze each term in this expression separately  as follows: 
\begin{eqnarray}
\nonumber
\left|\frac{  u_j -\hat  u_j}{ (\lambda^{1}_{j})^{1/2}}\right|
& \leq & 
\frac{  1}{ (\lambda^{1}_{j})^{1/2}}
\int_0^1 |\uvec(t)| |\phi^{1}_{j}(t) -\hat  \phi^{1}_{j}(t) | dt
\\
\nonumber
& \leq &
\frac{  \| \uvec\| \ \|\phi^{1}_{j} -\hat  \phi^{1}_{j}\| }{ (\lambda^{1}_{j})^{1/2}}
\\ \nonumber
&\leq &
8^{1/2} \| \uvec\|  \frac{\hat \Delta_{\bf X}}{(\lambda^{1}_{j})^{1/2}(\eta_n - 2\hat \Delta_{\bf X})} \\ 
\label{Eq.Bound_1term}
& \leq & 8^{1/2} \| \uvec\| {\hat \Delta_{\bf X}}{(\eta_n^{-3/2}+o_P(\eta_n^{-3/2}))},
\end{eqnarray}
where we have applied the Cauchy-Schwartz inequality, (\ref{Eq.Eigenf2}), (\ref{Eq.OrdenDelta}) and (\ref{Eq.CotaLambda}).
On the other hand, we have
\begin{eqnarray}
\nonumber
\left|
\hat  u_j
 \frac{  (\lambda^{1}_{j})^{1/2}-(\hat \lambda^{1}_{j})^{1/2}}{ (\lambda^{1}_{j}\hat \lambda^{1}_{j})^{1/2}}
 \right|
 & \leq &
 \int_0^1  | \uvec (t)| | \hat \phi^{1}_{j} (t)| dt
\frac{| \lambda^{1}_{j}-\hat \lambda^{1}_{j}|}{\left((\lambda^{1}_{j})^{1/2} + (\hat \lambda^{1}_{j})^{1/2}\right)(\lambda^{1}_{j}\hat \lambda^{1}_{j})^{1/2}}  
\\
\nonumber
& \leq &
\| \uvec \| 
\frac{\hat \Delta_{\bf X}}{\left((\lambda^{1}_{j})^{1/2} + (\hat \lambda^{1}_{j})^{1/2}\right)(\lambda^{1}_{j}\hat \lambda^{1}_{j})^{1/2}} 
\\
\label{Eq.Bound_2term}
&\leq &
\frac 1 2\| \uvec \| {\hat \Delta_{\bf X}}{(\eta_n^{-3/2}+o_P(\eta_n^{-3/2}))},
\end{eqnarray}
where we have applied (\ref{Eq.Rn}) and  (\ref{Eq.CotaLambda}). Concerning the final term, using (\ref{Eq.CotaLambda}) and  (\ref{Eq.Rn}) again, we obtain that
\begin{eqnarray}
\label{Eq.Bound_3term}
\left|
\frac{  u_j }{  (\lambda^{1}_{j})^{1/2}}
 + 
 \frac{\hat  u_j }{ (\hat \lambda^{1}_{j})^{1/2}}
 \right|
 & \leq &
\|{\bf u}\|
\left( \frac{1 }{  (\lambda^{1}_{j})^{1/2}}
 + 
 \frac{1}{ (\hat \lambda^{1}_{j})^{1/2}}
 \right)
 \leq \| \uvec\| (\eta_n^{-1/2}+o_P(\eta_n^{-1/2})).
  \end{eqnarray}
  
Now, if we define $C= 8^{1/2}+1$, combining  (\ref{Eq.Bound_1term}), (\ref{Eq.Bound_2term}), (\ref{Eq.Bound_3term}), (\ref{Eq.OrdenRn}) and (\ref{Eq.OrdenDelta}), we get the following:
\begin{eqnarray*}
\nonumber
\left| \hat{D}^{1}_{\hat R_n^{1}}(\Zvec,\bar{\Xvec}^1_n)
-
{D}^{1}_{\hat R_n^{1}}(\Zvec,\bar{\Xvec}^1_n)
\right|
& \leq &
\frac 1 {\hat R^{1}_{n}}\sum_{j=1}^{\hat R^{1}_{n}}
\left|
\frac{\left(  u_j\right)^2}{  \lambda^{1}_{j}}
 - 
 \frac{\left(\hat  u_j \right)^2}{\hat \lambda^{1}_{j}}
 \right|
 \\
 &\leq&
  C  \|\uvec\|^2 {\hat \Delta_{\bf X}} (\eta_n^{-2}+o_P(\eta_n^{-2})) 
  = O_P(n^{-1/2}\eta_n^{-2}).
\end{eqnarray*}
By construction, $\eta_n$ is such that $n\eta_n^5 \to \infty$. So, we have $|\hat D^{1}_{\hat R_n^{1}}({\bf Z},\bar{\Xvec}^1_n) - D^{1}_{\hat R_n^{1}}({\bf Z}, \bar{\Xvec}^1_n)| \stackrel{P}{\rightarrow} 0$ as $n \to \infty$, and this lemma is proved.
\FIN

\begin{Lemm} \label{Lemm.2Theo.Classif.2classes}
Under the assumptions in Theorem \ref{Theo.Classif.2classes}, it happens that the limits in probability of ${D}^{1}_{\hat R_n^{1}}(\Zvec,\bar{\Xvec}^1_n)$ and  ${D}^{2}_{\hat R_m^{2}}(\Zvec,\bar{\Xvec}^2_{m})$  coincide with that of ${D}^{1}_{\hat R_n^{1}}(\Zvec,\mu_{1})$ and  ${D}^{2}_{\hat R_m^{2}}(\Zvec,\mu_{2})$, respectively.
\end{Lemm}

\noindent
{\it Proof.}
We will first show that $D^{1}_{\hat R_n^{1}}({\bf Z},\bar{\Xvec}^1_n) - D^{1}_{\hat R_n^{1}}({\bf Z},\mu_{1}) \stackrel{P}{\to} 0$ as $n \to \infty$.
To prove this, let us denote $(z_1,\ldots, z_{\hat R_n^{1}})^T$ and $(\bar {x}_1,\ldots,\bar {x}_{\hat R_n^{1}})^T$ to be the projections of ${\bf Z}$ and $\bar{\Xvec}^1_n$ on the subspace generated by the first ${\hat R_n^1}$ eigenvectors of the matrix $\Sigma_{\hat R_n^{1}}^{1}$, when written in the basis generated by those eigenvectors. Since we are assuming that $ \mu_{1}=0$, we have
\begin{eqnarray}
\label{Eq.Diference.1}
D^{1}_{\hat R_n^{1}}({\bf Z},\bar{\Xvec}^1_n) - D^{1}_{\hat R_n^{1}}({\bf z},\mu_{1}) 
& = &
\frac 1 {\hat R_n^{1}} \sum_{j=1}^{\hat R_n^{1}}  \frac {(\bar{x}_j)^2} { \lambda^{1}_{j}} 
- \frac 2 {\hat R_n^{1}}
\sum_{j=1}^{\hat R_n^{1}}  \frac {z_j \bar{x}_j} {\lambda^{1}_{j}} .
\end{eqnarray}

The r.v.'s $\{z_j/(\lambda^{1}_{j})^{1/2}\}$ are i.i.d. with the standard normal distribution because  $\ProbZ= \Prob_1$. Moreover, they are independent from the i.i.d. variables $\{\bar {x}_j (n/\lambda^{1}_{j})^{1/2}\}$ whose distribution is also standard normal. The SLLN implies that, for any sequence $\{T_n\} \subset \Nat$, with $T_n \to \infty$
\[
\frac 1 {T_n} \sum_{j=1}^{T_n} \frac { z_j\bar { x}_j } { \lambda^{1}_{j}} \stackrel{a.s.}{\to} 0.
\]

According to (\ref{Eq.OrdenRn}), $\hat R_n^{1}\convp \infty$, and a not too complicated reasoning leads to 
\[
 \frac 2 {\hat R_n^{1}}
\sum_{j=1}^{\hat R_n^{1}}  \frac { z_j\bar { x}_j } { \lambda^{1}_{j}} \convp 0.
\]

Thus, the second term in the right hand side of (\ref{Eq.Diference.1}) converges to zero in probability.
The reasoning to prove that the first term in the right hand side in (\ref{Eq.Diference.1}) converges to zero in probability is similar to the previous one, taking into account that the variables $\left\{ {n (\bar{x}_j)^2} /{ \lambda^{1}_{j}}, 1 \leq j \leq \hat R_n^{1} \right\}$ are i.i.d. with $\chi^2$ distribution with one degree of freedom.

We will now show that $D^{2}_{\hat R_m^{2}}({\bf Z},\bar {\bf X}_m^2) -  D^{2}_{\hat R_m^{2}}({\bf Z},\mu_{2})\convp 0$ as $m \to \infty$. In this part, we will change the notation. Let us denote $(z_1,\ldots,z_{\hat R_n^{2}})^T$ and $(\bar {x}_1,\ldots,\bar {x}_{\hat R_n^{2}})^T$ to be the projections of ${\bf Z}$ and $\bar{\Xvec}^2_n$ on the subspace generated by the first ${\hat R_n^2}$ eigenvectors of the matrix $\Sigma_{\hat R_n^{2}}^{2}$, when written in the basis generated by those eigenvectors. The proof is split into two cases.

\subsubsection{$L_S^{12}$ and  $L_\mu^{12}$ are finite}

Let us consider a non-random sequence $\{T_m\} \subset \Nat$, going to infinity with exact order $m^{2/6}$. For every $\epsilon >0$, we have that
\begin{eqnarray*}
\lefteqn{
\Prob \left[
\left|
D^{2}_{\hat R_m^{2}}({\bf Z},\bar {\bf X}_m^2) - D^{2}_{\hat R_m^{2}}({\bf Z},\mu_{2}) \right| >\epsilon \right]
}
\\
&\leq&
\Prob  [\hat R_m^{2} > T_M]
+\Prob \left[
\frac 1{\hat R_m^{2}}
\sum_{j=1}^{T_m}
\left|
 \frac {({z}_j-\bar {x}_j)^2- ({z}_j- \mu_{2j})^2 } { \lambda^{2}_{j}} 
\right| >\epsilon \right].
\end{eqnarray*}

From  (\ref{Eq.OrdenRn}), we have that the first term here converges to zero. Therefore, to finish this step, we only need to prove that  
$$ \frac 1{\hat R_m^{2}}\sum_{j=1}^{T_m}
\left|
 \frac {({z}_j-\bar {x}_j)^2- ({z}_j- \mu_{2j})^2} { \lambda^{2}_{j}} 
\right| \convp 0  \mbox{ as } m \to \infty.$$

However,
\begin{eqnarray}
\nonumber
\sum_{j=1}^{T_m}
\left|
 \frac {({z}_j-\bar {x}_j)^2- ({z}_j- \mu_{2j})^2 } { \lambda^{2}_{j}} 
\right|
&\leq&
\sum_{j=1}^{T_m}  \frac {\left|(\bar { x}_j)^2  -(\mu_{2j})^2 \right|} { \lambda^{2}_{j}} 
+
2 \sum_{j=1}^{T_m}  \frac { \left|z_j(\bar { x}_j -\mu_{2j})\right|} { \lambda^{2}_{j}} 
 \\
 \label{Eq.Suma_Tn}
 &=&
 \sum_{j=1}^{T_m}  \frac {\left|\bar { x}_j  -\mu_{2j} \right| \left|\bar { x}_j  +\mu_{2j}\right|} { \lambda^{2}_{j}} 
+
 2\sum_{j=1}^{T_m}  \frac { \left|z_j(\bar { x}_j -\mu_{2j})\right|} { \lambda^{2}_{j}} .
\end{eqnarray}

If we take expectations and apply the Cauchy-Schwartz inequality, we have that
\begin{eqnarray*}
\Esp
\left[\sum_{j=1}^{T_m}  \frac {\left|\bar { x}_j  -\mu_{2j} \right| \left|\bar { x}_j  +\mu_{2j} \right|} { \lambda^{2}_{j}} \right]
&\leq&
\sum_{j=1}^{T_m}  
\left(
\frac {\Esp \left[(\bar { x}_j  -\mu_{2j})^2 \right] \Esp \left[(\bar { x}_j +\mu_{2j})^2 \right|} { (\lambda^{2}_{j})^2}
\right)^{1/2}.
\end{eqnarray*}

It happens that $\bar {x}_j  -\mu_{2j}$ and $\bar { x}_j +\mu_{2j}$ are one-dimensional normal variables, with means equal to $0$ and $2\mu_{2j}$, respectively, and variances equal to $\lambda^{2}_{j}/m$ for $1 \leq j \leq \hat R_m^{2}$. Taking this into account, applying Jensen's inequality and the fact that $L_\mu^{12}< \infty$, we have
\begin{eqnarray*}
\Esp
\left[\sum_{j=1}^{T_m}  \frac {\left|\bar { x}_j  -\mu_{2j} \right| \left|\bar { x}_j  +\mu_{2j} \right|} { \lambda^{2}_{j}} \right]
&\leq&
\sum_{j=1}^{T_m} \left(
\frac 1 m 
\left(
\frac 1 m + \frac{(2 \mu_{2j})^2}{\lambda^{2}_{j}}
\right)
\right)^{1/2}
\\
&\leq &
\left( \frac {T_m} m \sum_{j=1}^{T_m} \left(
\frac 1 m + \frac{(2 \mu_{2j})^2}{\lambda^{2}_{j}} 
\right)
\right)^{1/2}
\\
&= &
\left( \frac {T_m} m \left(
\frac {T_m} m + 4 \left\|  (\Sigma_{2,T_m})^{-1/2}( \mu_{2} -  \mu_{1})_{T_m} \right\|^2
\right)
\right)^{1/2}
\\
& = & 
\left( \frac {T_m} m \right)^{1/2}
 \left(
2 \left\|  (\Sigma_{2,T_m})^{-1/2}( \mu_{2} -  \mu_{1})_{T_m} \right\| + o(1)
\right)
\\
&=&
O(T_m m^{-1/2})=O(m^{-1/6}).
\end{eqnarray*}

Now, we consider the expectation of the second term in (\ref{Eq.Suma_Tn}). Given $d \in \Nat$, let us denote $S_{d}= (\Sigma_{2d})^{-1/2}  \Sigma_{1d} (\Sigma_{2d})^{-1/2}$. Having in mind that ${\bf Z} \sim \Prob_1$ and that $\bar { x}_j$ is a  normal r.v. with mean $\mu_{2j}$ and variance equal to $\lambda^{2}_{j}/m$ for $1 \leq j \leq \hat R_m^{2}$, a similar reasoning gives
\begin{eqnarray*}
\Esp\left[
\sum_{j=1}^{T_m}  \frac { \left|z_j(\bar { x}_j -\mu_{2j})\right|} { \lambda^{2}_{j}} 
\right]
& \leq &
\sum_{j=1}^{T_m}  
\left(
\left[
\frac {\Esp[(z_j)^2] \Esp[(\bar { x}_j -\mu_{2j})^2]} 
{( \lambda^{2}_{j})^2} 
\right]
\right)^{1/2}
\\
& = &
\frac 1 {m^{1/2}}\sum_{j=1}^{T_m}  
\left(
\frac {\Esp[(z_j)^2] } 
{\lambda^{2}_{j}} 
\right)^{1/2}
\\
& \leq &
\left( \frac {T_m} {m}
\sum_{j=1}^{T_m}  
\frac {\Esp[(z_j)^2] } 
{\lambda^{2}_{j}} 
\right)^{1/2}
\\
& = &
\left( \frac {T_m} {m}
\mbox{Trace}(S_{T_m})
\right)^{1/2}
=
O(T_m m^{-1/2}) =O(m^{-1/6}).
\end{eqnarray*}

Therefore, both terms in equation (\ref{Eq.Suma_Tn}) are $O_P(m^{-1/6})$. Since $\hat R_m^{2} \stackrel{P}{\to} \infty$, we have that $\frac 1{\hat R_m^{2}}\sum_{i=1}^{T_m}
\left|
 \frac {({z}_j-\bar {x}_j)^2- ({z}_j- \mu_{2j})^2} { \lambda^{2}_{j}} 
\right| \stackrel{P}{\to} 0$ as $m \to \infty$. 

\subsubsection{$L_S^{12}$ or $L_\mu^{12}$ is infinite}
Our problem is to show that if ${\bf Z} \sim \Prob_1$, then 
\begin{equation}\label{Eq.lim_infinito}
\frac 1 {\hat R_m^{2}}
\sum_{j=1}^{\hat R_m^{2}}  \frac { (z_j- \overline{x}_j )^2} { \lambda^{2}_{j}} \stackrel{P}{\to} \infty \mbox{ as } m \to \infty.
\end{equation}

This case is very similar to the last part of the proof of Theorem \ref{Thm1}. Here, we have
\[
\sum_{j=1}^{\hat R_m^{2}}  \frac { (z_j- \overline{x}_j )^2} { \lambda^{2}_{j}}
=
\| 
(\Sigma_{\hat R_m^{2}}^{2})^{-1/2} (\Zvec -\overline{\Xvec}_m^2) 
\|^2.
\]

Thus, if we denote ${\bf m}_{2d} = (\Sigma_{2d})^{-1/2}\mu_{2d}$,  $\uvec_{2d} =  \Sigma_{2d}^{-1/2} \Zvec_d$ and $\bar{\Xvec}_m^* = (\Sigma_{2\hat R_m^2})^{-1/2} (\bar {\Xvec}_m^2 - \mu_{2})_d$, then 
\begin{equation} \label{Eq.Des.Esc.sum}
\sum_{j=1}^{\hat R_m^{2}}  \frac { (z_j- \overline{x}_j )^2} { \lambda^{2}_{j}}  
=
 \|{\bf m}_{\hat R_m^{2}}^{2}\|^2 
 + 
 \|{\bf \uvec }_{\hat R_m^{2}}^{2}\|^2 
 +
 \|\bar{\Xvec}_m^*\|^2  
 -
 2 \langle {\bf m}_{\hat R_m^{2}}^{2}, { \uvec }_{\hat R_m^{2}}^{2} \rangle
 +
 2 \langle \bar {\Xvec}_m^*,{\bf m}_{\hat R_m^{2}}^{2} - { \uvec }_{\hat R_m^{2}}^{2}\rangle.
 \end{equation}
 
By assumption (\ref{A1.Thm1}) we know that
\[
 \lim_d \frac 1 d \|{\bf m}_d^{2}\|^2 = L_\mu^{12}.
 \]
 
If we assume that $L_\mu^{12}=\infty$, given $M>0$, there exists $D>0$ such that $d^{-1} \|{\bf m}_d\|^2>M$ for every $d \geq D$. Thus, from (\ref{Eq.OrdenRn}), we have that
\[
\Prob \left[
\frac 1 {\hat R_m^{2}} \|{\bf m}_{\hat R_M^{2}}^{2} \|^2
\leq M
\right]
\leq
\Prob [ {\hat R_M^{2}} < D]
\conv 0.
\]

Concerning the second term in (\ref{Eq.Des.Esc.sum}), it happens that the random vector \Zvec \ is independent from the sequence $\{\hat R_m^{2}\}$. Thus, conditionally to this sequence, the distribution of the sequence $\{\|{\bf u}_{\hat R_m^{2}}^{2}  \|^2\}$ coincides with  that of a subsequence of $\{\|{\bf u}_{d}^{2}  \|^2\}$. However, along the proof of Theorem \ref{Thm1}, we proved that 
\[
\frac 1 d \|{\bf u}_{d} ^{2} \|^2 \stackrel{P}{\to} L_S^{12} \mbox{ as } d \to \infty.
\]

From here, a proof similar to that one we developed for the sequence $\{\|{\bf m}_{\hat R_m^{2}}^{2} \|^2\}$ allows us to conclude that if $L_S^{12}=\infty$, then
\[
\frac 1 {\hat R_m^{2}} \|{\bf u}_{\hat R_m^{2}} ^{2} \|^2 \stackrel{P}{\to} \infty, \mbox{ as } m \to \infty.
\]

The same reasoning we employed in Theorem \ref{Thm1} is enough to prove that
\[
\frac 1 {\hat R_m^{2}}
\frac{ \langle {\bf m}_{\hat R_m^{2}}^{2}, { \uvec }_{\hat R_m^{2}}^{2} \rangle}{ \max(\| {\bf m}_{\hat R_m^{2}}^{2}\|^2 \|, \|{ \uvec }_{\hat R_m^{2}}^{2} \|^2) } \stackrel{P}{\to} 0 \mbox{ as } m \to \infty.
\]

 The distribution of $ \|\bar{\Xvec}_m^*\|^2$ is equal to that of a sum of ${\hat R_m^{2}}$ squares of centered, one-dimensional normal variables with variance equal to $m^{-1}$. Thus, taking again $\{T_m\}\subset \Nat$ going to infinity at exact rate $m^{2/6}$, we would have that for every $\epsilon >0$,
\begin{equation} \label{Eq.ConvMedia}
\Prob [\|\bar{\Xvec}_m^*  \|^2 > \epsilon]\leq 
\Prob[\hat R_m^{2} > T_m]
+
\frac 1 \epsilon m^{-4/6},
\end{equation} 
and consequently, $\|\bar{\Xvec}_m^*  \| \stackrel{P}{\conv} 0$ as $m \to \infty$.

Only the last term in (\ref{Eq.Des.Esc.sum}) remains to be analyzed. Here, the Cauchy-Schwartz inequality and equation (\ref{Eq.ConvMedia}) allow us to conclude that
\[
\left | \langle \bar {\Xvec}_m^*,{\bf m}_{\hat R_m^{2}}^{2} + { \uvec }_{\hat R_m^{2}}^{2}\rangle \right|
\leq
\|\bar {\Xvec}_M^*\| \left( \| {\bf m}_{\hat R_m^{2}}^{2} \| + \|{ \uvec }_{\hat R_m^{2}}^{2} \|\right)
= o_P\left(\sup(\| {\bf m}_{\hat R_m^{2}}^{2} \| , \|{ \uvec }_{\hat R_m^{2}}^{2} \|)^2\right).
\]

This proves that the leading terms in  (\ref{Eq.Des.Esc.sum}) are the first two. Thus, the fact that at least one of the sequences $\{\| {\bf m}_{\hat R_m^{2}}^{2}\|\}$, or $\{ \|{ \uvec }_{\hat R_m^{2}}^{2} \|\}$ goes to infinity gives us a proof of (\ref{Eq.lim_infinito}). 
\FIN

\section{Proof of Theorem \ref{Theo.Cluster.2classes}}

According to the proof of Theorem \ref{Theo.Classif.2classes} (see Subsection \ref{Proof.Theo.Classif.2classes}), we need to check conditions (\ref{Eq.Eigenv}) to (\ref{Eq.OrdenRn}) in Subsection \ref{Proof.Theo.Classif.2classes}. 

Note that (\ref{Eq.Eigenv}), (\ref{Eq.Eigenf}) and (\ref{Eq.OrdenRn}) hold trivially. The only thing to be done is to prove that (\ref{Eq.OrdenDelta}) holds. 
Lemma 1.1 in Appendix II shows that $\Esp[\hat \Delta^2]=O(N^{-1})$. From here, Markov inequality gives us that $\hat \Delta^2 = O_P(N^{-1})$, and consequently, we have (\ref{Eq.OrdenDelta}). \FIN

\section{Proof of Theorem \ref{Theo.Cluster}}

We will need the following lemma:
\begin{Lemm} \label{Lema.Provisional} 
Under the assumptions in Theorem \ref{Theo.Cluster}, we have that $\Prob[\hat R_N \geq R_N] \to 1$.
\end{Lemm}
\noindent
{\it Proof} : Let $N\in \Nat$. From (\ref{Eq.Eigenv}), we have that
\[
\inf_{j \leq  R_N} (\hat \lambda_j - \hat \lambda_{j+1}) \geq \inf_{j \leq  R_N} (\lambda_{j} -  \lambda_{j+1})  - 2\hat \Delta_{\bf Z} \geq (1+\delta)\eta_N - 2\hat \Delta_{\bf Z},
\]
and the proof ends because  (\ref{Eq.OrdenDelta}) and the  fact that $\eta_N \geq N^{-1/5}$ imply that $\Prob[ \delta\eta_N - 2\hat \Delta_{\bf Z} \geq 0] \to 1$.
\FIN

In this setting, recall that $L_\mu=0$ and $L_S < \infty$ (see (c) in Proposition \ref{Prop2}). 
We will only  prove part b); part a) being similar. W.l.o.g. we will assume that $h=1$ and $k=2$. Remember that, for every ${\bf Z}_{1} ,{\bf Z}_{2} $, we have that
$$D_{\hat R_N}({\bf Z}_{1} ,{\bf Z}_{2} ) =
\frac 1 {\hat R_N} \sum_{j=1}^{\hat R_N} \frac{\langle {\bf Z}_{1} -{\bf Z}_{2},\phi_j\rangle^2}{\lambda_j}
\mbox{ and }
\hat D_{\hat R_N}({\bf Z}_{1} ,{\bf Z}_{2} ) =
\frac 1 {\hat R_N} \sum_{j=1}^{\hat R_N} \frac{\langle {\bf Z}_{1} -{\bf Z}_{2},\hat \phi_j\rangle^2}{\hat \lambda_j}.
$$

We are going to consider the function
\[
\tilde D_{\hat R_N}({\bf Z}_{1} ,{\bf Z}_{2})
=
\frac 1 {\hat R_N} \sum_{j=1}^{\hat R_N} \frac{\langle {\bf Z}_{1} -{\bf Z}_{2}, \phi_j\rangle^2}{\hat \lambda_j}.
\]

Obviously,
\begin{eqnarray}
\nonumber
\sup_{{\bf Z}_{1} \in {\mathcal C}^N_1, {\bf Z}_{2} \in {\mathcal C}^N_2} \left|\hat D_{\hat R_N}({\bf Z}_{1} ,{\bf Z}_{2} )- L_S^{12}\right|  
& \leq &
\sup_{{\bf Z}_{1} \in {\mathcal C}^N_1, {\bf Z}_{2} \in {\mathcal C}^N_2} \left|\hat D_{\hat R_N}({\bf Z}_{1} ,{\bf Z}_{2}) - \tilde D_{\hat R_N}({\bf Z}_{1} ,{\bf Z}_{2} )\right| 
\\
\nonumber
&&
+ \sup_{{\bf Z}_{1} \in {\mathcal C}^N_1, {\bf Z}_{2} \in {\mathcal C}^N_2} \left|\tilde D_{\hat R_N}({\bf Z}_{1} ,{\bf Z}_{2} )- D_{\hat R_N} ({\bf Z}_{1} ,{\bf Z}_{2} )\right| 
\\
\nonumber
&&
+ \sup_{{\bf Z}_{1} \in {\mathcal C}^N_1, {\bf Z}_{2} \in {\mathcal C}^N_2} \left| D_{\hat R_N}({\bf Z}_{1} ,{\bf Z}_{2} )-  L_S^{12}\right| 
\\
\nonumber
&=:&
T_1 + T_2 + T_3.
\end{eqnarray}

Lemma \ref{Lema.Provisional}, and equations (\ref{Eq.OrdenRn}) and (\ref{Eq.Eigenv}) imply that there exists $C>0$ such that
\[
\Prob[ R_N \leq \hat R_N \leq CN^{1/5}] \to 1.
\]

Consequently, with probability going to $1$, it happens that
\[
\label{Eq.Lema.4.6_4}
D_{ R_N}({\bf Z}_{1} ,{\bf Z}_{2} ) \leq D_{\hat R_N}({\bf Z}_{1} ,{\bf Z}_{2} ) \leq D_{CN^{1/5}}({\bf Z}_{1} ,{\bf Z}_{2} ).
\]

Since, by assumption (\ref{Hip.Theo.Conssit.Clusytering.Ass.b}),
$\log N =o\left(\frac{ R_N}{\lambda^{}_{1}}\right)$
and trivially we have $\log N =o\left(\frac{CN^{1/5}}{\lambda^{}_{1}}\right)$,  b) in Theorem \ref{Theo.Conssit.Clusytering} gives that $T_3$ converges in probability to zero as $N \to \infty$. Since $L_S^{12} < \infty$, this fact implies that
\begin{equation}
\label{Eq.Proof.3.3_2}
\sup_{{\bf Z}_{1} \in {\mathcal C}^N_1, {\bf Z}_{2} \in {\mathcal C}^N_2} D_{\hat R_N}({\bf Z}_{1} ,{\bf Z}_{2} )= O_P(1).
\end{equation}

With respect to $T_2$, we have that
\begin{eqnarray*}
\nonumber
T_2 &\leq &
\sup_{{\bf Z}_{1} \in {\mathcal C}^N_1, {\bf Z}_{2} \in {\mathcal C}^N_2} \frac 1 {\hat R_N} \sum_{j=1}^{\hat R_N} \frac{\langle {\bf Z}_{1} -{\bf Z}_{2}, \phi_j\rangle^2}{ \lambda_j}\frac{|\lambda_j-\hat \lambda_j|}{\hat \lambda_j}
\\
\label{Eq.Proof.3.3_3}
&\leq&
\sum_{j=1}^{ \hat R_N} \frac{|\lambda_j-\hat \lambda_j|}{\hat \lambda_j} \sup_{{\bf Z}_{1} \in {\mathcal C}^N_1, {\bf Z}_{2} \in {\mathcal C}^N_2} D_{\hat R_N}({\bf Z}_{1} ,{\bf Z}_{2} )
=O_p(N^{-1/10}),
\end{eqnarray*}
where last equality follows from (\ref{Eq.Proof.3.3_2}), (\ref{Eq.Eigenv}), (\ref{Eq.OrdenDelta}), (\ref{Eq.OrdenRn}) and (\ref{Eq.Rn}). 

Finally, given ${{\bf Z}_{1} \in {\mathcal C}^N_1, {\bf Z}_{2} \in {\mathcal C}^N_2} $, the Cauchy-Schwartz inequality and the fact that $\|\hat \phi_j\|=\| \phi_j\|= 1 $ imply
\begin{eqnarray*}
\nonumber
\left| \hat D_{\hat R_n}({\bf Z}_{1} ,{\bf Z}_{2} ) - \tilde D_{\hat R_n}({\bf Z}_{1} ,{\bf Z}_{2} ) \right|
& \leq & 
\frac 1 {\hat R_N} \sum_{j=1}^{\hat R_N} \frac{\left| \langle {\bf Z}_{1} -{\bf Z}_{2},\hat \phi_j\rangle^2 -  \langle {\bf Z}_{1} -{\bf Z}_{2},\phi_j\rangle^2 \right|}{\hat \lambda_j}
\\
\nonumber
& = & 
\frac 1 {\hat R_N} \sum_{j=1}^{\hat R_N} \frac{\left| \langle {\bf Z}_{1} -{\bf Z}_{2},\hat \phi_j- \phi_j\rangle\right| \left|  \langle {\bf Z}_{1} -{\bf Z}_{2},\hat \phi_j + \phi_j \rangle \right|}{\hat \lambda_j}
\\
\nonumber
& \leq & 
 \|  {\bf Z}_{1} -{\bf Z}_{2}\|^2 \frac 1 {\hat R_N} \sum_{j=1}^{\hat R_N} \frac{\|\hat \phi_j- \phi_j\| \   \| \hat \phi_j + \phi_j \|}{\hat \lambda_j}
\\
\nonumber
& \leq & 
2  \|  {\bf Z}_{1} -{\bf Z}_{2}\|^2 \frac 1 {\hat R_N} \sum_{j=1}^{\hat R_N} \frac{\|\hat \phi_j- \phi_j\|}{\hat \lambda_j}
\\
\label{Eq.Lema.4.6_1}
& = & 
2 \left\|  {\bf Z}_{1} -{\bf Z}_{2}\right\|^2 H_N.
\end{eqnarray*}

Moreover, the application of (\ref{Eq.Eigenf}), (\ref{Eq.OrdenDelta}), (\ref{Eq.OrdenDeltan}) and (\ref{Eq.Rn}) gives that $H_N = O_P(N^{-1/10})$, which in turn is equivalent to saying that there exists $C>0$ such that $\Prob[H_n < CN^{-1/10}] \to 1$. This and the reasoning leading to (\ref{Eq.Cota.Uniforme}) give that to prove that $T_1 \convp 0$ is enough to show that for every $C>0$
\begin{equation} \label{Eq.Proof.3.3_4}
{N^2} \Prob\left[
\left\|  {\bf Z}_{1} -{\bf Z}_{2}\right\|^2 > CN^{1/10} \right] \to 0 \mbox{ as } N \to \infty,
\end{equation}
where ${\bf Z}_{1}$ and ${\bf Z}_{2}$ came from distributions $\Prob_1$ and $\Prob_2$, respectively.

To show (\ref{Eq.Proof.3.3_4}), notice that ${\bf Z}_{1} -{\bf Z}_{2}$ follows a Gaussian distribution whose mean function is $\mu_1 - \mu_2$ and its covariance  is $\Sigma_{12} = \Sigma_1+\Sigma_2$. Let us denote $\gamma_j$ with $j \in \mathbb{N}$ the ordered eigenvalues of $\Sigma_{12}$. Let us consider a basis composed by eigenfunctions of $\Sigma_{12}$, we denote by $(\mu_1-\mu_2)_j$ the components of $\mu_1-\mu_2$ in this basis and $\{u_j\}$ is a sequence of i.i.d. real standard normal variables  with $j \in \mathbb{N}$. Now, we have the following
\begin{eqnarray}
\nonumber
 \left\|  {\bf Z}_{1} -{\bf Z}_{2}\right\|^2
 & \sim &
 \sum_{j=1}^\infty \left(
 \gamma_j^{1/2}u_j + (\mu_1 - \mu_2)_j
 \right)^2
 \\
\nonumber
 & = &
 \sum_{j=1}^\infty \left(
 \gamma_j(u_j^2-1) + \gamma_j+ (\mu_1 - \mu_2)_j^2 + 2 (\mu_1 - \mu_2)_j \gamma_j^{1/2}u_j
 \right)
 \\
 \nonumber
 & = &
 \sum_{j=1}^\infty \left(
 \gamma_j(u_j^2-1)  + 2(\mu_1 - \mu_2)_j \gamma_j^{1/2}u_j 
 \right)
 + \mbox{trace}(\Sigma_{12}) + \|\mu_1 - \mu_2\|^2.
\end{eqnarray} 

Notice that $K:=\mbox{trace}(\Sigma_{12}) + \|\mu_1 - \mu_2\|^2<\infty$. Thus, 
\begin{eqnarray}
 \nonumber
\Prob\left[
 \left\|  {\bf Z}_{1} -{\bf Z}_{2}\right\|^2  > CN^{1/10} \right]
&=&
\Prob\left[
\sum_{j=1}^\infty \left(
 \gamma_j(u_j^2-1)  + 2(\mu_1 - \mu_2)_j \gamma_j^{1/2}u_j 
 \right)> CN^{1/10} - K\right]
 \\
 \nonumber
 &\leq&
 \Prob\left[
\sum_{j=1}^\infty 
 \gamma_j(u_j^2-1)
 > \frac 1 2 \left(CN^{1/10} - K\right)\right]
 \\
 \nonumber
 &&
 +
 \Prob\left[
\sum_{j=1}^\infty 
(\mu_1 - \mu_2)_j \gamma_j^{1/2}u_j 
 > \frac 1 4 \left(CN^{1/10} - K\right)\right]
 \\
\label{Eq.Proof.3.3_7}
 & =: &
 P_1 + P_2
\end{eqnarray}

Obviously, $\frac 1 4 \left(CN^{1/10} - K\right) \to \infty$. Thus, eventually $\frac 1 4 \left(CN^{1/10} - K\right)>1$  and, from Lemma \ref{Lem.Aux.Theo.Conssit.Clusytering}, we have that
\begin{equation} \label{Eq.Proof.3.3_5}
P_1 
\leq
\lim_{d\to \infty} \Prob\left[
\sum_{j=1}^d 
 \gamma_j(u_j^2-1)
 > \frac 1 2 \left(CN^{1/10} - K\right)\right] \leq 2 \exp\left(- \frac 1 {8\gamma_1} \left(CN^{1/10} - K\right)\right).
\end{equation}

\noindent
Concerning to $P_2$, first notice that, for every $d \in \mathbb{N}$, the real r.v. $\sum_{j=1}^d 
(\mu_1 - \mu_2)_j \gamma_j^{1/2}u_j $ is centered normal, with variance equal to $\sum_{j=1}^d 
(\mu_1 - \mu_2)_j^2 \gamma_j^{} \leq \gamma_1 \sum_{j=1}^d 
(\mu_1 - \mu_2)_j^2 \leq \gamma_1 \|\mu_1 - \mu_2\|^2$. Therefore,
\begin{eqnarray}
\nonumber
P_2 & \leq &
\lim_{d\to \infty}
\Prob\left[ \left|
\sum_{j=1}^d
(\mu_1 - \mu_2)_j \gamma_j^{1/2}u_j \right|
 > \frac 1 4 \left(CN^{1/10} - K\right)\right]
 \\
 \nonumber
 &\leq&
 \Prob\left[ |N(0,1)| > \frac 1 {4\gamma_1^{1/2} \|\mu_1 - \mu_2\|} \left(CN^{1/10} - K\right)\right]
 \\
\label{Eq.Proof.3.3_6}
 &\leq &
 \sqrt{\frac{ 2}{ \pi}} \exp\left(-\frac 1 {2\gamma_1 (4\|\mu_1 - \mu_2\|)^2} \left(CN^{1/10} - K\right)^2\right),
\end{eqnarray}
where last inequality comes from (\ref{Ineq.tail.normal}) because, eventually $1 <  \left(CN^{1/10} - K\right)/({4\gamma_1^{1/2} \|\mu_1 - \mu_2\|})$.
Finally,  (\ref{Eq.Proof.3.3_7}), (\ref{Eq.Proof.3.3_5}), and (\ref{Eq.Proof.3.3_6}) give (\ref{Eq.Proof.3.3_4}), and, consequently, that $T_1 \convp 0$ as $N \to \infty$.
\FIN


\newpage
\setcounter{section}{0}

\section*{Appendix II: Supplementary} \label{Sect.Supplementary}

\section{Extension of the procedure to non-Gaussian distributions} \label{NonGaussianExtension}

{
Obviously, non-Gaussian processes can also be mutually singular. In fact, Theorem 4.3 in \cite{Rao_Varadarajan_1963} contains a sufficient condition for this property to be satisfied. This allows us to consider the possibility to extend previous results to cover non-Gaussian distributions.} 
It is obvious that the developed proofs can cover non-Gaussian distributions as long as they satisfy the due properties. 
In this subsection, we state the properties a distribution should satisfy in order the proofs can be extended. Thus, let $\Prob_1$ and $\Prob_2$ be two probabilities on the Hilbert space \Hil. Here, \Zvec \ will denote a $L_2[0,1]$-valued random element with distribution $\Prob_1$, $\Prob_2$ or $\pi_1\Prob_1+\pi_s \Prob_2$ for some $\pi_1,\pi_2>0$ with $\pi_1+\pi_2=1$.

The basic assumption is the existence of a covariance  of $\Zvec$. We will also consider assumptions {\it A.1} and {\it A.2} (see Section 3 in the main paper) and $b \in L_2[0,1]$. Given a positive-definite $d\times d$ matrix $A_d$ and a $d$-dimensional subspace $V_d \subset L_2[0,1]$, we need to consider the $d$-dimensional random vector $\Uvec_d = (A_d)^{-1/2}(\Zvec - b)_d$ and the covariance matrix $S_d=A_d^{-1/2}\Sigma_dA_d^{-1/2}$, where $\Sigma_d$ is the covariance matrix of $\Zvec_d$ and $(\Zvec - b)_d$ is the projection on $V_d$ of $(\Zvec - b)$ with $d \in \Nat$.

Let us write $\Uvec_d- \Esp[\Uvec_d] =(u^1,\ldots,u^d)^T$ in the basis of the eigenvectors of $S_d$ and let $\alpha_1^d,\ldots, \alpha_d^d$ be the eigenvalues of $S_d$. Therefore, $u_i/\alpha_i^d$ for $i=1,\ldots,d$ are real standardised random variables which we need to assume i.i.d. Similar properties must hold for the decomposition of \Zvec \ in its eigenfunctions basis (also see \cite{DMY_2017}).
We finally need two exponential inequalities as those stated in Lemma 6.3 and equation (24) of the main paper.

\section{Discussion on location case for clustering using $D^{\Sigma_d,r}_{d}$} \label{Clust.Location.Discussion}

We have some work in progress in order to fix the problem with the `location only' case. Recall the notation used in Subsection 2.2.2 of the main paper. As stated there, the problem in this case is that
\[
D^{\Sigma_d^{\bf Z}}_{d} ({\bf u}, {\bf v}) = \frac 1 d \left \|\Sigma_d^{-1/2}({\uvec}-\vvec)_d \right \|^2 = \frac 1 d \sum_{i=1}^d \frac{(u_i-v_i)^2}{\lambda_i}
\stackrel{P}{\to} 0  \mbox{ as } d \to \infty. 
\]

Our idea is to replace the terms in the sum with some others going to $0$ slowly (or, if possible, not converging to zero at all). To use this idea, our proposal is as follows:
\[
D^{\Sigma_d^{\bf Z},r}_{d} ({\bf u}, {\bf v}) := \frac 1 d \left \|(\Sigma_d^{-1/2})^r({\uvec}-\vvec)_d \right \|^2 = \frac 1 d \sum_{i=1}^d \frac{(u_i-v_i)^2}{\lambda_i^r}, \ \mbox{with } r \in \mathbb{I},
\]
where $\mathbb{I}$ is the set of integers.
In this article, we have studied the case when $r=1$, i.e., $D^{\Sigma^{\bf Z}_d,1}_{d}$. However, this was not a strict requirement and we look into some possible scenarios below. 

\begin{itemize}
\item If $r \in \{0,-1,-2,\ldots\}$, assumption $\it A.2$ in the main paper trivially gives that $\frac{(u_i-v_i)^2}{\lambda_i^r} \leq \frac{(u_i-v_i)^2}{\lambda_i}$ eventually for large $i$, and consequently, $D^{\Sigma^{\bf Z}_d,r}_{d} ({\bf u}, {\bf v}) \stackrel{P}{\to} 0$ as $d \to \infty$.

\item When $r \in \{2,3,\ldots\}$, the transformation $D^{\Sigma^{\bf Z}_d,r}_{d}$ may be useful because $1/\lambda_i^r$ will start to take high values (recall assumption $\it A.2$) and this may lead to separation between the observations of corresponding to different clusters.
\end{itemize}

Keeping the viewpoint stated above in mind, we slightly modify the transformation $D^{\Sigma^{\bf Z}_d,r}_{d} ({\bf u}, {\bf v})$ in our practical implementation. 
Numerical results for the transformation $D^{\Sigma^{\bf Z}_d,4}_{d}$, using the same settings as in Section 4 for the difference in {\tt location only} case are reported below. We have excluded Example II from this comparison because, as stated earlier, the difference of means in this case is always null.

\begin{table}[!ht]
	\begin{center}
		\caption{One minus adjusted Rand indices for different GPs with difference in location (with standard error in brackets).} \label{Table.Clust.Supplementary}
		\vspace{.1in}
		
		\small
		\setlength{\tabcolsep}{10pt}
		\begin{tabular}{|c|c|c|c|c|c|} \hline
			GP $\downarrow$ & $k$-means & funclust & CL & DHP & CD \\ \hline
			
			I & {\bf 0.0001} & 0.0002 & {\bf 0.0001} & {\bf 0.0001} & {0.0012} \\
			& (0.0001) & (0.0001) & (0.0001) & (0.0000) & (0.0002) \\[.1cm] 
			
			
			III & {\bf 0.0646} & {\it 0.0795} & 0.0945 & {0.1480} & {0.1649} \\
			& (0.0015) & (0.0009) & (0.0045) & (0.0047) & (0.0017) \\[.1cm] 
			
			IV & 0.1606 & {\it 0.0318} & {0.1015} & {\bf 0.0134} & {0.1257} \\
			& (0.0007) & (0.0003) & (0.0000) & (0.0004) & (0.0019) \\[.1cm] 
			\hline
			
		\end{tabular}
	\end{center}
\end{table}

The performance of $k$-means is quite good in Examples I and III. Both DHP and CL also perform quite well, securing a first place in some cases. 
The proposed statistic $D^{\Sigma^{\bf Z}_d,4}_{d}$ shows significant improvement (recall from Proposition 2.5 
that $L_{\mu}=0$ for $D^{\Sigma^{\bf Z}_d,1}_{d}$), and this is reflected in the numerical figures of Table \ref{Table.Clust.Supplementary}. 
Clearly, there is room for further work with the proposed transformation $D^{\Sigma^{\bf Z}_d,r}_{d}$ with $r \in \{2,3,\ldots\}$, both theoretically as well as numerically.

\section{Discussion of {
literature} proposing methods of `perfect classification and clustering'} \label{DellaigleHall}

The results in this paper are quite related to those in \cite{Delaigle_Hall_2012}, \cite{Delaigle_Hall_2013}, \cite{Delaigle_Hall_Pham_2019} and \cite{TorrecillaEtA_l2020}. In this section, we analyze the relation between the existing work and this paper. Except \cite{Delaigle_Hall_Pham_2019}, all these papers are devoted to obtain perfect classification, while \cite{Delaigle_Hall_Pham_2019} is devoted to perfect clustering.

\begin{itemize}

\item 
In \cite{Delaigle_Hall_2012}, the authors propose a centroid classifier using principal component (PC) and partial least squares (PLS) scores. They prove the perfect classification property for this classifier under homoscedasticidity and the Gaussian assumption (Theorem 1), and then generalize it for non-Gaussian distributions (Theorem 2). The authors also study the asymptotic properties of this classifier under the heteroscedastic scenario later in Theorem 4. 

The proposed classifier obtains `perfect classification' if the series $\|({\Sigma^d_{1}})^{-1/2}\mu^{2}_d\|^2  = \sum_{i=1}^d (\mu^{2}_{d,i})^2/\lambda_i^{1}$
diverges as $d \to \infty$. They also prove that perfect classification is impossible if this series converges (Theorem 1). 
However, if this series converges, then $(\mu^{2}_{d,i})^2/\lambda_i^{2} \to 0$. In this case, it is easy to show that $\nu = 0$ and our procedure is useless. 
However, let us assume that $\lambda_i^{1}=i^{-2}$, and $\mu^{2}_{d,i} = i^{-3/2}$ with $i \geq 1$. Then $\sum_{i \geq 1} (\mu^{2}_{d,i})^2 /\lambda_i^{1} = \sum_{i \geq 1} i^{-1} = \infty$, and the classifier in \cite{Delaigle_Hall_2012} is perfect. On the other hand, we have that $\frac 1 d \|({\Sigma^d_{1}})^{-1/2}\mu^{2}_d\|^2 = \frac 1 d \sum_{i \leq d} i^{-1} \to 0$. 
{\tt Therefore, if both distributions have a common covariance operator\\ (homoscedastic), then our procedure does not improve the classifier\\ proposed by \cite{Delaigle_Hall_2012}}. This is coherent with the fact that the classifier proposed by the authors is optimal under homoscedasticity.

\item In Section 2.2 of \cite{Delaigle_Hall_2013}, the authors propose a quadratic classifier, and study its theoretical properties in Theorem 1 for general distributions. Given a test observation $\bf Z$, \cite{Delaigle_Hall_2013} take $j_n$ going to infinity and analyze the sign of the difference
	\[
	\left(
	\left\|
	\left(
	\hat \Sigma_{j_n}^{\bf X}
	\right)^{-1/2}
	({\bf Z} - \overline{\bf X}_n)
	\right\|^2 
	+
	\log |  \hat \Sigma_{j_n}^{\bf X}|
	\right)
	-
	\left(
	\left\|
	\left(
	\hat \Sigma_{j_n}^{\bf Y}
	\right)^{-1/2}
	({\bf Z} - \overline{\bf Y}_m)
	\right\|^2
	+
	\log |  \hat \Sigma_{j_n}^{\bf Y}|
	\right),
	\]
	where $|A|$ stands for determinant of the matrix $A$. 

In our paper, we select the random sequences $\{R_n^{1}\}$ and $\{R_m^{2}\}$ (see (18) in the main paper 
for further details) and analyze the behavior of the two-dimensional statistic:
\[
\left(
\frac 1 {R_n^{1}} \left\|
\left(
\hat \Sigma_{R_n}^{1}
\right)^{-1/2}
({\bf z} - \overline{\bf x}_n)
\right\|^2
,
\frac 1 {R_m^{2}} \left\|
\left(
\hat \Sigma_{R_m}^{2}
\right)^{-1/2}
({\bf z} - \overline{\bf y}_m)
\right\|^2
\right)^T.
\]

Note that these two statistics are not equivalent. If $\alpha_i^{1}$s denote the ordered eigenvalues of $\Sigma_{R_n}^{1}$, then the Abel series summation criteria gives
$$
\frac 1 {R_n^{1}} \log |  \Sigma_{R_n}^{1}| = \frac 1 {R_n^{1}} \sum_{i=1}^{R_n^{1}} \log(\alpha_i^{1})  \to - \infty, \mbox{ as } n \to \infty.
$$

Our procedure has more power (this comes from the factors $1/R_n^{1}$ and $1/R_m^{2}$, which imply that we take means while \cite{Delaigle_Hall_2013} only take sums). Specifically, {\tt we can handle cases in which {
assumptions (c)  and (d)} of (B.6) in the Supplemental of \cite{Delaigle_Hall_2013} (related with bounded eigenvalues) do not hold}. 


\item
The paper by \cite{TorrecillaEtA_l2020} aims to analyze the implications of HFp. As a consequence, the authors provide a procedure to determine if two GPs are mutually singular, or not. In the case of mutually singular processes, the paper includes a procedure giving asymptotically perfect classification. The results cover both homoscedastic and heteroscedastic situations. 

Theoretical comparison with the results in \cite{TorrecillaEtA_l2020} is difficult because those authors employ an approach based on the properties of the reproducing kernel Hilbert spaces which is quite different to the one that we use here. 
Moreover, \cite{TorrecillaEtA_l2020} does not provide a general classification method. The procedures proposed in this paper vary from case to case because they select the optimal procedure for specific parametric situations. 

\item 
The paper on functional clustering by \cite{Delaigle_Hall_Pham_2019} is based on finding a finite-dimensional subspace in which the data are projected, and clustering is done by applying a modification of the $k$-means algorithm on those projections. A theoretical result related to perfect clustering is stated in Theorem 1 of this paper. In the homoscedastic case, \cite{Delaigle_Hall_Pham_2019} gives an explicit expression of the subspace in which the data should be projected (see Theorem 2 of this paper). 

The technique proposed in this paper has some advantage over our proposal in the sense that they can handle the homoscedastic (differences only in location) case. However, it suffers from several limitations, the main one being that {\tt \cite{Delaigle_Hall_Pham_2019} is able to deal with mixtures involving only two components}. Moreover, {\tt on the technical side, the theory of \cite{Delaigle_Hall_Pham_2019} has some limitations.} It requires to fix, arbitrarily,  $p \in \Nat$; then, the data are projected on a $p$-dimensional subspace in which the clustering is to be done. New issues appear in the way in which the subspace should be chosen, and the way in which the clusters can be constructed. According to Theorem 1 of this paper, the generators of the subspace must be chosen in a finite set with cardinality $a_n \to \infty$ as the sample size $n \to \infty$. Moreover, the partition of the data set must be chosen  between those in a finite set of Voronoi tessellations of $\Rea^p$ with cardinality $b_n \to \infty$ as $n \to \infty$. 
Additionally, the result needs some technical conditions like the existence of some $c \in (0,1)$ such that for every $C>0$ it happens that $a_n^pb_n \exp(-Cn^c) \to \infty$ as $n \to \infty$.

\end{itemize}

\section{Lemma to prove Theorem 3.2}
\label{Auxiliar.Lemma}

\begin{Lemm} \label{Lemm.aux.Theor3.2}
Under assumptions in Theorem 3.2, we have that 
$\Esp[\hat\Delta^2]  = O\left(N^{-1}\right).$
\end{Lemm}

\noindent
{\it  Proof.} 
We begin by analyzing the integral in 
\begin{equation} \label{Eq.ECM_covariance_0}
\Esp[\hat\Delta^2] = \int_0^1\int_0^1  \Esp (\hat \Sigma (t,s) - \Sigma (t,s))^2 dt ds
\end{equation}

To this, let us fix $s,t \in [0,1]$. The Lebesgue measure of the set $\{(s,s); s\in [0,1]\}$ is zero. Thus, we can assume that $s\neq t$. This allows to conclude that
\[
\Esp \left(\hat \Sigma (t,s) - \Sigma (t,s)\right)^2 = O\left(N^{-1}\right).
\]

Moreover, w.l.o.g., we can assume that the involved variables are centred. In order to simplify the notation, we will denote $x_i = z_i(t)$ and $y_i=z_i(s)$ for $i=1,\ldots,N$. We will also denote $\sigma_{\bf xy}= \mbox{Cov}(x_1,y_1)$, $\sigma_{\bf xx}= \mbox{Var}(x_1)$ and $\sigma_{\bf yy}= \mbox{Var}(y_1)$. We write $\ox =N^{-1}\sum_{i\leq N} x_i$ and  $\oy =N^{-1}\sum_{i\leq N} y_i$. Notice that with this notation $\Sigma(t,s)=\sigma_{\bf  xy}$.
We have that
\begin{equation} \label{Eq.ECM_covariance_1}
 \Esp \left(\hat \Sigma (t,s) - \Sigma (t,s)\right)^2 = 
  \Esp \left(\hat \Sigma (t,s)\right)^2 -2 \sigma_{\bf  xy} \Esp \left(\hat \Sigma (t,s)\right) +  \sigma_{\bf  xy}^2.
\end{equation}

Let us consider each term here separately. The independence between observations, gives that if $i\neq j$, then $\Esp[x_i y_j]=0$. On one hand, we have  
\begin{eqnarray} 
\nonumber
\Esp \left(\hat \Sigma (t,s)\right)  &=&
\frac 1 N \Esp \left(
\sum_{i=1}^N (x_i - \ox) (y_i - \oy)
\right)
\\
\nonumber
&=&
\frac 1 N 
 \Esp\left(\sum_i x_i y_i  
 - \frac 2 N \sum_{i,j} x_i y_j
 + \frac 1 {N^2} \sum_{i,j} x_iy_j
\right)
\\
&=&
\label{Eq.ECM_covariance_2}
\frac 1 N 
 \left(N\sigma_{\bf  xy}  
 -2 \sigma_{\bf  xy}
 + \frac 1 {N} \sigma_{\bf  xy}
\right) 
= \sigma_{\bf  xy} + O\left(N^{-1}\right).
\end{eqnarray}

On the other hand, we have
\begin{eqnarray} 
\nonumber
\lefteqn{\Esp \left(\hat \Sigma (t,s)\right)^2  =
\frac 1 {N^2} \Esp \left(
\sum_{i=1}^N (x_i - \ox) (y_i - \oy)
\right)^2}
\\
\nonumber
&=&
\frac 1 {N^2} 
 \Esp\left(
 \sum_i x_i y_i  
+ \frac {1-2N} {N^2} \sum_{h,k} x_hy_k
\right)^2
\\
\label{Eq.ECM_covariance_3}
&=&
\frac 1 {N^2} 
 \Esp\left(
\left( \sum_i x_i y_i  \right)^2
+2\frac {1-2N} {N^2} \sum_{i,h,k} x_i y_i  x_hy_k
+\frac {(1-2N)^2} {N^4} \left(\sum_{h,k} x_hy_k\right)^2
\right) .
\end{eqnarray}

Next, we compute the expectation of the three sums involved in (\ref{Eq.ECM_covariance_3}) as follows:
\begin{eqnarray}
\Esp
\left( \sum_i x_i y_i  \right)^2 \nonumber
&=&
\Esp
\left( \sum_i  x_i^2 y_i^2 
+ \sum_{i\neq j} x_i y_i x_j y_j 
\right)
\\
[4mm]
&=&
\label{Eq.ECM_covariance_4}
N\Esp
\left(  x_1^2 y_1^2 \right)
+ 
N(N-1)\sigma_{\bf  xy} ^2.
\end{eqnarray}

\begin{eqnarray}
\Esp\left(
\sum_{i,h,k} x_i y_i x_h y_k
\right)
\nonumber
&=&
 \Esp\left(
  \sum_{i} x_i^2 y_i^2 
  +
  \sum_{i} x_i y_i 
  \hspace{-3.5mm}
  \sum_{
  {\tiny \begin{array}{c}h,k
  \\
   h  \mbox{ \hspace{-.6mm}or  \hspace{-.6mm}} k \neq i
   \end{array}
   }} 
   \hspace{-3.5mm}x_hy_k
\right)
\\
\nonumber
&=&
N \Esp\left(  x_1^2 y_1^2 \right)
  +
  \Esp\left(
    \sum_{i} x_i y_i \sum_{h \neq i} x_hy_h
\right)
\\
\label{Eq.ECM_covariance_5}
&=&
N \Esp\left(  x_1^2 y_1^2 \right)
  +
  N(N-1)\sigma_{\bf  xy} ^2.
\\
[4mm]
\nonumber
 \Esp\left(
\left(\sum_{i,j} x_iy_j\right)^2
\right)
&=&
 \Esp\left(
\sum_{i,j} \left(x_iy_j\right)^2
+
\sum_{(i,j) \neq (h,k)} x_iy_jx_hy_k
\right)
\\
\label{Eq.ECM_covariance_6}
&=&
N \Esp\left(  x_1^2 y_1^2 \right)
  +  N(N-1)\sigma_{\bf xx}\sigma_{\bf yy}
  +   
\sum_{(i,j) \neq (h,k)} \Esp\left( x_iy_jx_hy_k\right).
\end{eqnarray}
Concerning last term in (\ref{Eq.ECM_covariance_6}), note that if $i=j$, the expectation is null unless $h=k$, and in this case, its value is $\sigma_{\bf xy}^2$. However, if $i \neq j$, the term is null again unless $i=k$ and $j=h$, in which case, the value is $\sigma_{\bf xy}^2$ as before. Consequently, we obtain that 
\begin{equation}
\label{Eq.ECM_covariance_7}
 \Esp\left(
\left(\sum_{i,j} x_iy_j\right)^2
\right)
=
N \Esp\left(  x_1^2 y_1^2 \right)
   +  N(N-1)\sigma_{\bf xx}\sigma_{\bf yy}
  +   
N^2\sigma_{\bf xy}^2.
\end{equation}

Replacing in (\ref{Eq.ECM_covariance_3}) the values we have obtained in (\ref{Eq.ECM_covariance_4}), (\ref{Eq.ECM_covariance_5}) and (\ref{Eq.ECM_covariance_7}), we obtain
\begin{eqnarray}
\nonumber
\Esp \left(\hat \Sigma (t,s)\right)^2
&=&
\frac 1 N \Big(
 \Esp \left(x_1^2 y_1^2  \right) + (N-1)\sigma_{\bf  xy} ^2
\\
\nonumber
&&
\hspace{8mm}+2\frac {1-2N} {N^2}\left(\Esp\left(  x_1^2 y_1^2 \right)  +  (N-1)\sigma_{\bf  xy} ^2\right)
\\
\nonumber
&&
\hspace{8mm}+
\frac {(1-2N)^2} {N^4} \left(\Esp\left(  x_1^2 y_1^2 \right)   +  (N-1)\sigma_{\bf xx}\sigma_{\bf yy}   +   N^2\sigma_{\bf xy}^2
\right)
\Big)
\\
\label{Eq.ECM_covariance_8}
&=&
\sigma_{\bf  xy} ^2 + O\left(N^{-1}\right).
\end{eqnarray}

If we substitute  this expression  and (\ref{Eq.ECM_covariance_2}) in (\ref{Eq.ECM_covariance_1}), we obtain that
\[
\Esp \left(\hat \Sigma (t,s) - \Sigma (t,s)\right)^2 = O\left(N^{-1}\right).
\]

From (\ref{Eq.ECM_covariance_2}) and (\ref{Eq.ECM_covariance_8}),
it follows that all terms involved in the right hand side of this expression are related to $\Esp\left(x_1^2 y_1^2 \right), \sigma_{\bf xy}, \sigma_{\bf xx}$ and $\sigma_{\bf yy}$, which in turn depend on the values of $t$ and $s$. However, under assumption A.1, they are uniformly bounded and the result follows because the integrals in (\ref{Eq.ECM_covariance_0}) are on the unit interval $[0,1]$. \FIN

\newpage
\section{Full Numerical Results}
\label{Subs.Full.Numerical}

\subsubsection{Simulation Results for Classification}
\label{Subs.Full.Classification}

Full results for {\tt all three scenarios} are given below:

\begin{table}[!ht]
	\begin{center}
		\caption{Misclassification rates for different GPs with difference only in locations (with standard error in brackets).} \label{Table.7.2}
		\vspace{0.1in}
		
		\tiny
		\setlength{\tabcolsep}{10pt}
		\begin{tabular}{|@{\hspace{.05mm}}c@{\hspace{.05mm}}|*{3}{c}|*{3}{c}|*{2}{c}|*{2}{c}|} \hline
			&  & PC &  &  & PLS &  &  &  &  & \\
			GP $\downarrow$ & CD-1NN & CD-SVM & CD-CENT & CD-1NN & CD-SVM & CD-CENT & DH-PC & DH-PLS & NP1 & NP2 \\ \hline 
			
			I & 0.0198 & 0.0255 & 0.0189 & 0.0007 & 0.0003 & 0.0128 & 0.2284 & 0.0007 & 0.0865 & 0.1502 \\[.1cm] 
			& (0.0008) & (0.0011) & (0.0008) & (0.0001) & (0.0001) &  (0.0043) & (0.0058) & (0.0001) & (0.0021) & (0.0017) \\[.1cm] 
			
			
			II & 0.0368 & 0.0346 & 0.0343 & 0.1369 & 0.1203 & 0.1271 & 0.0706 & 0.1195 & 0.1501 & 0.2306 \\[.1cm] 
			& (0.0065) & (0.0063) & (0.0063) & (0.0110) & (0.0098) & (0.0107) & (0.0082) & (0.0096) & (0.0046) & (0.0038) \\[.1cm] 

			III & 0.1736 & 0.1602 & 0.1652 & 0.0702 & 0.0589 & 0.1383 & 0.3596 & 0.0561 & 0.1941 & 0.2736 \\[.1cm] 
			& (0.0026) & (0.0026) & (0.0027) & (0.0015) & (0.0013) & (0.0089) & (0.0044) & (0.0012) & (0.0042) & (0.0010) \\[.1cm] 
			
			IV & 0.1213 & 0.1574 & 0.1157 & 0.0442 & 0.0351 & 0.1060 & 0.3528 & 0.0356 & 0.1846 & 0.2665 \\[.1cm] 
			& (0.0029) & (0.0043) & (0.0026) & (0.0012) & (0.0010) & (0.0091) & (0.0046) & (0.0009) & (0.0034) & (0.0013) \\[.1cm] 
			
			
			
			
			\hline
		\end{tabular}
	\end{center}
\end{table}

\begin{table}[!ht]
	\begin{center}
		\caption{Misclassification rates for different GPs with difference in location and scale (with standard error in brackets).} \label{Table.7.3}
		\vspace{0.1in}
		
		\tiny
		\setlength{\tabcolsep}{10pt}
		\begin{tabular}{|@{\hspace{.05mm}}c@{\hspace{.05mm}}|*{3}{c}|*{3}{c}|*{2}{c}|*{2}{c}|} \hline
			&  & PC &  &  & PLS &  &  &  &  & \\
			GP $\downarrow$ & CD-1NN & CD-SVM & CD-CENT & CD-1NN & CD-SVM & CD-CENT & DH-PC & DH-PLS & NP1 & NP2 \\ 
		\hline			
			I & 0.0054 & 0.0042 & 0.0080 & 0.0240 & 0.0285 & 0.1504 & 0.2889 & 0.0256 & 0.0469 & 0.1781 \\[.1cm] 
			& (0.0004) & (0.0003) & (0.0006) & (0.0012) & (0.0026) & (0.0092) & (0.0050) & (0.0009) & (0.0014) & (0.0017) \\[.1cm] 
			
			
			II & 0.0056 & 0.0046 & 0.0059 & 0.0996 & 0.0892 & 0.1127 & 0.1110 & 0.1635 & 0.0475 & 0.2170 \\[.1cm] 
			& (0.0007) & (0.0006) & (0.0007) & (0.0087) & (0.0084) & (0.0082) & (0.0104) & (0.0109) & (0.0014) & (0.0025) \\[.1cm] 
			
			III & 0.0105 & 0.0103 & 0.0117 & 0.0450 & 0.0414 & 0.0858 & 0.3954 & 0.1576 & 0.0495 & 0.2437 \\[.1cm] 
			& (0.0006) & (0.0006) & (0.0006) & (0.0021) & (0.0016) & (0.0035) & (0.0038) & (0.0021) & (0.0014) & (0.0013) \\[.1cm] 
			
			IV & 0.0106 & 0.0103 & 0.0120 & 0.0285 & 0.0232 & 0.0438 & 0.3882 & 0.1304 & 0.0510 & 0.2358 \\[.1cm] 
			& (0.0006) & (0.0005) & (0.0006) & (0.0017) & (0.0016) & (0.0025) & (0.0045) & (0.0020) & (0.0015) & (0.0013) \\[.1cm] 
			
			
			
			
			\hline
		\end{tabular}
	\end{center}
\end{table}

\begin{table}[!ht]
	\begin{center}
		\caption{Misclassification rates for different GPs with difference only in scales (with standard error in brackets).} \label{Table.2}
		\vspace{.1in}
		
		\tiny
		\setlength{\tabcolsep}{10pt}
		\begin{tabular}{|@{\hspace{.05mm}}c@{\hspace{.05mm}}|*{3}{c}|*{3}{c}|*{2}{c}|*{2}{c}|} \hline
			&  & PC &  &  & PLS &  &  &  &  &  \\
			GP $\downarrow$ & CD-1NN & CD-SVM & CD-CENT & CD-1NN & CD-SVM & CD-CENT & DH-PC & DH-PLS & NP1 & NP2 \\  \hline 
			
			I & 0.0130 & 0.0114 & 0.0124 & 0.0314 & 0.0244 & 0.0464 & 0.4912 & 0.4564 & 0.0494 & 0.2424 \\[.1cm] 
			& (0.0006) & (0.0005) & (0.0006) & (0.0018) & (0.0014) & (0.0023) & (0.0025) & (0.0024) & (0.0014) & (0.0013) \\[.1cm] 
			
			II & 0.0230 & 0.0185 & 0.0207 & 0.4153 & 0.3675 & 0.3600 & 0.4898 & 0.4884 & 0.0479 & 0.2452 \\[.1cm] 
			& (0.0008) & (0.0006) & (0.0007) & (0.0029) & (0.0024) & (0.0022) & (0.0028) & (0.0026) & (0.0015) & (0.0012) \\[.1cm] 
			
			III & 0.0126 & 0.0118 & 0.0122 & 0.0439 & 0.0421 & 0.0890 & 0.4867 & 0.4542 & 0.0465 & 0.2425 \\[.1cm] 
			& (0.0007) & (0.0006) & (0.0006) & (0.0021) & (0.0078) & (0.0034) & (0.0025) & (0.0024) & (0.0014) & (0.0013) \\[.1cm] 
			
			IV & 0.0137 & 0.0128 & 0.0125 & 0.0260 & 0.0216 & 0.0415 & 0.4839 & 0.4570 & 0.0518 & 0.2379 \\[.1cm] 
			& (0.0006) & (0.0006) & (0.0006) & (0.0015) & (0.0011) & (0.0019) & (0.0028) & (0.0024) & (0.0016) & (0.0012) \\[.1cm] 
			
			
			
			
			\hline
		\end{tabular}
	\end{center}
\end{table}

\noindent
R codes for our classification methods are available from this link: \url{https://www.dropbox.com/sh/dug1n4ufxubqplr/AADxA1myR3K-krvZsAEh-KYwa?dl=0}.

\newpage
\subsubsection{Simulation Results for Clustering}
\label{Subs.Full.Clustering}

Full results for {\tt two scenarios} are given below:

\begin{table}[!ht]
	\begin{center}
		\caption{One minus adjusted Rand indices for different GPs with difference in location and scale (with standard error in brackets).} \label{Table.7.5}
		\vspace{.1in}
		
		\tiny
		\setlength{\tabcolsep}{10pt}
		\begin{tabular}{|c|*{3}{c}|*{2}{c}|*{2}{c}|} \hline
			GP $\downarrow$ & CD-$k$-means & CD-Spectral & CD-mclust & CL1 & CL2 & DHP1 & DHP2 \\ \hline
			
			I & 0.1678 & 0.0082 & 0.0001 & 0.0239 & 0.0554 & 0.8386 & 0.0818 \\[.1cm] 
			& (0.0010) & (0.0002) & (0.0001) & (0.0007) & (0.0005) & (0.0064) & (0.0025) \\[.1cm] 
			
			II & 0.9858 & 0.9847 & 0.4240 & 0.5767 & 0.9967 & 0.5470 & 0.5149 \\[.1cm] 
			& (0.0003) & (0.0004) & (0.0030) & (0.0045) & (0.0018) & (0.0047) & (0.0049) \\[.1cm] 
			
			III & 0.4191 & 0.9857 & 0.0625 & 0.2891 & 0.9962 & 0.4137 & 0.5613 \\[.1cm] 
			& (0.0042) & (0.0099) & (0.0006) & (0.0000) & (0.0000) & (0.0054) & (0.0060) \\[.1cm] 
			
			IV & 0.0316 & 0.0000 & 0.0000 & 0.1833 & 0.6660 & 0.1379 & 0.5211 \\[.1cm] 
			& (0.0001) & (0.0000) & (0.0000) & (0.0000) & (0.0000) & (0.0033) & (0.0061) \\[.1cm] 
			
			
			\hline
			
		\end{tabular}
	\end{center}
\end{table}

\begin{table}[!ht]
\label{Table.RandIndex.Scales}
	\begin{center}
		\caption{One minus adjusted Rand indices for different GPs with difference only in scales (with standard error in brackets).} \label{Table.7.6}
		\vspace{.1in}
		
		\tiny
		\setlength{\tabcolsep}{10pt}
		\begin{tabular}{|c|*{3}{c}|*{2}{c}|*{2}{c}|} \hline
			GP $\downarrow$ & CD-$k$-means & CD-Spectral & CD-mclust & CL1 & CL2 & DHP1 & DHP2 \\ \hline
			
			I & 0.0063 & 0.0001 & 0.0000 & 0.8269 & 1.0017 & 0.9989 & 0.9966 \\[.1cm] 
			& (0.0002) & (0.0001) & (0.0000) & (0.0000) & (0.0000) & (0.0005) & (0.0008) \\[.1cm] 
			
			II & 0.0091 & 0.0100 & 0.0084 & 0.9065 & 1.0019 & 1.0001 & 0.9999 \\[.1cm] 
			& (0.0010) & (0.0019) & (0.0003) & (0.0007) & (0.0005) & (0.0002) & (0.0003) \\[.1cm] 
			
			III & 0.5549 & 0.9805 & 0.0856 & 0.9994 & 0.9998 & 0.9984 & 0.9967 \\[.1cm] 
			& (0.0052) & (0.0099) & (0.0006) & (0.0000) & (0.0000) & (0.0004) & (0.0007) \\[.1cm] 
			
			IV & 0.0102 & 0.0014 & 0.0005 & 0.8464 & 0.9928 & 0.9994 & 0.9980 \\[.1cm] 
			& (0.0021) & (0.0013) & (0.0004) & (0.0000) & (0.0000) & (0.0003) & (0.0006) \\[.1cm] 
			
			
			\hline
			
		\end{tabular}
	\end{center}
\end{table}

\noindent
Full result for the {\tt location scenario} (using the transformation $D^{\Sigma_d,4}_{d}$ stated in Section \ref{Clust.Location.Discussion}) is given below:

\begin{table}[!ht]
	\begin{center}
		\caption{One minus adjusted Rand indices for different GPs with difference in locations (with standard error in brackets).} \label{Table.3}
		\vspace{.1in}
		
		\tiny
		\setlength{\tabcolsep}{10pt}
		\begin{tabular}{|c|*{3}{c}|*{2}{c}|*{2}{c}|} \hline
			GP $\downarrow$ & CD-$k$-means & CD-Spectral & CD-mclust & C1 & C2 & DHP1 & DHP2 \\ \hline
			
			I & 0.0012 & 0.0814 & 0.0016 & 0.0001 & 0.0001 & 0.9896 & 0.0001 \\[.1cm] 
			& (0.0002) & (0.0272) & (0.0003) & (0.0000) & (0.0000) & (0.0001) & (0.0001) \\[.1cm] 
			
			
			III & 0.1649 & 0.3660 & 0.3007 & 0.0945 & 0.9975 & 0.1480 & 0.1623 \\[.1cm] 
			& (0.0017) & (0.0076) & (0.0024) & (0.0010) & (0.0043) & (0.0047) & (0.0075) \\[.1cm] 
			
			IV & 0.1257 & 0.3777 & 0.1784 & 0.1015 & 0.9001 & 0.0134 & 0.1473 \\[.1cm] 
			& (0.0019) & (0.0085) & (0.0021) & (0.0000) & (0.0000) & (0.0004) & (0.0039)\\[.1cm] 
			\hline
			
		\end{tabular}
	\end{center}
\end{table}

\noindent
R codes for our clustering methods are available from this link: \url{https://www.dropbox.com/sh/ont3ggvz44g5j07/AAC1PRuzIWx9_yFUiR5mk4Zna?dl=0}.

\end{document}